\documentclass{article}
\usepackage{amsthm}

\theoremstyle{definition}
\newtheorem{lemma}{Lemma}[section]
\newtheorem{proposition}{Proposition}[section]
\newtheorem{remark}{Remark}[section]
\newtheorem{example}{Example}[section]

\usepackage[utf8]{inputenc}
\usepackage[T2A]{fontenc}
\usepackage[english]{babel}
\usepackage[intlimits]{amsmath}
\usepackage{amssymb}
\allowdisplaybreaks[4]
\usepackage{graphicx}

\newtheorem{definition}{Definition}

\newtheorem{theorem}{Theorem}[section]

\newtheorem{corollary}{Corollary}[section]

\makeatletter
\@addtoreset{equation}{section}
\makeatother

\begin{document}

\title{Perturbation expansions of signal subspaces for long signals}
\author{Vladimir Nekrutkin \\
    Department of Statistical Modelling, \\
    Saint Petersburg State University \\
    198504, Universitetskiy pr. 28, Saint Petersburg, Russia\\
    \texttt{vnekr@statmod.ru}}

\maketitle





\begin{abstract}
Singular Spectrum Analysis and many other subspace-based methods of signal processing are implicitly
relying on the assumption of close proximity of unperturbed and perturbed
signal subspaces extracted by the Singular Value Decomposition
of special ``signal'' and ``perturbed signal'' matrices.
In this paper, the analysis of the largest principal angle between these subspaces is performed in terms
of the perturbation expansions of the corresponding orthogonal projectors. Applicable upper bounds are derived.
The main attention is paid to the asymptotic case when the length of the time series tends to infinity.
Results concerning conditions
for convergence, the rate of convergence, and the main terms of
proximity are presented.
\end{abstract}

%

\section{Introduction}
Though the first publication on the subspace-based methods of signal processing
 traditionally runs back to \cite{de_Prony795}, the explosive development of these methods
 ocurred during the last 30 -- 35 years.
In numerous scientific and engineering areas from climatology  \cite{Fraedrich86}
and meteorology \cite{WeareN82}
to acoustics \cite{Cadzow88},
from petroleum geology \cite[pp 41 -- 51]{Belonin_etal_71} to marine sciences \cite{Colebrook78},
from human \cite{BasilevskyH79} and animal \cite{EfimovG83} population
dynamics  to noise reduction problems \cite{KumaresanT80},
many researches independently proposed similar
``signal-subspace'' ideas and applied them to time series of their own interest.

Omitting details that are specific and important for each version of the method, the general scheme of the
sig\-nal-sub\-space approach can be explained as follows. Consider
a time series $\mathrm{F}_N=(x_0,\ldots,x_{N-1})$ treated as  a ``signal''.
(Note that $\mathrm{F}_N$ can be one-dimensional or multidimensional, real-valued or
complex-valued.)
This series is linearly transformed into a certain
$L\times K$ ``signal matrix'' ${\bf H}$. (The entries of ${\bf H}$ can also be either real-valued or complex-valued.)
It is supposed that the signal
$\mathrm{F}_N$ and the transformation $\mathrm{F}_N\mapsto {\bf H}$ are such that
$d\stackrel{\rm def}={\rm rank}\, {\bf H}<\min(L,K)$  and the linear space
spanned by columns of matrix ${\bf H}$
contains important
information about the series $\mathrm{F}_N$. The latter space is called {\it signal subspace}.

In terms of the ``covariance'' matrix
${\bf A}\stackrel{\rm def}={\bf H}{\bf H}^*$ this means that
the eigenspace $\mathbb{U}_0$ corresponding to zero eigenvalue
of the matrix ${\bf A}$ is not degenerate and that we are interested in the orthogonal complement
$\mathbb{U}_0^\bot$
of $\mathbb{U}_0$. Note that the matrix ${\bf H}$ must be a certain structured matrix.
Depending on the data at hand and on the problem under study, Hankel, Toeplitz, block-Hankel,
and other similar matrices are used in practice.

Assume now that  the perturbed series $\mathrm{F}_N(\delta)=\mathrm{F}_N+\delta\mathrm{E}_N$ is observed
rather than the unperturbed series $\mathrm{F}_N$,
where $\mathrm{E}_N=(e_0,\ldots,e_{N-1})$ is a ``noise series'' and $\delta$ stands for a formal perturbation parameter.
Thus, instead of  the ``signal matrix'' ${\bf H}$ we observe the perturbed
matrix ${\bf H}(\delta)={\bf H}+\delta {\bf E}$, where the ``noise matrix'' ${\bf E}$ is constructed from the series
$\mathrm{E}_N$ in the same manner as ${\bf H}$ is built from the series $\mathrm{F}_N$.
Note that generally there is no a priori assumption on the structure of the noise series. For example, the origin of $\mathrm{E}_N$
can be  random or deterministic.

Consider the Singular Value Decomposition (briefly, SVD) of ${\bf H}(\delta)$. If $\delta$ is small,
then continuity considerations show that the linear space ${\mathbb{U}}_0^\perp(\delta)$ spanned by $d$ leading
left singular vectors of this SVD can serve as an approximation of $\mathbb{U}_0^\perp$.
Thus we can expect that  useful (though approximate) information
about the signal $\mathrm{F}_N$ can be extracted from ${\mathbb{U}}_0^\perp(\delta)$ with the help
of a certain special technique. Since the linear space ${\mathbb{U}}_0^\perp(\delta)$
is spanned by $d$ leading eigenvectors of the matrix ${\bf A}(\delta)\stackrel{\rm def}={\bf H}(\delta){\bf H}^*(\delta)$,
this procedure is usually
formulated in terms of the eigendecomposition of the perturbed  ``covariance'' matrix  ${\bf A}(\delta)$.

These or similar ideas under different names (for example, Eigenvector Filtering \cite{Colebrook78}, Extended Empirical
Orthogonal Functions \cite{WeareN82}, Karhunen-Lo\`{e}ve Time Series Analysis \cite{BasilevskyH79},
Singular System Analysis \cite{BroomheadK86a}, Dynamic Principal Component Analysis \cite{KuSG95})
or without any special names are used in many
publications. The so-called Singular Spectrum Analysis
(whether treated as in \cite{VautardYG92} or
as in \cite{GNZh01}) has the same origin.

The most well-known example of this scheme is the following. Suppose that a one-dimensional series $\mathrm{F}_N$
is governed by some minimal Linear Recurrent Formula (briefly, LRF)
\begin{gather*}
x_n=\sum_{k=1}^{d}b_kx_{n-k}, \quad n\geq d
\end{gather*}
of order $d$. Transform $\mathrm{F}_N$ into the
 $L\times K$ Hankel ``trajectory'' matrix ${\bf H}$ with rows $(x_j,\ldots,x_{K+j-1})$, where $j=0,\ldots,L-1$ and $L+K=N+1$.
 If $\min(L,K)>d$, then ${\rm rank}\,{\bf H}=d=\dim \mathbb{U}_0^\perp$ and $\dim \mathbb{U}_0=K-d>0$.

The knowledge of $\mathbb{U}_0^\perp$ (or, equivalently, of $\mathbb{U}_0$) provides
 an essential information about the series $\mathrm{F}_N$.
In particular, if we consider Schubert basis $\mathcal{Y}=(Y_1,\ldots,Y_{L-d})$ of the linear space $\mathbb{U}_0$
(see \cite{Buch94} or \cite[ch. 5]{GNZh01} for relations between
this basis and LRFs), then the first element $Y_1$ of $\cal Y$  determines the minimal LRF governing the
series $\mathrm{F}_N$. Finding the roots of the related characteristic polynomial is the usual goal in
practical investigations.

There is a large number of publications dedicated to methods and algorithms for estimation of these roots by the use of the
linear space ${\mathbb{U}}_0^\perp(\delta)$ in the case when $\mathrm{E}$ is a white-noise series
(see, for example, chapters 11 and 13 in \cite{Marple87}).
Methods named as MUSIC and ESPRIT are presently the  most popular.
Standard references for the origins of these approaches are  \cite{Schmidt79} and \cite{RoyKailath89}.
A modern description of ESPRIT  together with the first-order
perturbation analysis of its
performance can be found in \cite{BadeauRG08}.

Let us discuss the aim, the technique, and the results  of the paper.
The main goal is to estimate the difference between
${\mathbb{U}}_0^\perp$ and ${\mathbb{U}}_0^\perp(\delta)$ in the case of long signals.

We describe both unperturbed ${\mathbb{U}}_0^\perp$ and perturbed ${\mathbb{U}}_0^\perp(\delta)$ signal subspaces by
their orthogonal projectors ${\bf P}_0^\perp$ and ${\bf P}_0^\perp(\delta)$, respectively.
The spectral norm (also known as the operator norm) $\|\,\cdot\,\|$ is used to measure the difference between these projectors.
(Spectral norm of a matrix
can be defined as its maximal singular value.) Since the dimension of ${\mathbb{U}}_0^\perp$ and ${\mathbb{U}}_0^\perp(\delta)$
is the same,  the norm
$\big\|{\bf P}_0^\perp(\delta)-{\bf P}_0^\perp\big\|$ is the sine of the largest principal
angle between the linear spaces ${\mathbb{U}}_0^\perp(\delta)$ and ${\mathbb{U}}_0^\perp$. (For example, this follows
from \cite[ch. 3 \S 39]{Akhiez93}.)

To estimate this norm
from above we use the perturbation expansions of ${\bf P}_0^\perp(\delta)$. These expansions
(see formulas \eqref{eq:pert_0_B} and \eqref{eq:pert_0} of Section \ref{sect:Zero_eigen})
are valid under the condition
$|\delta|<\delta_0$ for some positive $\delta_0$ and
can be straightforwardly deduced from the  classical monograph \cite[ch. II \S 3]{Kato66}. It follows from
\eqref{eq:pert_0} that the projector ${\bf P}_0^\perp(\delta)$ is continuous in $\delta$ as $|\delta|<\delta_0$. In other
words, if $|\delta|<\delta_0$ then ${\mathbb{U}}_0^\perp(\delta)$ can be continuously traced back to
${\mathbb{U}}_0^\perp$.

Suitable upper bounds for $\big\|{\bf P}_0^\perp(\delta)-{\bf P}_0^\perp\big\|$ lead to the solution of several
problems related to the proximity of perturbed and unperturbed signal subspaces in the case
of long signals. Solving these problems, we restrict ourselves to real-valued signals governed by LRFs.

More precisely, we consider an infinite signal series $\mathrm{F}=(x_0,x_1,\ldots,x_n,\ldots)$ which is governed
by a minimal linear recurrent formula of order $d$.
As it was already described, finite segments $\mathrm{F}_N=(x_0,x_1,\ldots,x_{N-1})$ of the series
$\mathrm{F}$ are transformed into signal $L\times K$ Hankel matrices ${\bf H}$, where $L$ depends on $N$.
Thus the signal subspace  $\mathbb{U}_0^\perp$ with the related projector ${\bf P}_0^\perp$ depends on $N$ too.
The noise series is supposed to be infinite as well and the expansions \eqref{eq:pert_0_B} and \eqref{eq:pert_0} are valid
   for all $\delta$
 such that  $|\delta|<\delta_0=\delta_0(N)$.

If $N\rightarrow \infty$, then we come to several problems related to the asymptotic behavior of
$\big\|{\bf P}_0^\perp(\delta)-{\bf P}_0^\perp\big\|$.
One problem is the behavior of $\delta_0(N)$ as $N\rightarrow \infty$ for a certain signal and noise series.
If $\delta_0(N)\rightarrow \infty$, then  we can use expansions of ${\bf P}_0^\perp(\delta)$ for any fixed $\delta$
provided that $N$ is sufficiently large. If $\delta_0(N)$ is bounded from above (and is separated
from zero), then expansions are valid at a certain interval around zero.

Another problem is related to the proximity of ${\bf P}_0^\perp(\delta)$ and ${\bf P}_0^\perp$.
In the present paper, we indicate three
origins of the proximity ${\bf P}_0^\perp(\delta)\approx{\bf P}_0^\perp$.
 The first origin can be formulated only in terms of the extremal singular values of the matrices ${\bf H}$
and ${\bf E}$. Roughly speaking,
the difference $\|{\bf P}_0^{\perp}(\delta)-{\bf P}_0^{\perp}\|$ is small if
the ``signal-to-noise ratio'' $\|{\bf H}\|/\|{\bf E}\|$ is large and the
``pseudo-condition number'' of the matrix  ${\bf H}{\bf H}^{\rm T}$, which is defined as the ratio of the maximal
and minimal positive eigenvalues of this matrix, is not very large. If these two conditions hold  asymptotically, then
$\|{\bf P}_0^{\perp}(\delta)-{\bf P}_0^{\perp}\|\rightarrow 0$.

Nevertheless, the convergence $\|{\bf P}_0^{\perp}(\delta)-{\bf P}_0^{\perp}\|\rightarrow 0$ can take place even if the
``signal-to-noise ratio'' is bounded from above. This happens when both (column and row)
pairs of linear spaces produced by matrices ${\bf H}$ and ${\bf E}$ are asymptotically orthogonal.

The third origin is connected with the
condition ${\bf E}{\bf E}^{\rm T}/K\rightarrow {\bf I}$, where ${\bf I}$ stands for the identity
matrix. This case corresponds to the white-noise ``noise series''. Note that if the ``noise series'' is random, then the perturbed operator is also random.
Moreover, for fixed $\delta$ and $N$ the inequality $|\delta|<\delta_0(N)$ holds only with a certain probability.

One more interesting problem related to the convergence $\|{\bf P}_0^{\perp}(\delta)-{\bf P}_0^{\perp}\|\rightarrow 0$ is the form of the
main term of the approximation  ${\bf P}_0^\perp(\delta)\approx {\bf P}_0^\perp$. Since $\delta$ is fixed and $N\rightarrow \infty$,
this term does not necessarily have to be linear in $\delta$.

Let us explain the style of our results for these asymptotic proximity
problems. For this purpose we use   simplified but typical examples.
More general examples are considered in Section \ref{ssect:hankel_coarse}.

\begin{example}
\label{ex:intro}
Let the signal series be defined as $x_n=\cos(2\pi\omega n+\phi)$ with $\omega \in (0,1/2)$.
Then the pse\-udo-con\-di\-tion  number of the matrix ${\bf H}{\bf H}^{\rm T}$ is bounded from above as
$N\rightarrow \infty$.\\
a) Suppose that the ``noise series'' is a linear stationary random series defined on a probability space
$(\Omega, {\cal F}, \mathbb{P})$. (Exact definitions and restrictions can be
found in Lemma \ref{summ:series_rand} of Section \ref{ssect:hankel_coarse}.)
If $L/N\rightarrow\alpha\in(0,1)$ as $N\rightarrow \infty$, then for a certain $\Omega'\in {\cal F}$ with
$\mathbb{P}(\Omega')=1$, for any $\omega\in \Omega'$, and for any $\delta \in \mathbb{R}$
\begin{gather*}
\limsup_{N}\sqrt{N (\ln N)^{-1}}\ \big\|{\bf P}_0^{\perp}(\delta)-{\bf P}_0^{\perp}\big\|
\leq
b_1\,|\delta| \quad \mbox{and} \\
\limsup_{N}{N (\ln N)^{-1}}\,
\big\|{\bf P}_0^{\perp}(\delta)-{\bf P}_0^{\perp}-\delta {\bf V}_{0}^{(1)}\big\|
\leq
b_2\,\delta^2,
\end{gather*}
where  $\delta {\bf V}_0^{(1)}$ stands for
the linear in $\delta$ term of the difference ${\bf P}_0^{\perp}(\delta)-{\bf P}_0^{\perp}$.
It should be mentioned that the ``signal-to-noise ratio'' $\|{\bf H}\|/\|{\bf E}\|$ tends almost surely to infinity as $N\rightarrow \infty$
in this example.\\
b) Assume that the ``noise series'' is also oscillating and defined by $e_n=\cos(2\pi\omega'n+\phi')$ with $(0,1/2)\ni\omega'\neq \omega$.
If $\min(L,K)\rightarrow \infty$, then for any $\delta\in (-1/4,1/4)$ and for $N>N_0(\delta)$
\begin{gather*}
\begin{split}
&\big\|{\bf P}_0^{\perp}(\delta)-{\bf P}_0^{\perp}\big\|=|\delta|\,O\big(1/\min(L,K)\big) \ \ \mbox{and}\\
&\big\|{\bf P}_0^{\perp}(\delta)-{\bf P}_0^{\perp}-{\bf L}(\delta)\big\|=\delta^2\,O\big(1/\min(L^2,K^2)\big),
\end{split}
\end{gather*}
where ${\bf L}(\delta)$ is the non-linear in $\delta$ operator defined in Theorem \ref{theor:precise_proj}.
This case corresponds to the asymptotic biorthogonality
of matrices ${\bf H}$ and ${\bf E}$.\\
c) If white noise stands for the ``noise series'' and $L=L_0={\rm const}>d$, then for a certain $\delta_0>0$,  for some
 $\Omega'\in {\cal F}$ with $\mathbb{P}(\Omega')=1$, for any $|\delta|<\delta_0$, and for any $\omega\in \Omega'$
\begin{gather*}
\limsup_{N}\sqrt{N (\ln \ln N)^{-1}}\ \big\|{\bf P}_0^{\perp}(\delta)-{\bf P}_0^{\perp}\big\|\leq c_1\,|\delta|
\quad \mbox{and} \\
\limsup_{N}{N (\ln \ln N)^{-1}}\,\big\|{\bf P}_0^{\perp}(\delta)-{\bf P}_0^{\perp}-{\bf T}(\delta)\big\|
\leq c_2\delta^2,
\end{gather*}
where  the non-linear in $\delta$ operator ${\bf T}(\delta)$ is defined in Theorem \ref{theor:precise_proj_noise}.
The proof of this case is based on the
almost sure (briefly, a.s.) convergence ${\bf E}{\bf E}^{\rm T}/K\rightarrow {\bf I}$.
\end{example}

The present paper is organized as follows.
Section \ref{sect:Zero_eigen} contains basic theoretical results. General definitions and notation, as well as formulas for the perturbation
expansion of the perturbed projector ${\bf P}_0^{\perp}(\delta)$
are placed in Section \ref{ssect:zero_eigen_series}.

In Section \ref{ssect:appr_err} general upper bounds for the difference
$\big\|{\bf P}_0^{\perp}(\delta)-{\bf P}_0^{\perp}\big\|$ are constructed and discussed.
In Theorem \ref{theor:gen_upper} we derive the inequality \eqref{eq:gen_upper} such that
${\bf P}_0^{\perp}(\delta)={\bf P}_0^{\perp}$ iff the right-hand side of \eqref{eq:gen_upper} equals zero.
(The conditions for the equality of the perturbed and unperturbed projectors are discussed in Section \ref{ssect:zero_pert}.)
Corollaries \ref{cor:coarse_proj} and \ref{cor:coarse_orth} contain interpretable sufficient conditions for the proximity
of ${\bf P}_0^{\perp}(\delta)$ and ${\bf P}_0^{\perp}$.

The rest  of Section \ref{ssect:appr_err} is devoted to the main terms of the approximation
${\bf P}_0^{\perp}(\delta)\approx{\bf P}_0^{\perp}$.
These results are collected in Theorems \ref{theor:coarse_proj} - \ref{theor:precise_proj_noise} of
 Section \ref{sssect:main_terms}.

In Section \ref{sect:Hankel_pert}  we perform the asymptotic analysis of
$\|{\bf P}_0^{\perp}(\delta)-{\bf P}_0^{\perp}\|$ in the case when the size of $L\times K$ matrices ${\bf H}$ and
${\bf E}$ increases. More precisely, we  consider  conditions for
the convergence $\big\|{\bf P}_0^{\perp}(\delta)-{\bf P}_0^{\perp}\big\|\rightarrow 0$ as $\max(L,K)\rightarrow \infty$  and study  the rate
 of this convergence.

Two general assertions of this kind are presented in Section \ref{ssect:asympt}.
Starting from Section \ref{ssect:hankel_coarse}, the case of Hankel matrices ${\bf H}$ and ${\bf E}$ is examined more
thoroughly. Several more or less general time series of finite rank are taken as signal and noise series and a number
of inequalities similar to those of Example \ref{ex:intro} are derived. To obtain these inequalities we need
to study the asymptotic behavior for the
minimal and maximal positive singular values of the corresponding $L\times K$ Hankel matrices as $L+K\rightarrow \infty$.
Results for deterministic  series
 are collected in Lemma \ref{summ:series}.

The case of stationary random  ``noise series'' is of  special interest.
To estimate  norms of the related Hankel matrices, see Lemma \ref{summ:series_rand}, we use recent results published in
\cite{Meckes07} and \cite{Adamczak08}. Proposition \ref{prop:LP} of Appendix 2 contains
proofs of inequalities we apply in our considerations. In its turn, inequalities concerning white-noise
``noise series'' are based on results of Section \ref{ssect:WN_gen} which is also placed in Appendix 2.

In Section \ref{sect:hankel_precise}
we present several examples of ``signal'' and ``noise'' series such that both the exact rates of convergence
$\big\|{\bf P}_0^{\perp}(\delta)-{\bf P}_0^{\perp}\big\|\rightarrow 0$ and the related
 main terms are expressed in the explicit form. These  examples are rather simple, but they show that refined
 versions  of the inequalities of Section \ref{ssect:asympt} can provide precise results.

In general, the term ``signal-subspace method'' implies that a method provides a small error if the perturbed
 subspace $\mathbb{U}_0^\perp(\delta)$ is close to the unperturbed $\mathbb{U}_0^\perp$. Propositions
 related to the estimation
 of LRFs and to the so-called Least-Square ESPRIT confirm this statement.
 These propositions are placed in Section \ref{sect:appl}.

Singular Spectrum Analysis (briefly, SSA), which can also be considered as a signal-subspace method, shows a different effect.
Namely, the SSA reconstruction  may not converge to a signal  even
if $\big\|{\bf P}_0^{\perp}(\delta)-{\bf P}_0^{\perp}\big\|\rightarrow 0$ as $N\rightarrow \infty$.
The related
discussion can be found in Section \ref{ssect:Hankel}.

Proofs of assertions are placed in Appendix 1.
It its turn, Appendix 2 collects auxiliary results related to random ``noise series''. Note that some simple but
laborsome asymptotic considerations are used without proofs.

\section{Perturbations of signal subspaces}
\label{sect:Zero_eigen}
In what follows we do not distinguish linear operators from their
matrices. In particular, we use the sign of transposition ${\bf
B}^{\rm T}$ instead of the conjugation sign ${\bf B}^{*}$. To avoid misunderstanding, we sometimes
use the notation $0_M$ for the zero vector in $\mathbb{R}^M$ as well as
the notation ${\bf 0}_{MP}$
for the zero operator $\mathbb{R}^M\mapsto \mathbb{R}^P$.

\subsection{Perturbation series}
\label{ssect:zero_eigen_series}
Let ${\bf H}$ be a real-valued non-zero $\mathbb{R}^K\mapsto \mathbb{R}^L$ matrix. Then
${\bf A}\stackrel{\rm def}={\bf H H}^{\rm T}$
is a self-adjoint semi-definite operator ${\bf A}: \mathbb{R}^L\mapsto \mathbb{R}^L$.
We assume that $d\stackrel{\rm def}={\rm rank}\, {\bf H}<\min(L,K)$.
Denote $\Sigma$ the set of eigenvalues of ${\bf A}$. Then $0\in\Sigma\subset [0,\infty)$.

Consider one more real-valued non-zero matrix
${\bf E}: \mathbb{R}^K\mapsto \mathbb{R}^L$ and put ${\bf H}(\delta)={\bf H}+\delta {\bf E}$. Then
\begin{gather*}
{\bf A}(\delta)\stackrel{\rm def}={\bf H}(\delta){\bf H}(\delta)^{\rm T}=
{\bf A}+\delta{\bf A}^{(1)}+\delta^2{\bf A}^{(2)}={\bf A}+{\bf B}(\delta),
\end{gather*}
where ${\bf A}^{(1)}={\bf H} {\bf E}^{\rm T}+{\bf E}{\bf H}^{\rm T}$,
${\bf A}^{(2)}={\bf E} {\bf E}^{\rm T}$, and ${\bf B}(\delta)=\delta{\bf A}^{(1)}+\delta^2{\bf A}^{(2)}$.
 Note that both ${\bf A}^{(1)}$ and ${\bf A}^{(2)}$ are self-adjoint operators,
and ${\bf A}(\delta)$ is positive semi-definite for any $\delta\in \mathbb{R}$.

We interpret ${\bf H}$ as a ``signal matrix'' and ${\bf E}$ as a ``noise matrix''.
It is supposed that we observe the perturbed matrix ${\bf H}(\delta)$ and are interested in the ``signal subspace''
spanned by the columns of the matrix ${\bf H}$.

The signal subspace is  the orthogonal complement $\mathbb{U}_0^{\perp}$ to the eigenspace
$\mathbb{U}_0$, corresponding to the zero eigenvalue of the matrix ${\bf H}{\bf H}^{\rm T}$. If ${\bf P}_0$
stands for the orthogonal projector onto $\mathbb{U}_0$, then ${\bf P}_0^\perp= {\bf I}-{\bf P}_0$
is the orthogonal projector onto $\mathbb{U}_0^{\perp}$. (Here ${\bf I}$ is the identical operator
$\mathbb{R}^L\mapsto \mathbb{R}^L$.)

Our goal is to compare the perturbed projector ${\bf P}_0^{\perp}(\delta)$  with the unperturbed projector
${\bf P}_0^{\perp}$.
In what follows, we set
$\mu_{\min}=\min\{\mu\in \Sigma\ {\rm such\ that}\ \mu>0 \}.$

The following assertion can be easily deduced from general results given in  \cite[ch. II \S 3]{Kato66}.
Denote ${\bf S}_0$ the pseudoinverse to ${\bf H}{\bf H}^{\rm T}$. Set ${\bf S}_0^{(0)}=-{\bf P}_0$ and
${\bf S}_0^{(k)}={\bf S}_0^{k}$ for $k\geq 1$.
\begin{theorem}
\label{theor:full_decomp}
1. Let $\delta_0>0$ and assume that
\begin{gather}
\label{eq:cond_main}
\|{\bf B}(\delta)\|<\mu_{\min}/2
\end{gather}
for all $\delta\in (-\delta_0,\delta_0)$.
Then  the perturbed projector
${\bf P}^{\perp}_0(\delta)$   admits the expansion
\begin{gather}
\label{eq:pert_0_B}
{\bf P}^{\perp}_0(\delta)={\bf P}_0^{\perp}+\sum_{p=1}^\infty{\bf W}_p(\delta),
\end{gather}
where
\begin{gather}
\label{eq:Wp}
{\bf W}_p(\delta)
=
(-1)^{p} \sum_{l_1+\ldots+l_{p+1}=p,\  l_j\geq 0}
{\bf W}_p(l_1,\ldots,l_{p+1})
\end{gather}
and
\begin{gather*}
{\bf W}_p(l_1,\ldots,l_{p+1})=
{{\bf S}}_{0}^{(l_1)}{{\bf B}}(\delta){{\bf S}}_{0}^{(l_2)}\ldots
\,{\bf S}_{0}^{(l_{p})}{{\bf B}}(\delta){\bf S}_{0}^{(l_{p+1})}.
\end{gather*}
In addition,
\begin{gather}
\label{eq:pert_0}
{\bf P}^{\perp}_0(\delta)=
{\bf P}_{0}^{\perp}+\sum_{n=1}^\infty\delta^n {\bf V}_{0}^{(n)},
\end{gather}
where
\begin{gather}
\label{eq:P_0_delta_n}
{\bf V}_{0}^{(n)}=
\sum_{p=\lceil n/2 \rceil}^n (-1)^{p}
\sum_{\substack{s_1+\ldots+s_p=n,\ s_i=1,2
\\
l_1+\ldots+l_{p+1}=p,\  l_j\geq 0}}{\bf V}_{0}^{(n)}({\bf s},{\ell}),
\end{gather}
${\bf s}=(s_1,\ldots,s_p)$, ${\ell}=(l_1,\ldots,l_{p+1})$, and
\begin{gather*}
{\bf V}_{0}^{(n)}({\bf s},{\ell})=
{{\bf S}}_{0}^{(l_1)}{{\bf A}}^{(s_1)}{{\bf S}}_{0}^{(l_2)}\ldots\,{{\bf A}}^{(s_p)}{{\bf S}}_{0}^{(l_{p+1})}.
\end{gather*}
\end{theorem}

\begin{remark}
\label{rem:full_decomp}
1. Both series \eqref{eq:pert_0_B} and \eqref{eq:pert_0} converge in the spectral norm.\\
2. Denote
\begin{gather}
\label{eq:B_delta}
B(\delta)=|\delta|\, \|{\bf A}^{(1)}\|+\delta^2 \|{\bf A}^{(2)}\|.
\end{gather}
If $\delta_0>0$ and $B(\delta_0)=\mu_{\min}/2$, then the
inequality \eqref{eq:cond_main} is valid for all $\delta$ such that
$|\delta|<\delta_0$.\\
3. Since ${\bf H}{\bf H}^{\rm T}=\sum_{\mu>0}\mu {\bf P}_\mu$, then ${\bf S}_0=\sum_{\mu>0}{\bf P}_\mu/\mu$.
(Here and further we write $\sum_{\mu>0}$ instead of
$\sum_{\Sigma\ni\mu>0}$.) It is easy to show that $\big\|{\bf S}_0^{(k)}\big\|=1/\mu_{\min}^k$ for any $k\geq 0$.\\
4.
The coefficient ${\bf V}^{(1)}_0$ of the linear in $\delta$ term  of the right-hand side of \eqref{eq:pert_0} has the
form
\begin{gather}
\begin{split}
&{\bf V}^{(1)}_0={\bf P}_{0}{\bf A}^{(1)}{\bf S}_0+{\bf S}_{0}{\bf A}^{(1)}{\bf P}_0
\\
&=
{\bf P}_{0}{\bf E}{\bf H}^{\rm T}{\bf S}_0+{\bf S}_0{\bf H}{\bf E}^{\rm T}{\bf P}_{0}.
\end{split}
\label{eq:lin_proj}
\end{gather}
5. Accurate calculations show that
\begin{gather}
\begin{split}
&{\bf V}^{(2)}_0=
{\bf P}_0{\bf A}^{(2)}{\bf S}_0+{\bf S}_0{\bf A}^{(2)}{\bf P}_0
\\
&+{\bf P}_0{\bf A}^{(1)}{\bf P}_0{\bf A}^{(1)}{\bf S}_0^{2}
+{\bf P}_0{\bf A}^{(1)}{\bf S}_0^{2}{\bf A}^{(1)}{\bf P}_0
\\
&+{\bf S}_0^{2}{\bf A}^{(1)}{\bf P}_0{\bf A}^{(1)}{\bf P}_0
-{\bf P}_0{\bf A}^{(1)}{\bf S}_0{\bf A}^{(1)}{\bf S}_0
\\
&-{\bf S}_0{\bf A}^{(1)}{\bf P}_0{\bf A}^{(1)}{\bf S}_0
-{\bf S}_0{\bf A}^{(1)}{\bf S}_0{\bf A}^{(1)}{\bf P}_0.
\end{split}
\label{eq:sec_order}
\end{gather}
6. The term ${\bf W}_1(\delta)$ in the right-hand side of \eqref{eq:pert_0_B} can be expressed as
\begin{gather}
\begin{split}
&{\bf W}_1(\delta)={\bf P}_{0}{\bf B}(\delta){\bf S}_0+{\bf S}_0{\bf B}(\delta){\bf P}_{0}
\\
&=
\delta\,{\bf V}^{(1)}_0+
\delta^2\big({\bf P}_{0}{\bf A}^{(2)}{\bf S}_0+{\bf S}_0{\bf A}^{(2)}{\bf P}_{0}\big).
\end{split}
\label{eq:first_term}
\end{gather}
\end{remark}

\subsection{Approximation errors}
\label{ssect:appr_err}
It this section we derive
upper bounds
for the norm $\big\|{\bf P}_0^{\perp}(\delta)-{\bf P}_0^{\perp}\big\|$.
These bounds produce interpretable sufficient conditions for the proximity of linear spaces
$\mathbb{U}_0^\perp(\delta)$
and  $\mathbb{U}_0^\perp$. We also present operators
that can play a role of the main terms of the difference ${\bf P}_0^{\perp}(\delta)-{\bf P}_0^{\perp}$.

\subsubsection{Zero perturbation effects}
\label{ssect:zero_pert}
We start with the necessary and sufficient conditions for the equality ${\bf P}_0^{\perp}(\delta)={\bf P}_0^{\perp}$.
Consider the function $B(\delta)$ defined in \eqref{eq:B_delta}.

\begin{theorem}
\label{theor:zero_pert}
Let $\delta_0>0$ and assume  that $B(\delta_0)=\mu_{\min}/2$. Then the following assertions are equivalent:\\
1. ${\bf P}_0^{\perp}(\delta)={\bf P}_0^{\perp}$ for all $\delta\in (-\delta_0, \delta_0)$;\\
2.
\begin{gather}
\label{eq:zero_pert}
{\bf S}_0{\bf H}{\bf E}^{\rm T}{\bf P}_0={\bf S}_0{\bf E}{\bf E}^{\rm T}{\bf P}_0={\bf 0}_{LL};
\end{gather}
3.
\begin{gather}
\label{eq:zero_pert1}
{\bf S}_0{\bf B}(\delta){\bf P}_0={\bf 0}_{LL}
\end{gather}
for all $\delta$ from a certain neighbourhood of zero;\\
4. \begin{gather}
\begin{split}
&{\bf S}_0{\bf H}{\bf E}^{\rm T}{\bf P}_0 + {\bf P}_0{\bf H}{\bf E}^{\rm T}{\bf S}_0
\\
&=
{\bf S}_0{\bf E}{\bf E}^{\rm T}{\bf P}_0 + {\bf P}_0{\bf E}{\bf E}^{\rm T}{\bf S}_0={\bf 0}_{LL};
\end{split}
\label{eq:sum_SHEP}
\end{gather}
5.
\begin{gather}
\label{eq:zero_pert2}
{\bf H}{\bf E}^{\rm T}{\bf P}_0={\bf 0}_{LL}\ \ {\rm and} \ \ {\bf H}^{\rm T}{\bf E}{\bf E}^{\rm T}{\bf P}_0={\bf 0}_{LK}.
\end{gather}
\end{theorem}

\begin{remark}
1. The equality ${\bf S}_0{\bf E}{\bf E}^{\rm T}{\bf P}_0={\bf 0}_{LL}$ is equivalent to
${\bf P}_0^\perp{\bf E}{\bf E}^{\rm T}{\bf P}_0={\bf 0}_{LL}$.\\
2.  Let us discuss the conditions \eqref{eq:zero_pert2}.

Let $\mathbb{U}_E$ stand for the linear space spanned by the
columns of the matrix ${\bf E}$ and denote $s=\dim \mathbb{U}_E$.
Suppose that there exists an orthonormal basis $P_1,\ldots,P_s$ of the space $\mathbb{U}_E$ such that a)
each $P_i$ is an eigenvector
of the matrix ${\bf E}{\bf E}^{\rm T}$, b) $P_1,\ldots,P_l\in \mathbb{U}_0^\perp$, and
c) $P_{l+1},\ldots,P_s\in \mathbb{U}_0$
for some $0\leq l\leq s$. Then  ${\bf H}^{\rm T}{\bf E}{\bf E}^{\rm T}{\bf P}_0={\bf 0}_{LK}$.
(Note that the latter equality is equivalent to ${\bf P}_0^\perp{\bf E}{\bf E}^{\rm T}{\bf P}_0={\bf 0}_{LL}$.)

Indeed, consider Singular Value Decomposition ${\bf E}{\bf E}^{\rm T}=\sum_{i=1}^s\nu_iP_iP_i^{\rm T}$
 of the  matrix ${\bf E}{\bf E}^{\rm T}$. Then
${\bf P}_0 P_i\!=\!0_{L}$ for $1\!\leq\! i\!\leq \!l$, ${\bf P}_0^\perp P_i=0_{L}$ for $i>l$, and
\begin{gather*}
\begin{split}
&{\bf P}_0^\perp{\bf E}{\bf E}^{\rm T}{\bf P}_0=\sum_{i=1}^l\nu_i{\bf P}_0^\perp P_i\left(P_i^{\rm T}{\bf P}_0\right)
\\
&+
\sum_{i=l+1}^s\nu_i\left({\bf P}_0^\perp P_i\right)P_i^{\rm T}{\bf P}_0={\bf 0}_{LL}.
\end{split}
\end{gather*}

There are several important particular cases related to the situation under discussion.

If $l=s$, then $\mathbb{U}_E\subset \mathbb{U}_0^\perp$, ${\bf E}^{\rm T}{\bf P}_0={\bf 0}_{LK}$, and both conditions of
\eqref{eq:zero_pert2}
are fulfilled.
The example when ${\bf E}$ is proportional to ${\bf H}$ (then the equality
${\bf P}_0^{\perp}(\delta)={\bf P}_0^{\perp}$ becomes evident)
is just a particular case of the inclusion $\mathbb{U}_E\subset \mathbb{U}_0^\perp$.

If $l<s$, then the natural sufficient condition for the equality ${\bf H}{\bf E}^{\rm T}{\bf P}_0={\bf 0}_{LL}$ is
${\bf H}{\bf E}^{\rm T}={\bf 0}_{LL}$.
The analogous sufficient condition ${\bf H}^{\rm T}{\bf E}={\bf 0}_{KK}$ for the equality
${\bf H}^{\rm T}{\bf E}{\bf E}^{\rm T}{\bf P}_0={\bf 0}_{LK}$
corresponds to the case $l=0$ with $\mathbb{U}_E\subset \mathbb{U}_0$.

Lastly, suppose that $s=L$ and that all singular values of the matrix ${\bf E}$ coincide. If we take $l=d$
and define $P_1,\ldots, P_d$
as an orthonormal basis of $\mathbb{U}_0^\perp$, then we get the second equality of \eqref{eq:zero_pert2}. This case
corresponds to a special ``noise matrix'' ${\bf E}$ with ${\bf E}{\bf E}^{\rm T}$ proportional to ${\bf I}$.
\end{remark}

Consider more explicitly the case where both  ${\bf H}{\bf E}^{\rm T}$ and ${\bf H}^{\rm T}{\bf E}$ are zero matrices.
Then  ${\bf P}^{\perp}_0(\delta)={\bf P}_0^{\perp}$ under the conditions of Theorem
\ref{theor:zero_pert}.

We call  matrices ${\bf H}$ and ${\bf E}$ {\it right-orthogonal} if
${\bf H}{\bf E}^{\rm T}={\bf 0}_{LL}$, or, equivalently,
${\bf E}{\bf H}^{\rm T}={\bf 0}_{LL}$. If ${\bf H}{\bf E}^{\rm T}={\bf 0}_{LL}$ and
${\bf H}^{\rm T}{\bf E}={\bf 0}_{KK}$ (the latter equality means that ${\bf H}$ and ${\bf E}$ are {\it left-orthogonal}),
then the matrices are called {\it biorthogonal}. The biorthogonality condition corresponds to the notion
of weak separability in Singular
Spectrum Analysis (see \cite[Sections 1.5 and 6.1]{GNZh01}).

The following statement elucidates the notion of left orthogonality.
\begin{lemma}
\label{lem:weak_sep}
Let ${\bf H}$ and ${\bf E}$ be non-zero $L\times K$ matrices.\\
1. If ${\bf H}$ and ${\bf E}$ are left-orthogonal, then $0\in \Sigma$,
${\bf P}_0{\bf E}={\bf E}$, and ${\bf P_\mu}{\bf E}={\bf 0}_{KL}$ for all non-zero $\mu\in \Sigma$.\\
2. If ${\bf P_\mu}{\bf E}={\bf 0}_{KL}$ for all non-zero $\mu\in\Sigma$, then ${\bf H}$ and ${\bf E}$ are left-orthogonal.
\end{lemma}

\begin{remark}
Denote $\Sigma_{\bf H}$ and $\Sigma_{\bf E}(\delta)$ the sets of positive eigenvalues of the operators ${\bf A}={\bf H}{\bf H}^{\rm T}$ and
$\delta^2{\bf E}{\bf E}^{\rm T}=\delta^2{\bf A}^{(2)}$, respectively.

Suppose that $\bf H$ and $\bf E$ are biorthogonal and assume additionally that
$\Sigma_{\bf H}\cap \Sigma_{\bf E}(\delta)=\varnothing$. (Note that this corresponds to the strong separability in
SSA, see \cite[\S1.5]{GNZh01}).
Then it is easy to see that  SVD of the matrix ${\bf H}+\delta {\bf E}$ is the sum of SVDs of the matrices
${\bf H}$ and $\delta {\bf E}$. Consequently, both singular values and the related
singular vectors of the operator ${\bf H}$ do not change
under the perturbation ${\bf H}\mapsto {\bf H}+\delta{\bf E}$.

Imposing also the condition $\delta_0^2\|{\bf E}{\bf E}^{\rm T}\|=\mu_{\min}$, we obtain that
${\bf P}_0^{\perp}(\delta)={\bf P}_0^{\perp}$ for all $\delta\in (\delta_0,\delta_0)$.
This condition tells us that each positive eigenvalue of the matrix ${\bf H}{\bf H}^{\rm T}$
is greater than all eigenvalues of the matrix $\delta^2 {\bf E}{\bf E}^{\rm T}$.



This means that in the case of biorthogonality, we have
${\bf P}_0^{\perp}(\delta)={\bf P}_0^{\perp}$ for all $\delta\in (\delta_0,\delta_0)$
under the
necessary and sufficient condition $\delta_0^2\|{\bf E}{\bf E}^{\rm T}\|=\mu_{\min}$
rather than under the sufficient
condition $\delta_0^2\|{\bf E}{\bf E}^{\rm T}\|=\mu_{\min}/2$ of Theorem \ref{theor:zero_pert}.
\end{remark}

If the matrices ${\bf H}$ and ${\bf E}$ are right-or\-tho\-go\-nal, then
${\bf B}(\delta)=\delta^2 {\bf A}^{(2)}$, the condition $\delta^2\|{\bf E}{\bf E}^{\rm T}\|<\mu_{\min}/2$ provides the
validity of Theorem  \ref{theor:full_decomp},
and, due to \eqref{eq:pert_0_B}, the expansion of the perturbed projector takes the form
\begin{gather}
\begin{split}
&{\bf P}_0^{\perp}(\delta)={\bf P}_0^{\perp}
\\
&+\!
\sum_{p=1}^\infty  (-1)^p\,\delta^{2p}\!\sum_{\substack{l_1+\ldots+l_{p+1}=p\\  l_j\geq 0}}
\!{{\bf S}}_{0}^{(l_1)}{\bf A}^{(2)}\ldots
{\bf A}^{(2)}{{\bf S}}_{0}^{(l_{p+1})}.
\end{split}
\label{eq:pert_0_B_right}
\end{gather}

\subsubsection{General upper bounds}
Roughly speaking, Theorem \ref{theor:zero_pert} shows that under the conditions
$B(\delta_0)=\mu_{\min}/2$ and $|\delta|<\delta_0$, the
equalities ${\bf S}_0{\bf B}(\delta){\bf P}_0={\bf 0}$ and ${\bf P}_0^{\perp}(\delta)={\bf P}_0^{\perp}$ are equivalent.
Moreover,
\begin{gather*}
\big\|{\bf S}_0{\bf B}(\delta){\bf P}_0\big\|\leq \frac{\|{\bf B}(\delta)\|}{\mu_{\min}}\leq
\frac{B(\delta)}{\mu_{\min}}\leq 1/2.
\end{gather*}
These considerations give rise to the supposition that $\big\|{\bf S}_0{\bf B}(\delta){\bf P}_0\big\|$ can serve as a natural
measure of the proximity
${\bf P}_0^{\perp}(\delta)\approx{\bf P}_0^{\perp}$.

\begin{theorem}
\label{theor:gen_upper}
If $\delta_0>0$ and $\big\|{\bf B}(\delta)\big\|/\mu_{\min}<1/4$ for all $\delta\in (-\delta_0,\delta_0)$,
then
\begin{gather}
\label{eq:gen_upper}
\big\|{\bf P}_0^{\perp}(\delta)-{\bf P}_0^{\perp}\big\|\leq
4C\ \frac{\|{\bf S}_0{\bf B}(\delta){\bf P}_0\|}{1-4\|{\bf B}(\delta)\|/{\mu_{\min}}},
\end{gather}
where  $C=e^{1/6}/\sqrt\pi$.
\end{theorem}

Let us discuss the conditions which provide the right-hand side of \eqref{eq:gen_upper} to be small.

\vspace{2mm}
\noindent
{\it The case of big signal and small noise matrices.}
$\phantom{}$\\
We start with a condition formulated
in terms of eigenvalues of matrices
${\bf A}={\bf H}{\bf H}^{\rm T}$ and ${\bf A}^{(2)}={\bf E}{\bf E}^{\rm T}$. Denote
\begin{gather*}
\Theta_1= \sqrt{\frac{\nu_{\max}}{\mu_{\max}}} \quad {\mbox{and}} \quad \Theta_2=\frac{\mu_{\max}}{\mu_{\min}}\ ,
\end{gather*}
where $\nu_{\max}=\|{\bf A}^{(2)}\|$.
Note that $\Theta_1$ is a sort of  the ``noise-to-signal ratio''. Since $\Theta_2=\|{\bf A}\|\,\|{\bf S}_0\|$
and since ${\bf S}_0$ is the pseudo-inverse  to  ${\bf A}$, then  $\Theta_2$
can be called as  the ``pseudo-condition number'' of the matrix  ${\bf A}$.
\begin{corollary}
\label{cor:coarse_proj}
Under the conditions of Theorem \ref{theor:gen_upper},
\begin{gather}
\label{eq:coarse_proj}
\big\|{\bf P}_0^{\perp}(\delta)-{\bf P}_0^{\perp}\big\|\leq 4C\, \frac{\|{\bf B}(\delta)\|}{\mu_{\min}} \
\frac{1}{1-4\|{\bf B}(\delta)\|/{\mu_{\min}}}
\end{gather}
with
\begin{gather}
\label{eq:B/mu_coarse}
\frac{\|{{\bf B}}(\delta)\|}{\mu_{\min}}\leq \frac{B(\delta)}{\mu_{\min}}
\leq 2|\delta|\,{\Theta_1 \Theta_2}+\delta^2 \Theta_1^2\Theta_2.
\end{gather}
\end{corollary}

\begin{remark}
1. The inequality \eqref{eq:coarse_proj} shows that the difference ${\bf P}_0^{\perp}(\delta)-{\bf P}_0^{\perp}$ is small if
the norm of the perturbation operator ${\bf B}(\delta)$ is substantially smaller than the minimal positive eigenvalue of
the matrix ${\bf H}{\bf H}^{\rm T}$.\\
2. Inequalities \eqref{eq:coarse_proj},  \eqref{eq:B/mu_coarse} jointly yield
sufficient conditions for the close proximity of the projectors
${\bf P}_0^{\perp}(\delta)$ and ${\bf P}_0^{\perp}$ in terms of the eigenvalues of the matrices
${\bf H}{\bf H}^{\rm T}$ and ${\bf E}{\bf E}^{\rm T}$. Roughly speaking, for fixed $\delta$ the difference
${\bf P}_0^{\perp}(\delta)-{\bf P}_0^{\perp}$ is small if $\Theta_1\ll 1$
(this means that the signal matrix ${\bf H}$ is ``big'' and the noise matrix ${\bf E}$ is ``small'')
and if the positive spectrum of the matrix
${\bf H}{\bf H}^{\rm T}$ is not wide-spread in the sense that the quotient $\Theta_2$ is not very large.

In particular, if
$2|\delta|\,{\Theta_1 \Theta_2}+\delta^2\,\Theta_1^2\Theta_2 \leq \varepsilon<1/4$, then
\begin{gather}
\label{eq:bigsmall_eps}
\big\|{\bf P}^{\perp}_0(\delta)-{\bf P}_0^{\perp}\big\|\leq 4C\frac{\varepsilon}{1-4\varepsilon}\,.
\end{gather}
\end{remark}

\noindent
{\it The case of the approximate orthogonalities.}\\
The upper bound in the inequalities \eqref{eq:coarse_proj}, \eqref{eq:B/mu_coarse}
is rather rough, since it does not incorporate the
orthogonality properties of matrices ${\bf H}$ and ${\bf E}$.
(On the other hand, this upper bound shows
that ${\bf P}_0^{\perp}(\delta)$ can be close to ${\bf P}_0^{\perp}$ without any orthogonalities.)
Even
if these matrices are biorthogonal, still
the right-hand side of \eqref{eq:coarse_proj} remains positive: in this case \eqref{eq:coarse_proj} takes a form
\begin{gather*}
\big\|{\bf P}_0^{\perp}(\delta)-{\bf P}_0^{\perp}\big\|\leq 4C\delta^2\, \frac{\nu_{\max}}{\mu_{\min}} \
\frac{1}{1-4\delta^2\nu_{\max}/{\mu_{\min}}}\,.
\end{gather*}

To improve this inconvenience we start with the following assertion concerning minimal principal angles between linear
spaces.

\begin{proposition}
\label{prop:min_angle}
Consider matrices ${\bf M}_1$, ${\bf M}_2$ acting from $\mathbb{R}^K$ onto $\mathbb{R}^L$.
Denote
$\theta_{\min}$ the minimal principal angle between subspaces $\mathbb{U}_1$ and $\mathbb{U}_2$, that are spanned
by the columns
of matrices ${\bf M}_1$ and ${\bf M}_2$. Lastly, let $\sigma_1^{(\min)}$, $\sigma_2^{(\min)}$ stand for the minimal
singular values of ${\bf M}_1$, ${\bf M}_2$, respectively.
Then
\begin{gather}
\begin{split}
&\sigma_1^{(\min)}\sigma_2^{(\min)}\cos(\theta_{\min})
\\
&\leq
\big\|{\bf M}_1^{\rm T}{\bf M}_2\big\|\leq \big\|{\bf M}_1\big\|\,\big\|{\bf M}_2\big\|\,\cos(\theta_{\min}).
\end{split}
\label{eq:min_angle}
\end{gather}
\end{proposition}

The following assertion follows from  Theorem \ref{theor:gen_upper} and Proposition \ref{prop:min_angle}.
Denote $\theta_{r}$ ($\theta_{l}$)  minimal principal angles between subspaces spanned by rows (columns) of matrices
${\bf H}$ and ${\bf E}$.
\begin{corollary}
\label{cor:coarse_orth}
Under conditions of Theorem \ref{theor:gen_upper},
\begin{gather}
\label{eq:coarse_orth}
\big\|{\bf P}_0^{\perp}(\delta)-{\bf P}_0^{\perp}\big\|\leq
4C\ \frac{\|{\bf S}_0{\bf B}(\delta)\|}{1-4\|{\bf B}(\delta)\|/{\mu_{\min}}}
\end{gather}
and
\begin{gather}
\label{eq:SBdelta_upper}
\big\|{\bf S}_0{\bf B}(\delta)\big\|\leq |\delta|\,{\Theta_1\Theta_2}
\left(2\cos(\theta_{r})+|\delta|\,\Theta_1\,\cos(\theta_l)\right).
\end{gather}
\end{corollary}

\begin{remark}
\label{rem:appr_orth}
1. Both $\cos(\theta_r)$ and $\cos(\theta_l)$ serve as
proper measures of right and left orthogonalities for matrices ${\bf H}$ and ${\bf E}$.\\
2. Using inequalities  \eqref{eq:coarse_orth} and \eqref{eq:SBdelta_upper} we
obtain  sufficient conditions for the proximity ${\bf P}_0^{\perp}(\delta)\approx {\bf P}_0^{\perp}$ in the
case when the ``noise-to-signal ratio'' $\Theta_1$ is not small. Thus, the ``pseudo-condition number'' $\Theta_2$
of the matrix  ${\bf H}{\bf H}^{\rm T}$ should not be very large and both pairs of linear spaces produced by matrices
${\bf H}$ and ${\bf E}$ should be almost orthogonal.\\
\end{remark}

\subsubsection{Main terms of the approximations}
\label{sssect:main_terms}
In this section we discuss special refinements of inequalities \eqref{eq:gen_upper}, \eqref{eq:coarse_proj},
and \eqref{eq:coarse_orth}. More precisely, we present operators that can be considered as the main terms of the
difference ${\bf P}_0^{\perp}(\delta)-{\bf P}_0^{\perp}$ in the case when this difference is small by norm.
The results show that inequalities \eqref{eq:gen_upper}, \eqref{eq:coarse_proj}, and \eqref{eq:coarse_orth} produce different
(though connected with each other) versions of the main terms.

\vspace{2mm}
\noindent
{\it The case of  big signal and small noise matrices.}$\phantom{}$\\
Let us start with the inequality \eqref{eq:coarse_proj}.

\begin{theorem}
\label{theor:coarse_proj}
Under conditions of Theorem \ref{theor:gen_upper},
\begin{gather}
\begin{split}
&\big\|{\bf P}_0^{\perp}(\delta)-{\bf P}_0^{\perp}-{\bf W}_1(\delta)\big\|\\
&\leq 16\,C
\left(\frac{\|{\bf B}(\delta)\|}{\mu_{\min}}\right)^2
\frac{1}{1-4\|{\bf B}(\delta)\|/{\mu_{\min}}}\ ,
\end{split}
\label{eq:coarse_proj_mainterm}
\end{gather}
where ${\bf W}_1(\delta)$ is given by \eqref{eq:first_term} and
$C=e^{1/6}/\sqrt{\pi}$.
\end{theorem}

Let us compare \eqref{eq:coarse_proj_mainterm} and \eqref{eq:coarse_proj} with the help of the inequality \eqref{eq:B/mu_coarse}.
If ${\bf B}(\delta)/\mu_{\min}$ is small,
then the inequality \eqref{eq:coarse_proj_mainterm} shows that the operator
${\bf W}_1(\delta)$ can be considered as the main term of the difference ${\bf P}_0^{\perp}(\delta)-{\bf P}_0^{\perp}$.
In particular, $\big\|{\bf P}_0^{\perp}(\delta)-{\bf P}_0^{\perp}-{\bf W}_1(\delta)\big\|\ll
\big\|{\bf P}_0^{\perp}(\delta)-{\bf P}_0^{\perp}\big\|$
if the ``noise-to-signal ratio'' $\Theta_1$ is small and the ``pseudo-condition number'' $\Theta_2$
is not very large.

For example, if
$2|\delta|\,{\Theta_1 \Theta_2}+\delta^2\,\Theta_1^2\Theta_2 \leq \varepsilon<1/4$, then
\begin{gather*}
\big\|{\bf P}_0^{\perp}(\delta)-{\bf P}_0^{\perp}-{\bf W}_1(\delta)\big\|\leq 16C\frac{\varepsilon^2}{1-4\varepsilon}\ ,
\end{gather*}
while $\|{\bf P}_0^{\perp}(\delta)-{\bf P}_0^{\perp}\|$ satisfies \eqref{eq:bigsmall_eps}.

\vspace{2mm}
\noindent
{\it The case of approximate orthogonalities.}$\phantom{}$\\
Let us pass to the upper bound \eqref{eq:coarse_orth}.

We start with a simple remark. Denote ${\bf A}^{(2)}_0= {\bf P}_0{\bf A}^{(2)}{\bf P}_0$ and
assume that $B(\delta_0)<\mu_{\min}$.
Then it is easy to
check that the operator ${\bf I}-\delta^2 {\bf A}_0^{(2)}/\mu$ is invertible  for any positive $\mu\in \Sigma$ and
for any $\delta\in (-\delta_0,\delta_0)$ since
\begin{gather*}
\delta^2\, \frac{\big\|{\bf A}_0^{(2)}\big\|}{\mu}\leq \delta_0^2 \, \frac{\nu_{\max}}{\mu_{\min}}\leq
\frac{B(\delta_0)}{\mu_{\min}}<1.
\end{gather*}

\begin{theorem}
\label{theor:precise_proj}
Assume that $\delta_0>0$, $B(\delta_0)=\mu_{\min}/4$ and $|\delta|< \delta_0$. Denote
\begin{gather}
\label{eq:Ldelta_short}
{\bf L}_1(\delta)=\sum_{\mu>0}\frac{{\bf P}_\mu {\bf B}(\delta){\bf P}_0}{\mu}
\left({\bf I}-\delta^2{\bf A}^{(2)}_0/{\mu}\right)^{-1}
\end{gather}
and ${\bf L}(\delta)={\bf L}_1(\delta)+{\bf L}_1^{\rm T}(\delta)$.
Then
\begin{gather}
\label{eq:res_L_delta}
\big\|{\bf P}_0^{\perp}(\delta)\!-\!{\bf P}_0^{\perp}\!-\!{\bf L}(\delta)\big\|\!\leq\!
16C \frac{\|{\bf S}_0{\bf B}(\delta)\|\,\|{\bf S}_0{\bf B}(\delta){\bf P}_{0}\|}{1-4\|{\bf B}(\delta)\|/\mu_{\min}}
\end{gather}
where  $C=e^{1/6}/\sqrt\pi$.
\end{theorem}

The operator ${\bf L}(\delta)$ admits another representation. Denote
\begin{gather*}
{\bf K}_1(\delta)=\sum_{\mu>0}\frac{{\bf P}_\mu {\bf B}(\delta){\bf A}^{(2)}_0}{\mu^2}
\left({\bf I}-\delta^2{\bf A}^{(2)}_0/\mu\right)^{-1}
\end{gather*}
and put ${\bf K}(\delta)={\bf K}_1(\delta)+{\bf K}^{\rm T}_1(\delta)$.
\begin{proposition}
\label{prop:L_delta}
Under the conditions and notation of Theorem \ref{theor:precise_proj},
\begin{gather}
\label{eq:L_delta}
 {\bf L}(\delta)={\bf W}_1(\delta)+\delta^2{\bf K}(\delta),
\end{gather}
where  ${\bf W}_1(\delta)$ is defined in \eqref{eq:first_term}.
\end{proposition}

Let us discuss the result of Theorem \ref{theor:precise_proj}.
The upper bound \eqref{eq:coarse_orth} is reasonable in the situation when
$\|{\bf B}(\delta)\|/\mu_{\min}$ is not small and $\|{\bf S}_0{\bf B}(\delta)\|$ is small enough.

Assume that $\|{\bf S}_0{\bf B}(\delta)\|\approx\varepsilon$.
Since
\begin{gather*}
\|{\bf S}_0{\bf B}(\delta){\bf P}_0\|\!\leq\!\|{\bf S}_0{\bf B}(\delta)\|,
\end{gather*}
then the right-hand side of
\eqref{eq:res_L_delta}
is proportional to $\varepsilon^2$ while the right-hand side of \eqref{eq:coarse_orth}  is proportional to
$\varepsilon$. This means that  the operator ${\bf L}(\delta)$ can be considered as the main term of the difference
${\bf P}_0^\perp(\delta)-{\bf P}_0^\perp$.

\vspace{2mm}
\noindent
{\it General case.}\\
Let us study the general upper bound given in Theorem \ref{theor:gen_upper}.
It follows from the discussion of Section \ref{ssect:zero_pert} that the
inequality \eqref{eq:gen_upper} can give  weaker conditions for the
proximity ${\bf P}_0^{\perp}(\delta)\approx {\bf P}_0^{\perp}$ than inequalities \eqref{eq:coarse_proj}
and \eqref{eq:coarse_orth}.

Let us present a version of the main term of the difference ${\bf P}_0^\perp(\delta)-{\bf P}_0^\perp$ in this
general case.

\begin{theorem}
\label{theor:precise_proj_noise}
Under the conditions of Theorem \ref{theor:precise_proj},
\begin{gather}
\label{eq:res_T_delta}
\big\|{\bf P}_0^{\perp}(\delta)-{\bf P}_0^{\perp}-{\bf T}(\delta)\big\|\leq
16C\,\frac{\|{\bf S}_0{\bf B}(\delta){\bf P}_{0}\|^2}{1-4\|{\bf B}(\delta)\|/\mu_{\min}}\ ,
\end{gather}
where $C=e^{1/6}/\sqrt\pi$,  ${\bf T}(\delta)={\bf T}_1(\delta)+{\bf T}_1^{\rm T}(\delta)$,
\begin{gather}
\label{eq:T_1}
{\bf T}_1(\delta)= \sum_{i=0}^{\infty}(-1)^{i}\sum_{\mu_1,\ldots,\mu_i>0}\,
{\bf J}_i
\, {\bf L}_1(\delta)\, {\bf G}_i,
\end{gather}
$
{\bf J}_i={\bf J}_i(\mu_1,\ldots,\mu_i)=\prod_{k=1}^i \big({\bf P}_{\mu_k}{\bf B}(\delta)/{\mu_k}\big),
$
and
\begin{gather*}
{\bf G}_i={\bf G}_i(\mu_1,\ldots,\mu_i)=\prod_{k=1}^i\left({\bf I}-\delta^2{\bf A}^{(2)}_0/\mu_k\right)^{-1}.
\end{gather*}
\end{theorem}

\section{Subspace perturbations for  Hankel matrices of large size}
\label{sect:Hankel_pert}
Let $\mathrm{F}=(x_0,\ldots,x_{N-1},\ldots)$ and $\mathrm{E}=(e_0,\ldots,e_{N-1},\ldots)$. We treat $\mathrm{F}$
as a ``signal series'' and  $\mathrm{E}$ as a ``noise series''.
For fixed $N$ we choose the window length $L$
and construct two Hankel (``trajectory'') matrices
\begin{gather*}
{\bf H}={\bf H}_{K,L}=
\begin{pmatrix}
x_0      &x_1         &\ldots&x_{K-2}             &x_{K-1}     \\
\vdots  &\vdots    &\ddots&\vdots        &\vdots\\
x_{L-1} &x_L         &\ldots& x_{N-2}             &x_{N-1}     \\
\end{pmatrix}
\end{gather*}
and
\begin{gather*}
{\bf E}={\bf E}_{K,L}=
\begin{pmatrix}
e_0      &e_1         &\ldots&e_{K-2}              &e_{K-1}     \\
\vdots  &\vdots    &\ddots&\vdots        &\vdots\\
e_{L-1} & e_L        &\ldots& e_{N-2}              &e_{N-1}     \\
\end{pmatrix},
\end{gather*}
where $K=N-L+1$. In terms of Section \ref{sect:Zero_eigen},
${\bf H}$ serves as a signal matrix  and ${\bf E}$ is a noise matrix.

It is clear that ${\rm rank}\, {\bf H}\leq \min(L,K)$. As it was already mentioned, we are interested in
the case ${\rm rank}\, {\bf H}<\min(L,K)$.
To provide this condition, we assume that $\mathrm{F}$ is governed by a minimal Linear Recurrent Formula of order $d$.
Then ${\rm rank}\, {\bf H}=d$ for any $L$ and $K$ such that $\min(L,K)\geq d$.

Consider the perturbed series $\mathrm{F}(\delta)=\mathrm{F}+\delta \mathrm{E}$ and construct the corresponding
 Hankel matrix ${\bf H}(\delta)={\bf H}+\delta {\bf E}$.  Then we can apply all notation and
results of Section \ref{sect:Zero_eigen}
to this particular case of matrices ${\bf H}$ and ${\bf E}$.

In the present section we derive conditions providing
the convergence $\big\|{\bf P}_0^{\perp}(\delta)-{\bf P}_0^{\perp}\big\|\rightarrow 0$ as $N\rightarrow\infty$ and $L=L(N)$.
In what follows, we omit the dependence
of matrices, projectors, etc. on $N$ and $L$ in our notation.

\subsection{Two general asymptotic results}
\label{ssect:asympt}
We start with two assertions that follow from inequalities \eqref{eq:coarse_proj}, \eqref{eq:coarse_proj_mainterm} and
\eqref{eq:B/mu_coarse}.
Denote
\begin{gather*}
\Theta=\Theta_1\,\Theta_2 =\sqrt{\frac{\nu_{\max}}{\mu_{\max}}}\, \frac{\mu_{\max}}{\mu_{\min}}\, .
\end{gather*}

\begin{proposition}
\label{prop:big_NSR}
If $\Theta\rightarrow 0$ as $N\rightarrow \infty$, then
\begin{gather}
\label{eq:Diff_proj}
\limsup_{N}\,\Theta^{-1}\,\big\|{\bf P}_0^{\perp}(\delta)-{\bf P}_0^{\perp}\big\|
\leq 8C|\delta|
\end{gather}
and
\begin{gather}
\label{eq:Diff_proj_second}
\limsup_{N}\,\Theta^{-2}\,\big\|{\bf P}_0^{\perp}(\delta)-{\bf P}_0^{\perp}-\delta {\bf V}_{0}^{(1)}\big\|
\leq C'\delta^2,
\end{gather}
where $C=e^{1/6}/\sqrt{\pi}$  and $C'=2(32 C +1)$.
\end{proposition}

Thus we  see that $\big\|{\bf P}_0^{\perp}(\delta)-{\bf P}_0^{\perp}\big\|\rightarrow 0$
for any fixed $\delta$ if $\Theta$ tends to zero.
This means that we can treat $\Theta$ in the same sense as the perturbation parameter $\delta$ for fixed $N$.
Note that $\Theta\rightarrow 0$ only if $\Theta_1\rightarrow 0$.

At the same time the condition $\Theta_1\rightarrow 0$ is not necessary for the convergence
 $\big\|{\bf P}_0^{\perp}(\delta)-{\bf P}_0^{\perp}\big\|\rightarrow 0$.
 In the next proposition we show that this convergence occurs also
in the case $\Theta_1\asymp 1$
 if matrices ${\bf H}$ and ${\bf E}$ are asymptotically biorthogonal and
$\Theta_2$ is bounded from above.  (Note that for positive
sequences $a_n$ and $b_n$ we write  $a_n\asymp b_n$ iff $c_1\leq a_n/b_n\leq c_2$ for some constants $c_1,c_2>0$.)
Denote $\Delta=1/\limsup(\Theta \Theta_1)$.

\begin{proposition}
\label{prop:asymp_ort}
Suppose that $\Delta>0$, $\big\|{\bf H}{\bf E}^{\rm T}\big\|/\mu_{\min}\rightarrow 0$,
and $\big\|{{\bf S}}_{0}{\bf A}^{(2)}\big\|\rightarrow 0$ as $N\rightarrow \infty$.
Then for any
$\delta$ with $|\delta|<\delta_0\stackrel{\rm def}=\Delta/4$,
\begin{gather}
\label{eq:asymp_ort}
\big\|{\bf P}_0^{\perp}(\delta)\!-\!{\bf P}_0^{\perp}\big\|\!=\!
|\delta|\,O\!\left(
\big\|{\bf H}{\bf E}^{\rm T}\big\|/\mu_{\min}+|\delta| \big\|{\bf S}_0{\bf A}^{(2)}\big\|\right)
\end{gather}
and
\begin{gather}
\begin{split}
&\big\|{\bf P}_0^{\perp}(\delta)-{\bf P}_0^{\perp}-{\bf L}(\delta)\big\|
\\
&=
\delta^2\,O\!
\left(\left(\big\|{\bf H}{\bf E}^{\rm T}\big\|/\mu_{\min}+|\delta|\,\big\|{\bf S}_0{\bf A}^{(2)}\big\|\right)^2\right),
\label{eq:asymp_ort_L}
\end{split}
\end{gather}
provided that $N\geq N_0(\delta)$.

\end{proposition}

\begin{remark}
\label{rem:beorth}
1. The operator ${\bf L}(\delta)$ is defined in \eqref{eq:Ldelta_short}.\\
2. If $\Theta \Theta_1 \asymp 1$, then $\Theta \nrightarrow 0$ and $\delta_0<\infty$. Hence the result of Proposition
\ref{prop:asymp_ort} is valid only if  $|\delta|$ is bounded from above.
These restrictions are absent in the case
$\Theta\rightarrow 0$ considered in Proposition \ref{prop:big_NSR}.\\
3. Let $\Theta \Theta_1 \asymp 1$ and $\Theta_1\asymp 1$. Consequently, $\Theta_2$ is bounded from above.
Applying \eqref{eq:SBdelta_upper} to \eqref{eq:asymp_ort} we obtain that
$\big\|{\bf P}_0^{\perp}(\delta)-{\bf P}_0^{\perp}\big\|\rightarrow 0$ if both $\cos(\theta_r)$ and $\cos(\theta_l)$ tend to
zero. Moreover, \eqref{eq:asymp_ort} can be rewritten as
\begin{gather*}
\big\|{\bf P}_0^{\perp}(\delta)-{\bf P}_0^{\perp}\big\|=
|\delta|\,O\, \Big(\cos(\theta_r)+|\delta|\,\cos(\theta_l)\Big)
\end{gather*}
and the inequality \eqref{eq:asymp_ort_L} admits similar reformulation.
\end{remark}

\subsection{Examples. Rough upper bounds}
\label{ssect:hankel_coarse}
To illustrate propositions of the previous subsection we consider several types of ``signal'' and ``noise'' series. For a
series $\mathrm{F}=(f_0,\ldots,f_n,\ldots)$ we define the associated $L\times K$ trajectory matrix in the form
\begin{gather*}
{\bf F}={\bf F}_{K,L}=
\begin{pmatrix}
f_0      &f_1         &\ldots&f_{K-2}             &f_{K-1}     \\
\vdots  &\vdots    &\ddots&\vdots        &\vdots\\
f_{L-1} &f_L         &\ldots& f_{N-2}             &f_{N-1}     \\
\end{pmatrix}
\end{gather*}
and denote the maximal and minimal positive eigenvalues of the matrix ${\bf F}{\bf F}^{\rm T}$ as $\lambda_{\max}$ and
$\lambda_{\min}$, respectively.

In examples we consider the following four types of series.
\begin{enumerate}
\item
{\it A linear combination of increasing exponents}. In this case
\begin{gather}
\label{eq:exp_sum}
f_n=\beta_1a_1^n+\ldots+\beta_pa_p^n,
\end{gather}
where $\beta_j\neq 0$ for $1\leq j \leq p$ and $a_1>\ldots>a_{p}>1$. The series \eqref{eq:exp_sum} has rank $p$.
For short, we name this series {\it the series of exponential type}.
\item
{\it A polynomial series.} This series is defined by
\begin{gather}
\label{eq:polyn}
f_n=\gamma_pn^p+\gamma_{p-1}n^{p-1}+\ldots+\gamma_{1}n+\alpha_0,
\end{gather}
 where $\gamma_p\neq 0$. The rank of the series \eqref{eq:polyn} is $p+1$.
 \item
 {\it An oscillating series} with frequencies $\omega_l$, positive amplitudes $\gamma_l$
 and phases $\phi_l$.  The series has the form
 \begin{gather}
 \label{eq:osc}
 f_n=\sum_{l=1}^{p}\gamma_l \cos(2\pi\omega_l n+\phi_l),
  \end{gather}
where $\omega_l\in [0,1/2]$ and $\omega_l<\omega_j$ for $l<j$. The rank of the series is $2p$ if
$\omega_j\in (0,1/2)$ for all $j$. If $\omega_1=0$ and $\omega_p<1/2$ (or if $\omega_1>0$ and $\omega_p=1/2$),
then the rank is $2p-1$. If $\omega_1=0$ and $\omega_p=1/2$, then the rank is $2p-2$.
  \item
 {\it A linear stationary random series}. By definition,
\begin{gather}
\label{eq:linear_process_f}
f_n=\sum_{j=-\infty}^\infty c_{j}\,\varepsilon_{j+n},
 \end{gather}
 where $\varepsilon_n$  is the sequence
 of independent random variables with $\mathbb{E}\varepsilon_n=0$, $\mathbb{D}\varepsilon_n=1$,
  and $\mathbb{E}|\varepsilon_n|^{3}<\infty$.
  We assume that $S\stackrel{\rm def}=\sum_j |c_j|<\infty$ and $\sum_j c_j^2=1$. Then
 $\mathbb{E}f_n=0$ and $\mathbb{D}f_n=1$.
   Note that the series $(f_0,\ldots,f_N,\ldots)$ is not a series of finite rank.
\end{enumerate}

Let us discuss asymptotic
properties of eigenvalues
$\lambda_{\max}$ and $\lambda_{\min}$ corresponding to series
\eqref{eq:exp_sum}--\eqref{eq:linear_process_f} as $N\rightarrow \infty$ and $L=L(N)$.

\begin{lemma}
\label{summ:series}
1. 
For the series of exponential type
there exist positive $T_{\max}^{(a)},T_{\min}^{(a)}$ such that
$\lambda_{\max}/a_{1}^{2N}\rightarrow T_{\max}^{(a)}$ and $\lambda_{\min}/a_{p}^{2N}\rightarrow T_{\min}^{(a)}$ as $\min(L,K)\rightarrow \infty$.
The analogous result  holds if either $L=L_0={\rm const}>p$ or
$K=K_0={\rm const}>p$.\\
2.
For the  polynomial series,

a) if $L/N\rightarrow \alpha\in (0,1)$ then
$\lambda_{\max}/N^{2p+2}\!\rightarrow\! \Theta_{\max}$ and
 $\lambda_{\min}/N^{2p+2}\!\rightarrow\! \Theta_{\min}$ for some positive $\Theta_{\max}, \Theta_{\min}$;

b) if
 either $L$ or $K$ is a constant greater than $p+1$, then $\lambda_{\max}/N^{2p+1}\rightarrow \Theta_{\max}>0$ and
 $\lambda_{\min}/N\rightarrow \Theta_{\min}>0$.\\
3.
For the  oscillating series,
$\lambda_{\max}/LK\rightarrow \Lambda_{\max}>0$ and $\lambda_{\min}/LK\rightarrow \!\Lambda_{\min}\!>0$ as $\min(L,K)\rightarrow \!\infty$.
The analogous result  holds if either
$L=L_0={\rm const}>d$ or $K\!=\!K_0={\rm const}\!>\!d$.
\end{lemma}
We omit elementary but  laborsome proofs of these assertions.
\begin{remark}
\label{rem:lambda_amax}
1. Positive constants $T_{\max}^{(a)},T_{\min}^{(a)}$, $\Theta_{\max}, \Theta_{\min}$ and
$\Lambda_{\max},\Lambda_{\min}$ can be  found in the explicit form and
 depend both on the parameters of the series and on the behavior of $L=L(N)$.\\
2. For the oscillating series  \eqref{eq:osc} we have $\lambda_{\max}\asymp \lambda_{\min}\asymp N^2$
if $L/N\rightarrow \alpha \in (0,1)$ and $\lambda_{\max}\asymp \lambda_{\min}\asymp N$
 if either $L$ or $K$ does not depend on  $N$.\\
3. Parameters $a_j$ of the series \eqref{eq:exp_sum} are taken positive only for  convenience. All results
concerning the series of exponential type remain valid under assumption $|a_1|>\ldots>|a_p|$.
 \end{remark}
\begin{lemma}
\label{summ:series_rand}
Consider the linear stationary series \eqref{eq:linear_process_f}.\\
1.  If $L\rightarrow \infty$, then
 there exists an absolute constant $\gamma_0$ such that almost surely
 \begin{gather}
\label{eq:uneq_lm_max_stat0}
\limsup_N\frac{\lambda_{\max}(\omega)}{N\ln N}\leq \gamma_0\,S.
\end{gather}
2. Assume that  $L=L_0={\rm const}$ and denote
$R_f(\,\cdot\,)$ the covariance function of the series \eqref{eq:linear_process_f}.
Then
$\lambda_{\max}/N$ tends almost surely to the maximal eigenvalue $\sigma_{\max}$ of the matrix ${\bf \Sigma}=\big\{R_f(i-j)\big\}_{0\leq i,j\leq L_0-1}$ and
$\lambda_{\min}/N$ tends almost surely to the minimal eigenvalue $\sigma_{\min}$ of ${\bf \Sigma}$.
\end{lemma}
Both assertions of Lemma \ref{summ:series_rand} follow from results
  given in  Section \ref{sssect:LP} of Appendix 2.

\begin{remark}
\label{rem:gamma_prim}
1. The inequality \eqref{eq:uneq_lm_max_stat0} means that
for almost all $\omega$ and for any $\gamma'>\gamma_0$ we get
\begin{gather}
\label{eq:uneq_lm_max_stat}
\lambda_{\max}(\omega)<\gamma'S\,N\ln N
 \end{gather}
 provided that $N>N_0(\omega,\gamma')$.\\
2. If $L=L_0$, then the statement of Lemma \ref{summ:series_rand} holds without the
condition   $\sup_n\mathbb{E}|\varepsilon_n|^{3}<\infty$.\\
3. If  $L=L_0$ and $f_n\!=\!\varepsilon_n$ is the ``white noise'' series, then  ${\bf F}{\bf F}^{\rm T}/K\!\rightarrow\! {\bf I}_{L_0}$,
$\lambda_{\max}/N\rightarrow1$,
  and $\lambda_{\min}/N\rightarrow 1$ a.s.
\end{remark}
\begin{remark}
Since the series \eqref{eq:exp_sum}--\eqref{eq:osc}  are of finite rank, they can serve both as ``signal series'' and
as ``noise series''.
Then we use notation $x_n$ or $e_n$  instead of $f_n$. In the same manner we use notation $\mu_{\max}$ or
$\nu_{\max}$ instead of $\lambda_{\max}$ and so on.
The stationary series can be chosen only as a ``noise series''.
\end{remark}

\subsubsection{Signals of exponential type}
\label{ssect:sign_incr_exp}

Consider the signal series defined by \eqref{eq:exp_sum} and some ``noise series'' $\mathrm{E}$.
The following proposition gives a sufficient condition for the
convergence $\big\|{\bf P}_0^{\perp}(\delta)-{\bf P}_0^{\perp}\big\|\rightarrow 0$ as $N\rightarrow \infty$ in terms of
$\nu_{\max}$, $a_1$, and $a_p$.
Denote $\tau=a_1/a_p^2$.
\begin{proposition}
\label{prop:sum_exp}
Let $\Xi_a=\sqrt{T_{\max}^{(a)}}/T_{\min}^{(a)}$.
If  $\nu_{\max} \tau^{2N}=o(1)$ and $\delta\in \mathbb{R}$, then
\begin{gather*}
\limsup_N\,\nu_{\max}^{-1/2}\tau^{-N}\big\|{\bf P}_0^{\perp}(\delta)-{\bf P}_0^{\perp}\big\|
\leq 8C\,\Xi_a\, |\delta|
\quad
\rm{and}\\
\limsup_N\,\nu_{\max}^{-1}\tau^{-2N}
\big\|{\bf P}_0^{\perp}(\delta)-{\bf P}_0^{\perp}-\delta {\bf V}_0^{(1)}\big\|
\leq
C'\,\Xi_a^2\, \delta^2.
\end{gather*}
\end{proposition}

\begin{remark}
If we take $p=1$ in  \eqref{eq:exp_sum}, then $\mu_{\max}=\mu_{\min}\asymp a_1^{2N}$, $\tau=a_1^{-N}<1$,
and the condition of Proposition \ref{prop:sum_exp} reduces to $\nu_{\max}=o(a_1^{2N})$.
\end{remark}

Now let us present several examples related to Proposition \ref{prop:sum_exp}.
\begin{example}
\label{ex:sum_exp_sign}
\item[1.]
Let $e_n=\sum_{l=1}^{m}\gamma_l b_l^n$ with $\gamma_l\neq 0$ and let $b=\max_{1\leq l\leq m}\,|b_l|$. Assume that $b>1$.
Then $\nu_{\max}\sim T^{(b)}_{\max} b^{2N}$.
If $b\tau<1$, then for any $\delta$
\begin{gather*}
\limsup_N \,(b\tau)^{-N}\,\big\|{\bf P}_0^{\perp}(\delta)-{\bf P}_0^{\perp}\big\|\leq
8C\,\Xi_a\sqrt{T_{\max}^{(b)}}\ |\delta|
\quad \mbox{\rm and} \\
\limsup_N \,(b\tau)^{-2N}\,\big\|{\bf P}_0^{\perp}(\delta)-{\bf P}_0^{\perp}-\delta {\bf V}_0^{(1)}\big\|\leq
C'\,\Xi_a^2\,T_{\max}^{(b)}\ \delta^2.
\end{gather*}
\item[2.]
Consider a polynomial ``noise'' of order $m$ defined as in \eqref{eq:polyn}.

a)
In view of Lemma \ref{summ:series}, $\nu_{\max}\sim \Theta_{\max} N^{2m+2}$ in the case $L/N\rightarrow \alpha\in (0,1)$.
If $\tau<1$, then for any $\delta$
\begin{gather*}
\begin{split}
&\limsup_N\,N^{-m-1}\tau^{-N}\,\big\|{\bf P}_0^{\perp}(\delta)-{\bf P}_0^{\perp}\big\|\leq
8C\,\Xi_a\sqrt{\Theta_{\max}}\ |\delta|
\ \
\mbox{\rm and} \\
&\!\limsup_N N^{-2m-2}\tau^{-2N} \!\big\|{\bf P}_0^{\perp}(\delta)\!-\!{\bf P}_0^{(-)}\!-\!\delta {\bf V}_0^{(1)}\big\|\!\leq\!
C'\! \Xi_a^2 \Theta_{\max} \delta^2.
\end{split}
\end{gather*}

b)
If  $\min(L,K)={\rm const}$, then $\nu_{\max}\sim \Theta_{\max} N^{2m+1}$.
If $\tau<1$, then for any $\delta$
\begin{gather*}
\begin{split}
&\limsup_N\,N^{-m-1/2}\tau^{-N}\,\big\|{\bf P}_0^{\perp}(\delta)-{\bf P}_0^{\perp}\big\|\leq
8C\,\Xi_a\sqrt{\Theta_{\max}}\ |\delta|
\ \
\mbox{\rm  and} \\
&\!\limsup_N N^{-2m-1}\tau^{-2N} \!\big\|{\bf P}_0^{\perp}(\delta)\!-\!{\bf P}_0^{(-)}\!-\!\delta {\bf V}_0^{(1)}\big\|\!
\leq\!
C'\! \Xi_a^2 \Theta_{\max} \delta^2.
\end{split}
\end{gather*}
\item[3.]
\label{ex3:exp_cos_coarse}
Let the oscillating ``noise series'' be defined  as in of \eqref{eq:osc}.
Then $\nu_{\max}\sim \Lambda_{\max} LK$.
If $\tau<1$,  then for any $\delta$
\begin{gather*}
\begin{split}
&\limsup_N\,(LK)^{-1/2} \tau^{-N}
\big\|{\bf P}_0^{\perp}(\delta)-{\bf P}_0^{\perp}\big\|\leq
8C\,\Xi_a\sqrt{\Lambda_{\max}}\ |\delta|
\ \
\mbox{\rm and} \\
&\limsup_N (LK)^{-1} \tau^{-2N}
\big\|{\bf P}_0^{\perp}(\delta)\!-\!{\bf P}_0^{\perp}\!-\!\delta {\bf V}_0^{(1)}\big\|
\!\leq\!
C' \Xi_a^2 \Lambda_{\max} \delta^2.
\end{split}
\end{gather*}
\item[4.]
The case of random stationary ``noise series'' is studied in Example  \ref{ex:WN_LeqK}  and Proposition
\ref{prop:gen_wight_noise_coarse}.
\end{example}

\begin{remark}
\label{rem:reduced_rate}
It is easy to see that Proposition \ref{prop:sum_exp} (and therefore, Proposition \ref{prop:big_NSR}) can give
only an upper bound of the true rate of convergence ${\bf P}_0^{\perp}(\delta)\rightarrow {\bf P}_0^{\perp}$.
For example, for $x_n=a^n$ with $a>1$ and the constant ``noise series''
$e_n\equiv 1$,
we obtain the equality
$\big\|{\bf P}_0^{\perp}(\delta)-{\bf P}_0^{\perp}\big\|=O\big(Na^{-N}\big)$ in the case $L\sim K$.

On the other hand, Proposition \ref{prop:exp_const_proj_res} of Section \ref{ssect:exp_const_prec} affirms that
$\big\|{\bf P}_0^{\perp}(\delta)-{\bf P}_0^{\perp}\big\|\asymp\sqrt{N} a^{-N}$.
 The reason of this drawback is
that Proposition \ref{prop:big_NSR} can ignore the existing asymptotic orthogonalities of signal and noise matrices.
\end{remark}

\subsubsection{Oscillating signal series}
\label{sssect:osc_ser}
Consider the signal of oscillating  type
 with different frequencies $\omega_{l1}\in [0,1/2]$, positive amplitudes $\alpha_l$,
 and phases $\phi_{l1}$, $l=1,\ldots,p$.

 \begin{example}
\label{ex:osc_ser}
\item[1.]
Assume that  the ``noise series'' is also oscillating. Namely, let
\begin{gather*}
e_n=\sum_{j=1}^{m}\beta_j \cos(2\pi\omega_{j2} n+\phi_{j2}),
\end{gather*}
where $\beta_j>0$ and $\omega_{j1}\in [0,1/2]$.
 In addition, let $\omega_{l1}\neq \omega_{j2}$ for all $l,j$.

If $\min(L,K)\rightarrow \infty$, then $\mu_{\max}\sim \Lambda_{\max}^{(s)}LK$,
$\mu_{\min}\sim \Lambda_{\min}^{(s)}LK$, $\nu_{\max}\sim \Lambda_{\max}^{(n)} LK$, and
$\Theta\Theta_1=\nu_{\max}/\mu_{\min}\rightarrow 1/\Delta\stackrel{\rm def}=\Lambda_{\max}^{(n)}/\Lambda_{\min}^{(s)}$.
Consequently, the result of Proposition \ref{prop:asymp_ort} is valid for any
$\delta$ such that $|\delta|<\delta_0=\Delta/4$  provided that $N$ is  big enough.


  Calculations  show that
$
 \|{\bf H}{\bf E}^{\rm T}\|/\mu_{\min}=O(1/K)
 $
 and
 $
 \big\|{\bf S}_0{\bf A}^{(2)}\big\|=O(1/L).
 $
Therefore, it follows from \eqref{eq:asymp_ort}  and  \eqref{eq:asymp_ort_L} that
\begin{gather*}
\begin{split}
&\big\|{\bf P}_0^{\perp}(\delta)-{\bf P}_0^{\perp}\big\|=|\delta|\,O\big(1/\min(L,K)\big) \ \ \mbox{and}\\
&\big\|{\bf P}_0^{\perp}(\delta)-{\bf P}_0^{\perp}-{\bf L}(\delta)\big\|=\delta^2\,O\big(1/\min(L^2,K^2)\big).
\end{split}
\end{gather*}
\item[2.]
As in Example \ref{ex:sum_exp_sign},
the important case  of oscillation signal and  random stationary ``noise series'' is studied in
 Example \ref{ex:WN_LeqK}  and Proposition
\ref{prop:gen_wight_noise_coarse}.
\end{example}

\subsubsection{Random stationary  ``noise series''}
\label{ssect:WN_N}
Let $\varepsilon_n$ ($n\geq 0$) be a sequence of independent random variables defined on a probability space $(\Omega, {\cal F}, \mathbb{P})$.
Assume that $\mathbb{E}\varepsilon_n=0$,
$\mathbb{D}\varepsilon_n=1$, and $\mathbb{E}|\varepsilon_n|^{3}<\infty$.

In the present section we consider several signal series $x_n$
and the linear stationary series $e_n$ as the ``noise series''. (The series $e_n$  is
determined by the right-hand side of \eqref{eq:linear_process_f}.)
In this case
${\bf P}_0^{\perp}(\delta)$ is a random operator dependent on $\omega\in \Omega$.
Moreover, for fixed $(\delta,N,L)$ the
condition \eqref{eq:cond_main} of Theorem \ref{theor:full_decomp} holds only with a certain probability.

As in Section \ref{ssect:sign_incr_exp}, we further present   general statements and then consider several examples.
\begin{proposition}
\label{prop:gen_wnoise_LeqK}
1. Suppose that $L\rightarrow \infty$.
If $\mu_{\max}/\mu_{\min}^2=o\big(1/(N\ln{N})\big)$, then
there exist $\Omega_0\in {\cal F}$ with $\mathbb{P}(\Omega_0)=1$ and
an absolute constant $\gamma'$ such that for any $\omega \in \Omega_0$ and any $\delta$
\begin{gather*}
\limsup_{N}\, \frac{\mu_{\min}}{\sqrt{\mu_{\max}N\ln N}}\ \big\|{\bf P}_0^{\perp}(\delta)-{\bf P}_0^{\perp}\big\|
\!\leq\! 8C\sqrt{\gamma'S} |\delta|
\ \, \mbox{\rm and}
\\
\limsup_{N}\ \frac{\mu^2_{\min}}{{\mu_{\max}N\ln N}}\,\big\|{\bf P}_0^{\perp}(\delta)-{\bf P}_0^{\perp}-\delta {\bf V}_{0}^{(1)}\big\|
\leq C'\gamma'\,S\delta^2,
 \end{gather*}
where $C,C'$ are defined in Proposition \ref{prop:big_NSR} and $S\!=\!\sum_j |c_j|$.\\
2. If $L=L_0={\rm const}$ and  $\mu_{\max}/\mu_{\min}^2=o\big(1/N\big)$,
then
there exists a certain $\Omega'\in {\cal F}$ with $\mathbb{P}(\Omega')=1$ such that
for any $\omega \in \Omega'$ and for
any $\delta$
\begin{gather*}
\limsup_{N}\,\frac{\mu_{\min}}{\sqrt{N\mu_{\max}}}\ \big\|{\bf P}_0^{\perp}(\delta)-{\bf P}_0^{\perp}\big\|
\leq 8C\sqrt{\sigma_{\max}}\,|\delta| \ \ \mbox{\rm and}\\
\limsup_{N}\,\frac{\mu^2_{\min}}{N\mu_{\max}}\
\big\|{\bf P}_0^{\perp}(\delta)-{\bf P}_0^{\perp}-\delta {\bf V}_{0}^{(1)}\big\|
\leq C'\sigma_{\max}\,\delta^2,
\end{gather*}
where $\sigma_{\max}$ is the maximal eigenvalue of the matrix ${\bf \Sigma}=\big\{R_e(i-j)\big\}_{0\leq i,j<L_0}$ and
$R_e(\,\cdot\,)$ stands for the covariance function of $e_n$.
\end{proposition}
The following examples illustrate Proposition \ref{prop:gen_wnoise_LeqK}.
\begin{example} Let $N\rightarrow \infty$.
\label{ex:WN_LeqK}
\item[1.]
Consider the exponential signal series \eqref{eq:exp_sum} with  $\beta_k\neq 0$ and decreasing $a_k>1$.
Denote $\tau=a_1/a_p^2$ and
suppose that $\tau<1$.

a)  If $L\rightarrow \infty$, then $\mu_{\max}/\mu_{\min}^2\sim \Xi_a^2\,\tau^{2N}$ and
almost surely
\begin{gather*}
\begin{split}
\!&\!\!\limsup_{N}\, (N\ln N)^{-1/2} \tau^{-N}\big\|{\bf P}_0^{\perp}(\delta)-{\bf P}_0^{\perp}\big\|
\!\leq\! 8C\,\Xi_a \sqrt{\gamma_0S}\,|\delta|,\\
&\!\limsup_{N}\, (N\ln N)^{-1} \tau^{-2N}\!\big\|{\bf P}_0^{\perp}(\delta)\!-\!{\bf P}_0^{\perp}
\!-\!\delta {\bf V}_{0}^{(1)}\!\big\|
\!\leq\! C' \Xi_a^2 \gamma_0 S\delta^2
\end{split}
\end{gather*}
for any $\delta$.

b)
For the same  signal series \eqref{eq:exp_sum} and $L=L_0={\rm const}$
\begin{gather*}
\limsup_{N}\,N^{-1/2}\,\tau^{-N}\big\|{\bf P}_0^{\perp}(\delta)-{\bf P}_0^{\perp}\big\|
\leq 8C\sqrt{\lambda_{\max}}\,|\delta|,\\
\limsup_{N}\,N^{-1}\,\tau^{-2N}\big\|{\bf P}_0^{\perp}(\delta)-{\bf P}_0^{\perp}-\delta {\bf V}_{0}^{(1)}\big\|
\leq C'\, \lambda_{\max}\,\delta^2.
\end{gather*}
with probability 1 for any $\delta$.
\item[2.]
Consider the the polynomial signal series \eqref{eq:polyn} and suppose that $L/N\rightarrow \alpha\in(0,1)$. Then
Lemma \ref{summ:series} shows that
$\mu_{\max}/\mu_{\min}^2\sim \Psi_p^{2}\, N^{-2p-2}$ with $\Psi_p=
\sqrt{\Theta_{\max}}/\Theta_{\min}$.
Therefore, almost surely
\begin{gather*}
\begin{split}
\!&\limsup_{N}\ (\ln N)^{-1/2} N^{p+0.5}\big\|{\bf P}_0^{\perp}(\delta)-{\bf P}_0^{\perp}\big\|
\leq 8C \Psi_p \sqrt{\gamma_0S} |\delta|,\\
\!&\limsup_{N}\ (\ln N)^{-1} N^{2p+1}\big\|{\bf P}_0^{\perp}(\delta)\!-\!{\bf P}_0^{\perp}
\!-\!\delta {\bf V}_{0}^{(1)}\!\big\|
\!\leq\! C' \Psi_p^2 \gamma_0 S\delta^2
\end{split}
\end{gather*}
for any $\delta$.
\item[3.]
Consider the oscillating signal \eqref{eq:osc} and assume that $L/N\rightarrow \alpha\in (0,1)$.
Then
$\mu_{\max}/N^2\rightarrow \alpha(1-\alpha)\Lambda_{\max}$ and $\mu_{\min}/N^2\rightarrow \alpha(1-\alpha)\Lambda_{\min}$.
Therefore,
\begin{gather*}
\mu_{\max}/\mu_{\min}^2\sim N^{-2}\Upsilon_{\alpha}^2 =o\big(1/(N\ln{N})\big)
\end{gather*}
with $\Upsilon_{\alpha}=\sqrt{\big(\alpha(1-\alpha)\big)^{-1}\Lambda_{\max}}/\Lambda_{\min}$
and
\begin{gather*}
\limsup_{N}\, \sqrt{N (\ln N)^{-1}}\,\big\|{\bf P}_0^{\perp}(\delta)-{\bf P}_0^{\perp}\big\|
\leq 8C\,\Upsilon_{\alpha} \sqrt{\gamma_0S}\,|\delta|, \\
\limsup_{N}\, {N (\ln N)^{-1}} \big\|{\bf P}_0^{\perp}(\delta)-{\bf P}_0^{\perp}-\delta {\bf V}_{0}^{(1)}\big\|
\!\leq\! C'\Upsilon_\alpha\Psi \gamma_0 S\delta^2
\end{gather*}
with probability 1 for any $\delta$.
\end{example}

The case of  oscillating signal  series
 and white-noise ``noise series'' $e_n=\varepsilon_n$ is of  particular importance.

\begin{proposition}
\label{prop:gen_wight_noise_coarse}
Consider the  oscillating signal  series \eqref{eq:osc} and i.i.d. ``noise series'' $e_n=\varepsilon_n$ defined on a
probability space $(\Omega, {\cal F}, \mathbb{P})$. Assume that $\mathbb{E}\varepsilon_n=0$,
$\mathbb{D}\varepsilon_n=1$, and $\mathbb{E}|\varepsilon_n|^3<\infty$.

If $L=L_0={\rm const}$,  $\delta_0^2<L_0\Lambda_{\min}/4$,   and $N\rightarrow\infty$, then there exists
$\Omega'\in {\cal F}$ with $\mathbb{P}(\Omega')=1$ such that for any $\delta\in (-\delta_0,\delta_0)$ and any
$\omega\in \Omega'$
\begin{gather}
\label{eq:gen_wight_noise_coarse1}
\limsup_{N}\sqrt{N (\ln \ln N)^{-1}}\,\big\|{\bf P}_0^{\perp}(\delta)-{\bf P}_0^{\perp}\big\|<c'\, |\delta|
\end{gather}
and
\begin{gather}
\label{eq:gen_wight_noise_coarse2}
\limsup_{N}{N (\ln \ln N)^{-1}}\big\|{\bf P}_0^{\perp}(\delta)\!-\!{\bf P}_0^{\perp}\!-\!{\bf T}(\delta)\big\|
\!<\!c{''}\delta^2
\end{gather}
where positive constants $c'$ and $c''$  depend on $L_0$, $\delta_0$, and parameters of the series \eqref{eq:osc}.
\end{proposition}

\section{Examples. Precise asymptotic results}
\label{sect:hankel_precise}
The results of Section~\ref{sect:Hankel_pert}
can provide
overestimated upper bound for the rate of convergence
$\big\|{\bf P}_0^{\perp}(\delta)-{\bf P}_0^{\perp}\big\|\rightarrow 0$, see  Remark \ref{rem:reduced_rate}.
Morover,  operators ${\bf W}_1(\delta)$,  ${\bf L}(\delta)$, and ${\bf T}(\delta)$ determined by
\eqref{eq:first_term}, \eqref{eq:Ldelta_short}, and \eqref{eq:T_1} respectively can serve only as candidates for
the main term of this convergence.

In this section we present several examples of ``signal'' and ``noise'' series for which both true rates of convergence
$\big\|{\bf P}_0^{\perp}(\delta)-{\bf P}_0^{\perp}\big\|\rightarrow 0$ and the related
 main terms are obtained explicitly.

These examples are rather simple and  have minor practical value. In particular,
 we consider signal series of rank 1. Then all relevant operators are derived in accessible form.

In addition, the choice of examples corresponds to different
conditions for the convergence $\big\|{\bf P}_0^{\perp}(\delta)-{\bf P}_0^{\perp}\big\|\rightarrow 0$.
More precisely, the example of Section \ref{ssect:exp_const_prec} describes the case of ``big signal'' and
``small noise'' while the Section \ref{ssect:const_saw} is devoted to  asymptotic orthonormalities.
The case of white noise ``noise series'' is  considered in Section \ref{sssect:white_noise_noise}.

\subsection{Exponential signal series and constant ``noise series''}
\label{ssect:exp_const_prec}
Consider the ``signal series'' $x_n=a^n$ with $a>1$ and
the ``noise series'' $e_n\equiv 1$.
Let $W_j=\big(1,a,\ldots,a^{j-1}\big)^{\rm T}$, ${E}_j=(1,\ldots,1)^{\rm T}\in \mathbb{R}^j$, and $\beta_j=W_j^{\rm T}E_j$.
Then $\|E_j\|=\sqrt j$ and
 $\|W_j\|=\sqrt{(a^{2j}-1)/(a^2-1)}$.

Since matrices ${\bf H}=W_LW_K^{\rm T}$ and ${\bf E}=E_LE_K^{\rm T}$ have rank 1, then
$\nu_{\max}=\|E_L\|^2\|E_K\|^2=LK$ and
$\mu_{\max}=\mu_{\min}=\|W_L\|^2\|W_K\|^2\asymp a^{2N}$ as $N\rightarrow\infty$.
Thus, applying Proposition \ref{prop:big_NSR} with
$\Theta \sim \sqrt{LK}/a^{N}\rightarrow 0$ and using
\eqref{eq:Diff_proj} and \eqref{eq:Diff_proj_second}, we obtain
that
for any $\delta$ and $N>N_0(\delta)$
\begin{gather}
\label{eq:sec_or_exp1}
\big\|{\bf P}_0^{\perp}(\delta)-{\bf P}_0^{\perp}\big\|=|\delta|O\big(\sqrt{LK}a^{-N}\big)
\end{gather}
and
\begin{gather}
\label{eq:sec_or_exp}
\big\|{\bf P}_0^{\perp}(\delta)-{\bf P}_0^{\perp}-\delta {\bf V}_{0}^{(1)}\big\|=\delta^2 O\big(LKa^{-2N}\big),
\end{gather}
where ${\bf V}_{0}^{(1)}$ is defined in \eqref{eq:lin_proj}.

The following proposition presents precise versions of \eqref{eq:sec_or_exp1} and \eqref{eq:sec_or_exp}.
Denote
\begin{gather*}
H(a,L)=\frac{(a+1)}{a}\,\frac{a^L\sqrt{L}\,\|W_L\|^2-\beta_L^2}{\|W_L\|^2}\ .
\end{gather*}

\begin{proposition}
\label{prop:exp_const_proj_res}
Let ${\bf Z}_0^{(1)}$ stand for ${\bf E}{\bf H}^{\rm T}{\bf S}_0+{\bf S}_0{\bf H}{\bf E}^{\rm T}$.\\
1. If $L/N\rightarrow \alpha\in (0,1)$, then
\begin{gather*}
\frac{a^N}{\sqrt N}\,\big\|{\bf P}_0^{\perp}(\delta)-{\bf P}_0^{\perp}\big\|\rightarrow
|\delta|\,\frac{\, \alpha\,(a+1)\sqrt{a^2-1}}{a}\ \ \ \mbox{\rm and}\\
a^N\big\|{\bf P}_0^{\perp}(\delta)-{\bf P}_0^{\perp}-
\delta\, {\bf Z}_0^{(1)}\big\|\rightarrow |\delta|\,\frac{2\,(a+1)^2}{a}\ .
\end{gather*}
2. If $K=K_0={\rm const}$ and $N\rightarrow \infty$, then
\begin{gather*}
\frac{a^N}{\sqrt N}\ \big\|{\bf P}_0^{\perp}(\delta)-{\bf P}_0^{\perp}\big\|
\rightarrow |\delta|\,\frac{\,(a+1)\sqrt{a^2-1}}{a(1-a^{-K_0})}\ \ \ \mbox{\rm and}\\
a^N\big\|{\bf P}_0^{\perp}(\delta)-{\bf P}_0^{\perp}-
\delta\, {\bf Z}_0^{(1)}\big\|\rightarrow |\delta|\,\frac{2\,(a+1)^2}{a(1-a^{-K_0})}\ .
\end{gather*}
3. If $L=L_0={\rm const}$ and $N\rightarrow \infty$, then
\begin{gather*}
\begin{split}
&a^N\,\big\|{\bf P}_0^{\perp}(\delta)-{\bf P}_0^{\perp}\big\|\rightarrow |\delta| H(a,L_0) \ \ \mbox{\rm while} \\
&\big\|{\bf P}_0^{\perp}(\delta)-{\bf P}_0^{\perp}-\delta {\bf V}_{0}^{(1)}\big\|=\delta^2O\big(Na^{-2N}\big).
\end{split}
\end{gather*}
\end{proposition}

\begin{remark}
\label{rem:exp_const_prec}
1. Matrices ${\bf V}_0^{(1)}$ and ${\bf Z}_0^{(1)}$ have the following explicit form:
${\bf V}_0^{(1)}={\bf Z}_0^{(1)}-{\bf Z}_0^{(2)}$, where
\begin{gather*}
\begin{split}
&{\bf Z}_0^{(1)}=
\beta_K\frac{E_LW_L^{\rm T}+W_LE_L^{\rm T}}{\|W_L\|^2\|W_K\|^2}\quad \mbox{and } \\
&{\bf Z}^{(2)}_0=
\beta_L\beta_K\frac{W_LW_L^{\rm T}}{\|W_L\|^4\|W_K\|^2}\,.
\end{split}
\end{gather*}
2. Proposition \ref{prop:exp_const_proj_res} shows that the best rate of convergence of
${\bf P}_0^{\perp}(\delta)$ to ${\bf P}_0^{\perp}$ is achieved under the choice $L=L_0={\rm const}$. This rate is equal to
$a^{-N}$ and the main term of approximation ${\bf P}_0^{\perp}(\delta)\approx{\bf P}_0^{\perp}$ is the linear term
$\delta {\bf V}_0^{(1)}$. In this case the norms of operators ${\bf Z}_0^{(1)}$ and ${\bf Z}_0^{(2)}$
have the same order of growth
and therefore we cannot reduce ${\bf V}_0^{(1)}$ to ${\bf Z}_0^{(1)}$
as it was  done in the case $L\rightarrow \infty$.
\end{remark}

\subsection{Constant series as a signal}
In this section we consider the constant series $x_n=1$ as a signal.
Let
$W_j=(1,\ldots,1)^{\rm T}\in \mathbb{R}^j$. Then the $L\times K$ trajectory matrix of the series has the form
${\bf H}=W_LW_K^{\rm T}$ and $\big\|{\bf H}\big\|=\sqrt{LK}$. Therefore, $d=1$ and $\mu_{\max}=\mu_{\min}=LK$.

\subsubsection{Saw series as a ``noise''}
\label{ssect:const_saw}
Consider the saw series $e_n=(-1)^n$ as a ``noise''.
Then ${\bf E}=E_LE_K^{\rm T}$, where
$E_j\!=\!(e_0,\ldots,e_{j-1})^{\rm T}$.
As earlier, let $\beta_j\!=\!W_j^{\rm T}E_j$. Note that $\beta_j\!=\!0$ for even $j$ and $\beta_j\!=\!1$ for odd $j$.

Since $\nu_{\max}=LK$, then
$\Theta_1=\Theta_2=1$ and we cannot apply  Proposition \ref{prop:big_NSR}.
Nevertheless we can use Theorem \ref{theor:precise_proj} and Proposition \ref{prop:asymp_ort} to obtain conditions for
the convergence
$\big\|{\bf P}_0^{\perp}(\delta)-{\bf P}_0^{\perp}\big\|\rightarrow 0$ for small $\delta$ and $N\rightarrow \infty$.
Since $d=1$ we impose $\min(L,K)>1$.

\begin{proposition}
\label{prop:ort_const_saw}
Let $|\delta|<1/2$.
If $L$ and $K$ are both even, then ${\bf P}_0^{\perp}(\delta)={\bf P}_0^{\perp}$, otherwise
\begin{gather*}
\begin{split}
&\big\|{\bf P}_0^{\perp}(\delta)-{\bf P}_0^{\perp}\big\|
\!=\! |\delta|
\begin{cases}
O\big(K^{-1}\big)\!&
\begin{aligned}
&\mbox{\rm for odd} \,  K, \mbox{\rm even} \,  L
\\
&\mbox{\rm and } \,  K \rightarrow \infty,
\end{aligned}\\
O\left(L^{-1}\right)\!&
\begin{aligned}
&\mbox{\rm for odd} \, L , \mbox{\rm even} \,  K
\\
&\mbox{\rm and} \, L\rightarrow \infty,
\end{aligned}\\
O\big(L^{-1}\!+\!K^{-1}\big)\!&
\begin{aligned}
&\mbox{\rm for} \,  L,K\, \mbox{\rm both odd}\\
&\mbox{\rm and} \,   \min(L,K)\!\rightarrow\! \infty.
\end{aligned}
\end{cases}
\end{split}
\end{gather*}
%
\end{proposition}

\begin{remark}
Since the series $E$ is periodic with period 2, the different results for odd and even $L,K$ are not amazing.  In view of
 \cite[ch. 1 \S 1.6.2]{GNZh01}, the choice of window length as an integer multiple of   the period of periodical
component of the series  can essentially improve separability. Here we have got the formal affirmation of
this principle.

The case for even $L,K$ exactly corresponds to biorthogonality of matrices ${\bf H}$ and ${\bf E}$, see
the discussion in Section \ref{ssect:zero_pert}.
Therefore, no additional restrictions on $L,K$ are needed.
\end{remark}

Now our goal is to study the main term of the difference $\!{\bf P}_0^{\perp}(\delta)\!-\!{\bf P}_0^{\perp}\!$ by applying
Theorem\! \ref{theor:precise_proj} and Proposition~\ref{prop:L_delta}.

Consider fixed $\delta$ such that  $|\delta|<1/2$. Denote
\begin{gather}
\begin{split}
& E_{L,w}=E_L-\beta_L\,W_L/L,
 \label{eq:ELw}
 \end{split}
\end{gather}
$W_E=W_LE_{L,w}^{\rm T}+E_{L,w}W_L^{\rm T}$,
and
\begin{gather}
\begin{split}
{\bf M}(\delta)
\!=\!
\begin{cases}
\displaystyle
\frac{\delta}{1-\delta^2}\,\frac{W_E}{LK}&
\begin{aligned}
&\mbox{for odd} \,  K \\
&\mbox{and even} \,  L,
\end{aligned}\\
\displaystyle
\frac{\delta^2}{1-\delta^2}\,\frac{W_E}{L\sqrt{L^2-1}}&
\begin{aligned}
&\mbox{for odd} \,  L \\
&\mbox{and even} \,  K,
\end{aligned}\\
\displaystyle
\frac{\delta}{1-\delta^2}\!\left(\!\frac{1}{K}\!+\!\frac{\delta}{L}\!\right)\!
\frac{W_E}{\sqrt{L^2\!-\!1}}&
\begin{aligned}
&\mbox{for} \,  L,K \\
&\mbox{both odd}.
\end{aligned}
\end{cases}
\label{eq:M(delta)}
\end{split}
\end{gather}

\begin{proposition}
\label{prop_const_saw}
The norm of the operator ${\bf M}(\delta)$ is
\begin{gather}
\big\|{\bf M}(\delta)\big\|\!=\!
\begin{cases}
\displaystyle
\frac{|\delta|}{1-\delta^2}\frac{1}{K}&
\! \mbox{\rm for odd} \,  K \,  \mbox{\rm and even} \,  L,\\
\displaystyle
\frac{\delta^2}{1-\delta^2} \frac{1}{L}&
\! \mbox{\rm for odd} \,  L \,  \mbox{\rm and even} \,  K,\\
\displaystyle
\frac{|\delta|}{1-\delta^2}\left|\frac{1}{K}\!+\!\frac{\delta}{L}\right|&
\! \mbox{\rm for} \,  K \,  \mbox{\rm and}\  L \,  \mbox{\rm both  odd},
\end{cases}
\label{eq:norm_M_const_saw}
\end{gather}
while
\begin{gather}
\begin{split}
&\big\|{\bf P}_0^{\perp}(\delta)-{\bf P}_0^{\perp}-{\bf M}(\delta)\big\|\\
&\!=\!
\begin{cases}
O\big(K^{-2}\big)\!&
\begin{aligned}
&\mbox{\rm for odd} \,  K, \mbox{\rm even} \,  L\\
&\mbox{\rm and } \,  K \rightarrow \infty,
\end{aligned}\\
O\left(L^{-2}\right)\!&
\begin{aligned}
&\mbox{\rm for odd} \, L , \mbox{\rm even} \,  K\\
&\mbox{\rm and } \,  L\rightarrow \infty,
\end{aligned}\\
O\big(L^{-2}\!+\!K^{-2}\big)\!&
\begin{aligned}
&\mbox{\rm for} \,  L,K\, \mbox{\rm both odd}\\
&\mbox{\rm and} \,   \min(L,K)\!\rightarrow\! \infty.
\end{aligned}
\end{cases}
\label{eq:res_prec_const_saw}
\end{split}
\end{gather}
\end{proposition}

\begin{remark}
1. Unlike to Proposition \ref{prop:exp_const_proj_res},
Proposition \ref{prop_const_saw} shows that the main term ${\bf M}(\delta)$ is
nonlinear in~$\delta$.\\
2. The norm \eqref{eq:norm_M_const_saw} indicates the distinctions between positive and negative $\delta$ in the
case of odd $L,K$. For example, if $\delta<0$ and  $L=-\delta K$, then
$\big\|{\bf P}_0^{\perp}(\delta)-{\bf P}_0^{\perp}\big\|\asymp L^{-2}$ whereas this norm has the order $L^{-1}$ for any positive $\delta$ and $L\sim K$.
\end{remark}

\subsubsection{White noise as the ``noise series''}
\label{sssect:white_noise_noise}
In this section we study the difference ${\bf P}_0^{\perp}(\delta)-{\bf P}_0^{\perp}$ for the constant signal $x_n=1$
and the white-noise ``noise series'' $e_n=\varepsilon_n$ where $e_n$ are i.i.d. random variables $\varepsilon_n$ such that
$\mathbb{E}\varepsilon_n=0$, $\mathbb{E}\varepsilon_n^2=1$ and
 $\mathbb{E}\varepsilon_n^4<\infty$.

The a.s. behavior of $\big\|{\bf P}_0^{\perp}(\delta)-{\bf P}_0^{\perp}\big\|$ as $N\rightarrow \infty$ is already
studied in Section \ref{ssect:hankel_coarse}.
Now we derive the main term of the difference ${\bf P}_0^{\perp}(\delta)-{\bf P}_0^{\perp}$ as
$N\rightarrow \infty$ and $L=L_0={\rm const}$ in terms of a suitable version of the central limit theorem.

Denote $\Psi_{L_0}=\{\psi_{ij}\}_{i,j=0}^{L_0-1}$
a random symmetrical $L_0\times L_0$ Toeplitz  matrix where $\psi_{0j}$ ($j=0,\ldots,L_0-1$)
are independent,
$\psi_{00}\in \mathrm{N}(0,\mathbb{E}\varepsilon^4-1)$, and $\psi_{0j}\in \mathrm{N}(0,1)$ for $j\geq 1$.
(Note that $\xi\in \mathrm{N}(a,\sigma^2)$ means that $\xi$ has a normal distribution, $\mathbb{E}\xi=a$, and $\mathbb{D}\xi=\sigma^2$.)
\begin{proposition}
\label{prop:const_wnoise_main_term_L_0}
For fixed $L_0,\delta$ such that $\delta^2\!<\!L_0/4$ denote
\begin{gather*}
\Omega_N=\{\omega\in \Omega \  \mbox{\rm such   that} \ \|{\bf B}(\delta)\|<\mu_{\min}/2\}.
\end{gather*}
Then $\mathbb{P}(\Omega_N)\rightarrow 1$ as $N\rightarrow \infty$  and
\begin{gather}
\begin{split}
&{\cal L}\left(\sqrt{N}\left({\bf P}_0^{\perp}(\delta)-{\bf P}_0^{\perp}\right)\,\big|\,\Omega_N\right)
\\
&
\Longrightarrow
{\cal L}\left(\frac{\delta^2}{L_0}\big({\bf P}_0^\perp \Psi_{L_0}\,{\bf P}_0
+{\bf P}_0 \Psi_{L_0}\,{\bf P}_{0}^\perp\big)\right),
\end{split}
\label{eq:const_witenoise_CLT}
\end{gather}
where $``\Rightarrow$'' stands for weak convergence of distributions and ${\cal L}(\xi)$ means the distribution of the
random vector $\xi$.
\end{proposition}

\begin{remark}
1. Roughly speaking, Proposition \ref{prop:const_wnoise_main_term_L_0} shows that the  main term of
the difference ${\bf P}_0^{\perp}(\delta)-{\bf P}_0^{\perp}$ has the  asymptotic form
\begin{gather*}
\frac{\delta^2}{L_0\sqrt{N}}\,\big({\bf P}_0^\perp \Psi_{L_0}\,{\bf P}_0+{\bf P}_0 \Psi_{L_0}\,{\bf P}_0^\perp\big).
\end{gather*}
2. Since $L_0$ is fixed and the norm $\|\,\cdot\,\|$ is the continuous functional on the space of
$L_0\times L_0$ matrices, then
\begin{gather*}
\begin{split}
&{\cal L}\left(\sqrt{N}\,\big\|{\bf P}_0^{\perp}(\delta)-{\bf P}_0^{\perp}\big\|\,\big|\,\Omega_N\right)
\\
&
\Longrightarrow
{\cal L}\left(\frac{\delta^2}{L_0}\, \big\|{\bf P}_0^\perp \Psi_{L_0}\,{\bf P}_0+
{\bf P}_0 \Psi_{L_0}\,{\bf P}_0^\perp\big\|\right)
\end{split}
\end{gather*}
under conditions of Proposition \ref{prop:const_wnoise_main_term_L_0}.
\end{remark}

\section{On the way to applications}
\label{sect:appl}
In the present section we briefly describe the application of the previous results to several methods of Signal Subspace
Analysis. Section \ref{ssect:LRFs} is devoted to the approximation of linear recurrent formulas governing the signal
and in Section \ref{ssect:LS_ESPRIT} we study the real-valued variant of Least-Square ESPRIT. We demonstrate that both methods
asymptotically  converge under assumption that $\|{\bf P}_0^\perp(\delta)-{\bf P}_0^\perp\|\rightarrow 0$.

In Section \ref{ssect:LS_ESPRIT} we discuss the reconstruction stage of Singular Spectrum Analysis (briefly, SSA).
Though the
precision of SSA can not be described only in terms of proximity of the perturbed and unperturbed projectors,
the obtained results help to reformulate the problem in a more transparent form.

\subsection{Approximations of linear recurrent formulas}
\label{ssect:LRFs}
Let the signal series $\mathrm{F}=(x_0,\ldots,x_n,\ldots)$ be governed by a linear recurrent formula (LRF)
\begin{gather}
\label{eq:LRF_d}
x_n=\sum_{k=1}^{d}b_kx_{n-k}, \quad n\geq d
\end{gather}
and suppose that \eqref{eq:LRF_d} is the minimal LRF for the series $\mathrm{F}$. In particular, this means that $b_d\neq 0$.

For $L,K>d$ let $\bf H$ stand for the trajectory $L\times K$ matrix  of the series $\mathrm{F}$. Then
$\mathrm{rank}\, {\bf H}=d$.  Consider the signal subspace $\mathbb{U}_0^\perp$ and the corresponding projector
${\bf P}_0^\perp$.
Denote $\mathfrak{e}_L=(0,0,\ldots,0,1)^{\rm T}\in \mathbb{R}^L$.
As it is proved in \cite[sect. 5.2]{GNZh01}, ${\bf P}_0 \mathfrak{e}_L\neq 0$.
(Note that $\|{\bf P}_0\mathfrak{e}_L\|$ is a cosine between vector $\mathfrak{e}_L$ and  linear space
$\mathbb{U}_0$.)

Let ${\bf G}_L$ stand for the
$(L-1)\times L$ matrix
\begin{gather}
{\bf G}_L=
\begin{pmatrix}
1&0&\ldots&0&0&0\\
0&1&\ldots&0&0&0\\
\vdots&\vdots&\ddots&\vdots&\vdots&\vdots\\
0&0&\ldots&1&0&0\\
0&0&\ldots&0&1&0
\end{pmatrix}
\label{eq:G_L}
\end{gather}
and denote
\begin{gather}
\label{eq:LRF_coef_L}
R= (a_{L-1},\ldots,a_1)^{\rm T}=-\, \frac{1}{\|{\bf P}_0\mathfrak{e}_L\|^2}\ {\bf G}_L{\bf P}_0 \mathfrak{e}_L.
\end{gather}
Due to \cite[th. 5.2]{GNZh01},
\begin{gather}
\label{eq:LRF_L}
x_n=\sum_{k=1}^{L-1}a_kx_{n-k}, \quad n\geq L.
\end{gather}
Note that analogous expressions are known from early 80s, see for example,
\cite{KumaresanT83}.
It is clear that \eqref{eq:LRF_L} coincides with \eqref{eq:LRF_d} for $L=d+1$.

The formula \eqref{eq:LRF_coef_L} can be used to derive the approximation of $R$ in the case when the signal series $\mathrm{F}$
is corrupted by an additive noise series $\mathrm{E}$. In other words, if $\mathrm{F}(\delta)=\mathrm{F}+\delta\mathrm{E}$, then
the natural form of this approximation has the form
\begin{gather*}
R(\delta)=-\, \frac{1}{\|{\bf P}_0(\delta) \mathfrak{e}_L\|^2}\ {\bf G}_L{\bf P}_0(\delta) \mathfrak{e}_L.
\end{gather*}
If ${\bf P}_0(\delta)$ is close to ${\bf P}_0$, then $R(\delta)$ must be close to $R$.

\begin{proposition}
\label{prop:close_R}
Let $\Delta {\bf P}(\delta)$ stand for $\|{\bf P}_0^\perp(\delta)-{\bf P}_0^\perp\|$. Denote
$\vartheta=\|{\bf P}_0^{\perp}\mathfrak{e}_L\|$ and
suppose that $\Delta {\bf P}(\delta)<\|{\bf P}_0 \mathfrak{e}_L\|=\sqrt{1-\vartheta^2}$.
Then
\begin{gather}
\label{eq:rec_DeltaR}
\|R(\delta)\!-\!R\|\!\leq\!
\frac{\Delta {\bf P}(\delta)}{1\!-\!\vartheta^2}\!
\left(\!1\!-\! \frac{\Delta {\bf P}(\delta)}{\sqrt{1\!-\!\vartheta^2}}\!\right)^{-2}\!
\left(\!1+\frac{2}{\sqrt{1\!-\!\vartheta^2}}\!\right).
\end{gather}
\end{proposition}

Proposition \ref{prop:close_R}
gives sufficient conditions for the convergence $\|R(\delta)-R\|\rightarrow 0$ in terms of
$\|{\bf P}_0^\perp(\delta)-{\bf P}_0^\perp\|$ and
$\|{\bf P}_0\mathfrak{e}_L\|$. Namely, if $\|{\bf P}_0^\perp(\delta)-{\bf P}_0^\perp\|\rightarrow 0$ and if
$\|{\bf P}_0\mathfrak{e}_L\|$ is separated from zero, then $\|R(\delta)-R\|\rightarrow 0$.

The second condition automatically holds in the case $L={\rm const}>d$ (then both $\mathfrak{e}_L$ and ${\bf P}_0$
do not depend on $N$).
It is not difficult to show (we omit proofs for brevity)
that $\|{\bf P}_0\mathfrak{e}_L\|$ is
separated from zero as $L\rightarrow \infty$ for signals of exponential type,
polynomial signals and oscillating signals. Therefore, in  these cases \eqref{eq:rec_DeltaR}
takes the form $\|R(\delta)-R\|\!=\!
O\big(\|{\bf P}_0^\perp(\delta)\!-\!{\bf P}_0^\perp\|\big)$.
Examples of
Section \ref{ssect:hankel_coarse}  refine this general assertion for a number of ``signal'' and
``noise'' series.

\subsection{LS-ESPRIT for real-valued signals}
\label{ssect:LS_ESPRIT}
As in Section \ref{ssect:LRFs}, we suppose that the signal series
$\mathrm{F}=(x_0,\ldots,x_n,\ldots)$ is governed by the linear recurrent formula \eqref{eq:LRF_d} with $b_d\neq 0$.
It is well-known that the general solution of \eqref{eq:LRF_d} is expressed through the roots
of the characteristic polynomial $P_d(\lambda)=\lambda_d-\sum_{k=1}^d b_k\lambda^{d-k}$.

Assume that $L,K> d$  and consider the $L\times K$ trajectory matrix ${\bf H}$
of the series $\mathrm{F}$. As in previous sections,
let $\mathbb{U}_0^\perp$ stand for the linear space spanned by columns of the matrix $\bf H$. For a certain basis
$U_1,\ldots,U_d$ of $\mathbb{U}_0^\perp$ denote ${\bf U}=[U_1:\ldots:U_d]$. Lastly, let the matrix
${\bf G}_L$ be defined by \eqref{eq:G_L} and
 denote
the $(L-1)\times L$ matrix ${\bf G}^{(L)}$ by
\begin{gather*}
{\bf G}^{(L)}=
\begin{pmatrix}
0&1&0&\ldots&0&0\\
0&0&1&\ldots&0&0\\
\vdots&\vdots&\vdots&\ddots&\vdots&\vdots\\
0&0&0&\ldots&1&0\\
0&0&0&\ldots&0&1
\end{pmatrix}
.
\label{eq:G^L}
\end{gather*}

Note that
\begin{gather*}
\begin{split}
&{\bf F}_1\stackrel{\rm def}={\bf G}_L^{\rm T}{\bf G}_L=
\begin{pmatrix}
1&0&\ldots&0&0\\
0&1&\ldots&0&0\\
\vdots&\vdots&\ddots&\vdots&\vdots\\
0&0&\ldots&1&0\\
0&0&\ldots&0&0
\end{pmatrix}
\ \ \mbox{\rm and}\\
&{\bf F}_2\stackrel{\rm def}={\bf G}_L^{\rm T}{\bf G}^{(L)}=
\begin{pmatrix}
0&1&0&\ldots&0&0\\
0&0&1&\ldots&0&0\\
\vdots&\vdots&\vdots&\ddots&\vdots&\vdots\\
0&0&0&\ldots&0&1\\
0&0&0&\ldots&0&0
\end{pmatrix}
.
\end{split}
\end{gather*}
(Both matrices are of the size $L\times L$.) It is easy to see that $\|{\bf F}_1\|=\|{\bf F}_2\|=1$.

The method called Least Square (briefly, LS) ESPRIT  for the
analysis of the series $\mathrm{F}$ is based on the following facts
(see discussion in \cite{Golyandina_new} and references within):
for any basis $U_1,\ldots,U_d$ of the linear space $\mathbb{U}_0^\perp$
\begin{enumerate}
\item
the matrix
${\bf U}^{\rm T}{\bf F}_1{\bf U}$ is invertible;
\item
the set of eigenvalues of the matrix
\begin{gather}
\label{eq:D_theor}
{\bf D}=\Big({\bf U}^{\rm T}{\bf F}_1{\bf U}\Big)^{-1}{\bf U}^{\rm T}{\bf F}_2{\bf U}
\end{gather}
coincides with the set of roots of the polynomial $P_d(\lambda)$ subject to  multiplicities of roots and
eigenvalues.
\end{enumerate}
In practice, leading left singular vectors of the matrix ${\bf H}$ usually stand
for the basis $U_1,\ldots,U_d$ of linear space
$\mathbb{U}_0^\perp$. The similar choice gives rise to approximation of the matrix ${\bf D}$ in the case when the series $\mathrm{F}$
is corrupted by an additive ``noise series'' $\mathrm{E}$ multiplied by a formal perturbation parameter $\delta$.

Denote ${\bf E}$ the $L\times K$ trajectory  matrix of the series $\mathrm{E}$.
If $U_j(\delta)$ $(j=1,\ldots,d)$ are leading left singular vectors of the matrix ${\bf H}(\delta)={\bf H}+\delta {\bf E}$
and ${\bf U}(\delta)=[U_1(\delta):\ldots:U_d(\delta)]$, then we can use the matrix
\begin{gather}
\label{eq:D_approx}
{\bf D}(\delta)=
\Big({\bf U}^{\rm T}(\delta){\bf F}_1{\bf U}(\delta)\Big)^{-1}{\bf U}^{\rm T}(\delta){\bf F}_2{\bf U}(\delta)
\end{gather}
to approximate  ${\bf D}$.

The following assertion helps to express ${\bf D}(\delta)$ through the perturbed projector ${\bf P}^\perp_0(\delta)$.
\begin{lemma}
\label{lem:lin_indep}
Let $U_1,\ldots,U_d$ be a basis of a linear space $\mathbb{U}\subset \mathbb{R}^L$. Denote ${\bf P}$ the orthogonal
projector on $\mathbb{U}$.\\
1. If ${\bf Q}: \mathbb{R}^L\mapsto \mathbb{R}^L$ and
$\|{\bf Q}-{\bf P}\|<1$, then vectors ${\bf Q}U_1,\ldots,{\bf Q} U_d$ are linearly independent.\\
2. Consider a linear space $\mathbb{V}\subset \mathbb{R}^L$ of dimension $d$ and denote ${\bf Q}$ the orthogonal
projector on $\mathbb{V}$. If $\|{\bf Q}-{\bf P}\|<1$, then there exist linearly independent vectors
$V_1,\ldots,V_d\in \mathbb{V}$ such that $U_j={\bf P}V_j$.
\end{lemma}

Lemma \ref{lem:lin_indep} shows that
if
$\|{\bf P}_0^\perp(\delta)-{\bf P}_0^\perp\|<1$,
then
the matrix
\begin{gather}
\label{eq:D_approx_th}
\widehat{\bf D}(\delta)\!\!=\!\!
\Big(\!{\bf U}^{\rm T}{\bf P}_0^\perp(\delta){\bf F}_1{\bf P}_0^\perp(\delta){\bf U}\!\Big)^{\!-1}\!
{\bf U}^{\rm T}{\bf P}_0^\perp(\delta){\bf F}_2{\bf P}_0^\perp(\delta){\bf U},
\end{gather}
where ${\bf U}=[U_1,\ldots,U_d]$,
has the same eigenvalues as the matrix \eqref{eq:D_approx} for any choice of the basis $U_1,\ldots,U_d$
of the linear space $\mathbb{U}_0^\perp$. Moreover $\widehat{\bf D}(\delta)={\bf D}(\delta)$ under a certain choice
of $U_1,\ldots,U_d$.

Suppose now that $\|{\bf P}_0^\perp(\delta)\!-\! {\bf P}_0^\perp\|\!\rightarrow\! 0$ as
$N\!=\!K\!+\!L\!-\!1\!\rightarrow\!\infty$.
Then it is natural to suppose that $\|\widehat{\bf D}(\delta)-{\bf D}\|\rightarrow 0$ and therefore asymptotically
as $N\rightarrow \infty$
we get all roots of the polynomial $P_d(\lambda)$ at least in the case $L={\rm const}$.

Let us formalize these considerations. As earlier, consider a certain basis $U_1,\ldots,U_d$
of the linear space $\mathbb{U}_0^\perp$ and set ${\bf U}=[U_1,\ldots,U_d]$.

\begin{proposition}
\label{prop:Delta_P1}
If $\Delta {\bf P}(\delta)\stackrel{\rm def}=
\|{\bf P}_0^\perp(\delta)-{\bf P}_0^\perp\|<\upsilon/2$ with $\upsilon$ standing for $\|{\bf U}^{\rm T}{\bf F}_1{\bf U}\|/\|{\bf U}\|^2$,
then
\begin{gather}
\label{eq:Delta_P1}
\big\|\widehat{\bf D}(\delta)-{\bf D}\big\|\leq \frac{2 \Delta {\bf P}(\delta)}{\upsilon}\
\left(1+\frac{1}{1-2\Delta {\bf P}(\delta)/\upsilon}\right).
\end{gather}
\end{proposition}

\begin{remark}
Denote $\mathfrak{e}_L=(0,0,\ldots,0,1)^{\rm T}\in \mathbb{R}^L$ and $\vartheta=\|{\bf P}_0^{\perp}\mathfrak{e}_L\|$.
Then ${\bf F}_1={\bf I}-\mathfrak{e}_L\mathfrak{e}_L^{\rm T}$
and
\begin{gather*}
\begin{split}
&{\bf U}^{\rm T}{\bf F}_1{\bf U}={\bf U}^{\rm T}{\bf U}-{\bf U}^{\rm T}\mathfrak{e}_L\mathfrak{e}_L^{\rm T}{\bf U}\\
&={\bf U}^{\rm T}{\bf U}-{\bf U}^{\rm T}\big({\bf P}_0^{\perp}\mathfrak{e}_L\big)
\big(\mathfrak{e}_L^{\rm T}{\bf P}_0^{\perp}\big){\bf U}.
\end{split}
\end{gather*}
Therefore, $\|{\bf U}^{\rm T}{\bf F}_1{\bf U}\|\!\geq\! \|{\bf U}\|^2\big(1\!-\!\vartheta^2\big)$ and under restriction
$\Delta {\bf P}(\delta)\!<\!(1\!-\!\vartheta^2)\!/2$ the inequality
\eqref{eq:Delta_P1} can be transformed~to
\begin{gather}
\label{eq:]Delta_P2}
\big\|\widehat{\bf D}(\delta)\!-\!{\bf D}\big\|\!\leq\! \frac{2 \Delta {\bf P}(\delta)}{1\!-\!\vartheta^2}
\left(1+\frac{1}{1\!-\!2\Delta {\bf P}(\delta)/(1\!-\!\vartheta^2)}\right).
\end{gather}
Note that the upper bound \eqref{eq:]Delta_P2} does not depend on a basis of the linear space $\mathbb{U}_0^\perp$.
\end{remark}

The inequality \eqref{eq:]Delta_P2} shows that $\|\widehat{\bf D}(\delta)-{\bf D}\|\rightarrow 0$ under the same
conditions as for linear recurrent formulas of Section \ref{ssect:LRFs}. More precisely,
 $\|{\bf P}_0^\perp(\delta)-{\bf P}_0^\perp\|$
must tend to zero and $\|{\bf P}_0\mathfrak{e}_L\|^2=1-\vartheta^2$ must be  separated from zero.

Moreover,
in this case $\|\widehat{\bf D}(\delta)-{\bf D}\|=O\big(\|{\bf P}_0^\perp(\delta)-{\bf P}_0^\perp\|\big)$.
Therefore, examples of Section \ref{ssect:hankel_coarse} provide the corresponding upper bounds for
$\|\widehat{\bf D}(\delta)-{\bf D}\|$.

\subsection{On the reconstruction stage of Singular Spectrum Analysis}
\label{ssect:Hankel}
If the aim of Singular Spectrum Analysis is treated as an interpretable decomposition of time
series $\mathrm{G}_N=(g_0,\ldots,g_{N-1})$ onto 2 or more additive components, then the whole SSA procedure can be expressed
as follows. (See \cite[ch. 1]{GNZh01} for details; for our goals it is sufficient
to decompose $\mathrm{G}_N$ onto 2 components.)

The decomposition stage consists of the choice of the ``window length'' $L$, construction of $L\times K$ ``trajectory matrix''
${\bf G}_N$ of the series $\mathrm{G}_N$  and Singular Value Decomposition of ${\bf G}$
onto ``elementary'' rank-one
matrices ${\bf G}_N^{(j)}$.

The reconstruction stage consists of the summation of a certain number of ${\bf G}_N^{(j)}$ (then we get the ``reconstructed''
matrix ${\bf G}'_N$) and ``hankelization'' of ${\bf G}'_N$. Formally, the result $\mathcal{S}{\bf G}'_N$ of hankelization of
the matrix ${\bf G}'_N=\{g'_{ij}\}_{{i=0},\,{j=0}}^{{L-1,\,K-1}}$ is the $L\times K$ Hankel matrix with elements
$\widetilde{g}\,'_{ij}$ equal to the average of $g'_{kl}$ such that $k+l=i+j$.
Since each $L\times K$ Hankel matrix is in natural one-to-one correspondence with a series of length $N=L+K-1$, then
we obtain the decomposition $\mathrm{G}_N=\widetilde{\mathrm{G}}'_N+(\mathrm{G}_N-\widetilde{\mathrm{G}}'_N)$
of the initial series $\mathrm{G}_N$.

Suppose now that the series $\mathrm{G}_N$ is the sum of the ``signal'' $\mathrm{F}_N$ governed by LRF \eqref{eq:LRF_d}
with $b_d\neq 0$ and the series $\delta \mathrm{E}_N$, where $\mathrm{E}_N$ is a ``noise series'' and $\delta$ is a formal
perturbation parameter. Then it is  natural to state the problem of an (approximate) extraction of the signal $\mathrm{F}_N$ from the sum
$\mathrm{G}_N=\mathrm{F}_N+\delta \mathrm{E}_N$.

For small $\delta$ the standard approach to this problem is expressed in terms discussed in previous sections. Namely,
under the choice of
$d<L<N-d-1$ the trajectory matrix ${\bf G}_N={\bf H}(\delta)={\bf H}+\delta {\bf E}$ is processed by SVD and $d$
leading elementary matrices
are summed to get the approximation $\widetilde{\bf H}={\bf G}'_N$ of ${\bf H}$. Then the hankelization
procedure yields the approximation $\widetilde{\mathrm{F}}_N(\delta)$ of the series $\mathrm{F}_N$.

Note that ${\bf H}={\bf P}_0^{\perp}{\bf H}$. Then
$\widetilde{\bf H}={\bf P}_0^{\perp}(\delta)\big({\bf H}+\delta{\bf E}\big)$ and therefore the approximation
error $\Delta_\delta({\bf H})=\widetilde{\bf H}-{\bf H}$
of the reconstructed trajectory matrix  has the form
\begin{gather}
\begin{split}
&\Delta_\delta({\bf H})={\bf P}_0^{\perp}(\delta){\bf H}(\delta)-{\bf P}_0^{\perp}{\bf H}\\
&=\left({\bf P}_0^{\perp}(\delta)-{\bf P}_0^{\perp}\right){\bf H}(\delta)+
\delta\,{\bf P}_0^{\perp}{\bf E}.
\end{split}
\label{eq:rec_matr_error}
\end{gather}
\begin{remark}
\label{rem:recon_zero}
It follows from results of Section \ref{ssect:zero_pert}
that $\widetilde{\bf H}={\bf H}$ provided that ${\bf H}$ and ${\bf E}$
are biorthogonal.
\end{remark}

To measure the difference between $\widetilde{\mathrm{F}}_N(\delta)$ and $\mathrm{F}_N$ we must introduce a
convenient metric. In many practical cases the proper choice is
\begin{gather*}
\big\|\widetilde{\mathrm{F}}_N(\delta)-\mathrm{F}_N\big\|_{\max}=\max_{\,0\leq n<N}|\widetilde{f}_n(\delta)-f_n|.
\end{gather*}

Assume now that for any $N$  both $\mathrm{F}_N$ and $\mathrm{E}_N$ are segments of infinite series $\mathrm{F}$ and
$\mathrm{E}$.
If  $\|\widetilde{\mathrm{F}}_N(\delta)-\mathrm{F}_N\|_{\max}$ tends to zero as $L=L(N)$ and $N\rightarrow \infty$, then
SSA {\it asymptotically reconstructs} the (infinite) series $\mathrm{F}$ from the (infinite)
perturbed series $\mathrm{F}+\delta \mathrm{E}$.
It is easy to see that
\begin{gather}
\begin{split}
&\big\|\widetilde{\mathrm{F}}_N(\delta)-\mathrm{F}_N\big\|_{\max}=
\big\|\mathcal{S}\big({\bf P}_0^{\perp}(\delta){\bf H}(\delta)\big)-{\bf H}\big\|_{\max}\\
&=\big\|\mathcal{S}\Delta_\delta({\bf H})\big\|_{\max}\,,
\label{eq:hankel_matrix_series}
\end{split}
\end{gather}
where $\|{\bf A}\|_{\max}=\max_{i,j}|a_{ij}|$ for the matrix
${\bf A}$ with entries~$a_{ij}$ and
$\mathcal{S}$ stands for the hankelization operator. It is clear that $\|\mathcal{S}{\bf A}\|_{\max}\leq \|{\bf A}\|_{\max}$.

The spectral norm $\|{\bf A}\|$ and uniform norm $\|{\bf A}\|_{\max}$ are equivalent, but this
equivalence is lost as the size of a matrix tends to infinity.

It is well-known that
$\|{\bf A}\|_{\max}\leq \|{\bf A}\|$. This means that if the spectral norm of the matrix ${\bf A}$ is small, then
all entries of the Hankel matrix $\mathcal{S}\mathbf{A}$ are small too.
The opposite inequality, which is also well-known (see \cite{GolubvanL96} for both inequalities)  has the form
$\|{\bf A}\|\leq \sqrt{LK}\|{\bf A}\|_{\max}$ for an $L\times K$ matrix ${\bf A}$. This inequality gives a hint that
a large-size matrix with small entries can have a large spectral norm. Indeed, the $n\!\times \!n$ Hankel matrix
${\bf G}_n$ with equal entries $g^{(n)}_{ij}\!=\!n^{-1/2}$ has the spectral norm $\|{\bf G}_n\|\!=\!\sqrt{n}\rightarrow \!\infty$
while  $g^{(n)}_{ij}\!\rightarrow\! 0$ as~$n\!\rightarrow\! \infty$.

In general, this means that even if the spectral norm $\|\Delta_\delta({\bf H})\|$ does not tend to infinity as
$N\rightarrow \infty$, still the convergence
$\|\widetilde{\mathrm{F}}_N(\delta)-\mathrm{F}_N\|_{\max}\rightarrow 0$ can occur.
 Therefore, taking into account the last term of the right-hand side of
\eqref{eq:rec_matr_error}, we see that general upper bounds
of the kind \eqref{eq:rec_DeltaR} or \eqref{eq:Delta_P1} can hardy be valid for
$\|\widetilde{\mathrm{F}}_N(\delta)-\mathrm{F}_N\|_{\max}$.

However, equalities \eqref{eq:rec_matr_error} and \eqref{eq:hankel_matrix_series} help to simplify the problem.
Suppose that
\begin{gather}
\label{eq:cond_inf}
\left\|\left({\bf P}_0^{\perp}(\delta)-{\bf P}_0^{\perp}-{\bf N}(\delta)\right){\bf H}(\delta)\right\|
\rightarrow 0
\end{gather}
as $N\rightarrow \infty$
for certain operators ${\bf N}(\delta)={\bf N}_N(\delta)$. Then
\begin{gather*}
\begin{split}
&\mathcal{S}\Delta_\delta({\bf H})=
\mathcal{S}\!\left(\big({\bf P}_0^{\perp}(\delta)-{\bf P}_0^{\perp}-{\bf N}(\delta)\big){\bf H}(\delta)\right)\\
&+
\mathcal{S}\!\left({\bf N}(\delta){\bf H}(\delta)+\delta {\bf P}_0^{\perp}{\bf E}\right)
\end{split}
\end{gather*}
and
\begin{gather*}
 \big\|\widetilde{\mathrm{F}}_N(\delta)-\mathrm{F}_N\big\|_{\max}=
 \big\|\mathcal{S}\big({\bf N}(\delta){\bf H}(\delta)+\delta {\bf P}_0^{\perp}{\bf E}\big)\big\|_{\max}
 +o(1)
\end{gather*}
as $N\rightarrow \infty$.
Thus, the problem reduces to the investigation of the asymptotic behavior of entries  of Hankel matrices
$\mathcal{S}\big({\bf N}(\delta){\bf H}(\delta)+\delta {\bf P}_0^{\perp}{\bf E}\big)$.

Let us consider several examples of signal and  noise series when the convergence
\eqref{eq:cond_inf} occurs. In all examples we straightforwardly apply inequalities of Section \ref{ssect:hankel_coarse},
therefore results are given without comments.
The norm
$\left\|\left({\bf P}_0^{\perp}(\delta)-{\bf P}_0^{\perp}-{\bf N}(\delta)\right){\bf H}(\delta)\right\|$ is estimated~as
\begin{gather*}
\begin{split}
&\left\|\left({\bf P}_0^{\perp}(\delta)-{\bf P}_0^{\perp}-{\bf N}(\delta)\right){\bf H}(\delta)\right\|\\
&\leq
\big\|{\bf P}_0^{\perp}(\delta)-{\bf P}_0^{\perp}-{\bf N}(\delta)\big\|\,\left(\big\|{\bf H}\big\|+
|\delta|\, \big\|{\bf E}\big\|\right).
\end{split}
\end{gather*}
In the  most of examples ${\bf N}(\delta)=\delta{\bf V}^{(1)}_0$.
Due to \eqref{eq:lin_proj}, we~have
\begin{gather*}
{\bf V}^{(1)}_0 {\bf H}(\delta)={\bf P}_{0}{\bf E}{\bf H}^{\rm T}{\bf S}_0{\bf H}
+\delta\big({\bf P}_{0}{\bf E}{\bf H}^{\rm T}{\bf S}_0{\bf E}+{\bf S}_0{\bf H}{\bf E}^{\rm T}{\bf P}_{0}{\bf E}\big).
\end{gather*}
For short, denote
\begin{gather*}
\Lambda({\bf N})=
{\bf N}(\delta){\bf H}(\delta)+\delta {\bf P}_0^{\perp}{\bf E}
 \quad {\rm and} \quad
\lambda({\bf N})=\big\|\mathcal{S}\Lambda({\bf N})\big\|_{\max}\, .
\end{gather*}
\begin{example}
\label{ex:sum_exp_sign_SSA}
{\it Signals of exponential type}. (See Example \ref{ex:sum_exp_sign} of Section \ref{ssect:sign_incr_exp}.)\\
Consider the signal series \eqref{eq:exp_sum} with  $\beta_k\!\neq\!0$ and decreasing $a_k\!>\!\!1$. Then $\|{\bf H}\|\asymp a_1^N$.
Denote $\theta=a_1^{3/2}/a_p^2$\,.
\item[1.]
Let the ``noise series'' be defined by $e_n=\sum_{l=1}^{m}\gamma_l b_l^n$ with $\gamma_l\neq 0$ and let $b\stackrel{\rm def}=\max_{1\leq l\leq m}\,|b_l|>1$.
Thus $\|{\bf E}\|\asymp b^{N}$.
If $b\theta<1$, then $\|{\bf E}\|=o(\|{\bf H}\|)$,
\begin{gather*}
\big\|{\bf P}_0^{\perp}(\delta)-{\bf P}_0^{\perp}-\delta\,{\bf V}_0^{(1)}\big\|\, \big\|{\bf H}(\delta)\big\|
=O(b\theta)^{2N}
\end{gather*}
for large $N$,  and
\begin{gather*}
 \big\|\widetilde{\mathrm{F}}_N(\delta)-\mathrm{F}_N\big\|_{\max}=\lambda\big(\delta{\bf V}_0^{(1)}\big)+
 O\left((b\theta)^{2N}\right).
\end{gather*}
\item[2.]
Consider a polynomial ``noise'' of order $m$ defined in the manner of \eqref{eq:polyn}.

a) If $L/N\rightarrow \alpha\in (0,1)$, then $\|{\bf E}\|\asymp N^{m+1}=o(\|{\bf H}\|)$.
If $\theta<1$, then
\begin{gather*}
\big\|\widetilde{\mathrm{F}}_N(\delta)-\mathrm{F}_N\big\|_{\max}=\lambda\big(\delta{\bf V}_0^{(1)}\big)+
O\big(N^{2m+2}\theta^{2N}\big).
\end{gather*}

b) If $N\rightarrow \infty$ and either $L$ or $K$ is a constant, then $\|{\bf E}\|\asymp N^{m+1/2}=o(\|{\bf H}\|)$.
 If  $\theta<1$, then
\begin{gather*}
\big\|\widetilde{\mathrm{F}}_N(\delta)-\mathrm{F}_N\big\|_{\max}=\lambda\big(\delta{\bf V}_0^{(1)}\big)+
O\big(N^{2m+1}\theta^{2N}\big).
\end{gather*}
\item[3.]
Let the oscillating ``noise series'' be defined  in the manner of \eqref{eq:osc}.
Then $\|{\bf E}\|\asymp \sqrt{LK}=o(\|{\bf H}\|)$. If $\theta<1$, then
\begin{gather}
\label{eq:hanel_exp_osc}
\big\|\widetilde{\mathrm{F}}_N(\delta)-\mathrm{F}_N\big\|_{\max}=\lambda\big(\delta{\bf V}_0^{(1)}\big)+
O\big(LK \theta^{2N}\big).
\end{gather}
\end{example}

\begin{example} {\it Oscillating signals.} (See Example \ref{ex:osc_ser} of Section \ref{sssect:osc_ser}.)\\
Let the oscillating signal be defined by \eqref{eq:osc}. Then $\|{\bf H}\|\asymp \sqrt{LK}$.
If the ``noise series'' is also oscillating but with frequencies different  from those of a signal, then
$\|{\bf E}\|\asymp \sqrt{LK}\asymp \|{\bf H}\|$ and
\begin{gather*}
\big\|\widetilde{\mathrm{F}}_N(\delta)-\mathrm{F}_N\big\|_{\max}=\lambda\big({\bf L}(\delta)\big)+
O\big(\sqrt{LK}/\min(L^2,K^2)\big)
\end{gather*}
under assumption that $|\delta|<\delta_0$ for some $\delta_0>0$. Note that the operator ${\bf L}(\delta)$ is defined
by \eqref{eq:L_delta}.
\end{example}

\begin{example}{\it Random stationary ``noise'' series.} (See Example \ref{ex:WN_LeqK} and proposition
\ref{prop:gen_wnoise_LeqK},
\ref{prop:gen_wight_noise_coarse} of Section \ref{sssect:osc_ser}.)
\item[1.]
For a signal series of exponential type  and under related conditions the following
results hold.

a) If $N\rightarrow \infty$ and $L\rightarrow \infty$, then  almost surely $\|{\bf E}\|/\sqrt{N\ln{N}}<c$ for $N>N_0(\omega)$
and  some constant $c>0$. Therefore, $\|{\bf E}\|=o(\|{\bf H}\|)$ with probability 1.
If $\theta<1$, then
\begin{gather*}
\big\|\widetilde{\mathrm{F}}_N(\delta)-\mathrm{F}_N\big\|_{\max}=\lambda\big(\delta{\bf V}_0^{(1)}\big)+
O\left({N\ln N}\,\theta^{2N}\right)
\end{gather*}
almost surely for $N>N_0(\omega,\delta)$.

b) If  $L=L_0={\rm const}$, then $\|{\bf E}\|/\sqrt{N}$ tends a.s. to the positive constant $c_0$ and
 \begin{gather*}
\big\|\widetilde{\mathrm{F}}_N(\delta)-\mathrm{F}_N\big\|_{\max}=\lambda\big(\delta{\bf V}_0^{(1)}\big)+
O\big(N \theta^{2N}\big)
 \end{gather*}
 with probability 1 for $N>N_0(\omega,\delta)$.
\item[2.]
 Let a polynomial signal of order $p$ be defined by \eqref{eq:polyn}.
If $L/N\rightarrow \alpha\in (0,1)$, then $\|{\bf H}\|\asymp N^{p+1}$ while $\|{\bf E}\|$ has a.s. the order of growth
$\sqrt{N\ln{N}}$. Therefore,
\begin{gather*}
\big\|\widetilde{\mathrm{F}}_N(\delta)-\mathrm{F}_N\big\|_{\max}=\lambda\big(\delta{\bf V}_0^{(1)}\big)+
O\left({\ln N}/N^{p}\right)
\end{gather*}
almost surely for $N>N_0(\omega,\delta)$.
\item[3.]
Consider the oscillating signal defined by \eqref{eq:osc} and the white-noise ``noise series'' as in Proposition
\ref{prop:gen_wight_noise_coarse}. If $L=L_0={\rm const}$, then $\|{\bf H}\|\asymp \sqrt{N}$ and $\|{\bf E}\|/\sqrt{N}
\rightarrow 1$ a.s. Thus $\|{\bf H}(\delta)\|=O(\sqrt N)$ with probability~1. Since almost surely
\begin{gather*}
\limsup_{N>N_0(\omega)}{N (\ln \ln N)^{-1}}\,\big\|{\bf P}_0^{\perp}(\delta)-{\bf P}_0^{\perp}-{\bf T}(\delta)\big\|
<c{''}\delta^2
\end{gather*}
with a positive constant $c{''}$, then
\begin{gather*}
\big\|\widetilde{\mathrm{F}}_N(\delta)-\mathrm{F}_N\big\|_{\max}=\lambda\big({\bf T}(\delta)\big)+
O\big({\ln \ln N}/\sqrt{N}\big)
\end{gather*}
almost surely for $N>N_0(\omega,\delta)$ and $|\delta|<\delta_0$.
\end{example}

Let us present two simple examples where these considerations are followed to the end.

\subsubsection{Reconstruction: constant signal and saw ``noise series''}
\label{sssect:rec_const_saw}
Consider the signal $x_n=1$ and the ``noise'' $e_n=(-1)^n$. This example was studied in detail in
Section \ref{ssect:const_saw}, where the main term ${\bf M}(\delta)$
of the difference ${\bf P}_0^\perp(\delta) - {\bf P}_0^\perp$
was derived in explicit form, see \eqref{eq:M(delta)},
under condition $|\delta|<1/2$. Here we use results and notation of this section.

\begin{proposition}
\label{prop:LK_even_odd}
1. If both $L$ and $K$ are even, then $\widetilde{\mathrm{F}}_N(\delta)=\mathrm{F}_N$.\\
2. If $L$ is even, $K$ is odd, $K\rightarrow \infty$ and $L=o(K^3)$, then
\begin{gather}
\label{eq:rec_Leven}
\big\|\widetilde{\mathrm{F}}_N(\delta)-\mathrm{F}_N\big\|_{\max}=
O\left(\max\big(K^{-1},\sqrt{L/K^3}\big)\right).
\end{gather}
3. If $K$ is even, $L$ is odd, $L\rightarrow \infty$ and $K=o(L^3)$, then
\begin{gather}
\label{eq:rec_Keven}
\big\|\widetilde{\mathrm{F}}_N(\delta)-\mathrm{F}_N\big\|_{\max}=
O\left(\max\big(L^{-1},\sqrt{K/L^3}\big)\right).
\end{gather}
4. If both $L$ and $K$ are odd, $L=o(K^3)$ and $K=o(L^3)$, then
\begin{gather}
\label{eq:rec:LKodd}
\big\|\widetilde{\mathrm{F}}_N(\delta)-\mathrm{F}_N\big\|_{\max}=
O\left(\min\big(\sqrt{L/K^3},\sqrt{K/L^3}\big)\right).
\end{gather}
\end{proposition}

\begin{remark}
The best order of convergence to zero of the right-hand side of \eqref{eq:rec_Leven} is $N^{-1}$. It is achieved under the
choice $L\!=\!L_0\!=\!{\rm const}$. The same result for \eqref{eq:rec_Keven} is attained by the choice $K\!=\!K_0$ and for
\eqref{eq:rec:LKodd} by the choice $L\!\sim\!N/2$.
\end{remark}

\subsubsection{Reconstruction: exponential signal series and constant ``noise series''}
\label{sssect:rec_exp_const}
Consider the signal series $x_n\!=\!a^n$ with $a>1$ and the ``noise series'' $e_n\equiv 1$. Then $\theta=\sqrt a$ and applying
\eqref{eq:hanel_exp_osc} we see~that
\begin{gather*}
\big\|\widetilde{\mathrm{F}}_N(\delta)-\mathrm{F}_N\big\|_{\max}=\lambda\big(\delta{\bf V}_0^{(1)}\big)+
O\big(N^{2}a^{-N}\big)
\end{gather*}
in the case $L/N\rightarrow \alpha\in (0,1)$ and that
\begin{gather*}
\big\|\widetilde{\mathrm{F}}_N(\delta)-\mathrm{F}_N\big\|_{\max}=\lambda\big(\delta{\bf V}_0^{(1)}\big)+
O\big(Na^{-N}\big).
\end{gather*}
in the case when $N\rightarrow \infty$ and either $L$ or $K$ is a constant.

\begin{proposition}
\label{prop:exp_const_hankel}
If $L/N\rightarrow \alpha\in (0,1)$ or if $N\rightarrow \infty$ and either $L$ or $K$ is a constant, then
\begin{gather*}
\limsup_N\big\|\widetilde{\mathrm{F}}_N(\delta)-\mathrm{F}_N\big\|_{\max}>0.
\end{gather*}
\end{proposition}

\begin{remark}
\label{rem:err_exp_saw}
The accurate analysis of the matrix
$\mathcal{S}\Lambda_1\big(\delta{\bf V}_0^{(1)}\big)$ (we omit calculations
due to their technical character) gives much more information on the
behavior of the series $\widetilde{\mathrm{F}}_N(\delta)$ as
$N\rightarrow \infty$. Let us consider the simplest case of odd $N$ with
$L=K$.
Denote $b=(a+1)/(a-1)$ and $\Delta f_n(\delta)=\widetilde{f}_n(\delta)-f_l(\delta)$.
Then the following asymptotic results hold.
\begin{enumerate}
\item
If $l<L$, then
\begin{gather}
\label{eq:small_L}
\max_{0\leq l< L}\left|\Delta {f}_l(\delta)- 2\delta\, b\,
\frac{a^{-L}\big(a^{l+1}\!-\!1\big)}{l+1}\right|
\!=\!|\delta|O\big(a^{-2L}\big)
\end{gather}
as $L\rightarrow \infty$. In particular, $\Delta {f}_{0}(\delta)=2\delta(a+1)a^{-L}+O(a^{-2L})$,
$\Delta {f}_{L-1}(\delta)\sim 2\delta\, b L^{-1}$, and
\begin{gather*}
 \Delta {f}_{l}(\delta)\sim 2\delta b\,\frac{a}{\lambda_0L}\ a^{-(1-\lambda_0)L}
\end{gather*}
in the case $l/L\rightarrow\lambda_0\in (0,1)$.
\item
If $l=2L-k\geq L$, then
\begin{gather}
\begin{split}
&\max_{1<k\leq L}\bigg|\Delta {f}_{2L-k}(\delta)\!-\!2\delta b\,
\frac{a(a^{k-1}\!-\!1)\!-\!(a^2-1)(k-1)}{(k-1)a^k}\bigg|\\
&=|\delta|\,O\big(La^{-L}\big)
\end{split}
\label{eq:big_L}
\end{gather}
as $L\rightarrow \infty$.
This means that for any small $\varepsilon > 0$
\begin{gather*}
  \max_{l\in [L,\ (2-\varepsilon)L]}\left|\Delta {f}_{l}(\delta)-2\delta b\,\frac{1}{2L-l}\right|=
  |\delta|\,O\big(L^{-2}\big)\,.
\end{gather*}
In particular, if $l/L\rightarrow\lambda_0\in [1,\ 2)$, then $L\Delta f_{l}^{(1)}\rightarrow 2\delta\, b(2-\lambda_0).$
\end{enumerate}
These results improve the accuracy of Proposition \ref{prop:exp_const_hankel} in the case $L=K$. Namely, though the sequence
$\max_{\,0\leq l<N}|\widetilde{f}_l(\delta)-f_l|$ does not tend to zero,
still
$\max_{\,0\leq l<(1-\varepsilon)N}|\widetilde{f}_l(\delta)-f_l|\rightarrow 0$ as
$N\rightarrow \infty$ for any small $\varepsilon$.
\end{remark}

\begin{figure}[h]
    \center
   \includegraphics[width=12.0cm, height=6.0cm]{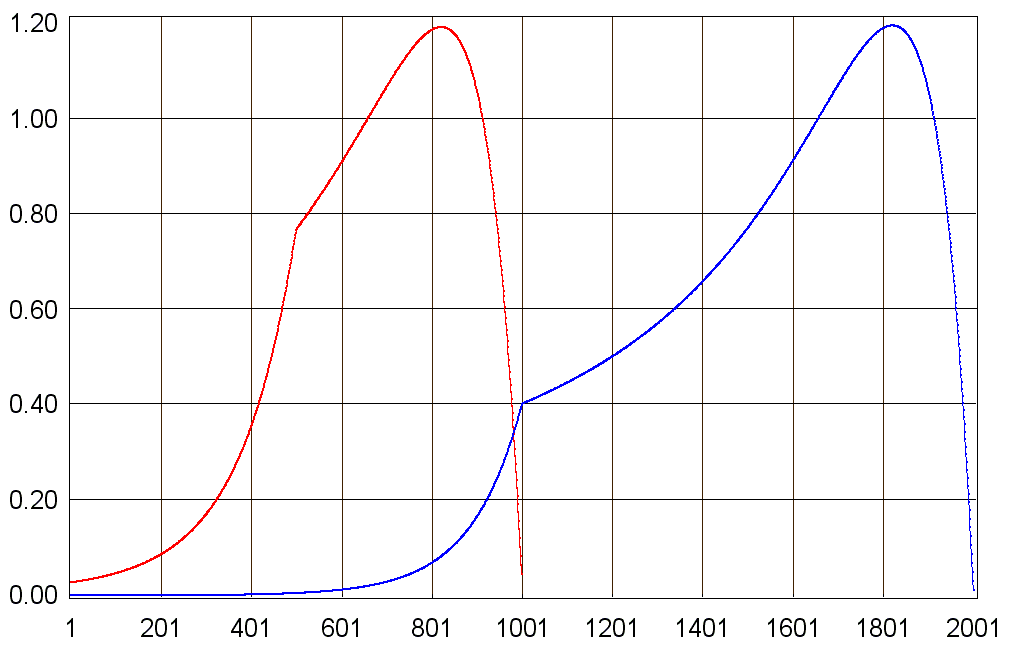}
    \caption{Errors of reconstruction.}
    \label{pic:exp_const}
\end{figure}

Computer experiments confirm these considerations. For example, two lines of Fig. \ref{pic:exp_const} present SSA-computed
reconstruction errors $\Delta {f}_n(\delta)$ for the signal $f_n=(1.01)^n$ with the noise series $e_n=1$ and
$\delta=1$. Red line corresponds to $N=999$ and $L=K=500$, blue line describes the case
$N=1999$ and $L=K=1000$.

Both curves illustrate results of Remark \ref{rem:err_exp_saw}. In particular, analytical main terms of the
differences $\Delta {f}_n(\delta)$ presented in
\eqref{eq:small_L} and \eqref{eq:big_L} perfectly correspond to the computed data of Fig. \ref{pic:exp_const}.

\section{Appendix 1: Proofs}
{\bf Proof of Theorem \ref{theor:zero_pert}}

\vspace{1mm}
\noindent
1. {\it The equivalence between \eqref{eq:zero_pert} and \eqref{eq:zero_pert1}.}
Since
\begin{gather*}
{\bf S}_0{\bf B}(\delta){\bf P}_0=\delta\,{\bf S}_0{\bf H}{\bf E}^{\rm T}{\bf P}_0+
\delta^2\,{\bf S}_0{\bf E}{\bf E}^{\rm T}{\bf P}_0,
\end{gather*}
then \eqref{eq:zero_pert1} automatically follows from \eqref{eq:zero_pert}. On the other hand, if
${\bf S}_0{\bf H}{\bf E}^{\rm T}{\bf P}_0+\delta\,{\bf S}_0{\bf E}{\bf E}^{\rm T}{\bf P}_0={\bf 0}_{LL}$ for all
$\delta\in (-\varepsilon,\varepsilon)$, then we obtain the equalities \eqref{eq:zero_pert}.

2. {\it The equivalence between \eqref{eq:zero_pert} and \eqref{eq:sum_SHEP}.}
It is easy to see that \eqref{eq:sum_SHEP}
follows from \eqref{eq:zero_pert}. Let us demonstrate that the equality
${\bf S}_0{\bf H}{\bf E}^{\rm T}{\bf P}_0={\bf 0}_{LL}$
follows from the equality
\begin{gather}
\label{eq:sum_SHEP2}
{\bf S}_0{\bf H}{\bf E}^{\rm T}{\bf P}_0 + {\bf P}_0{\bf H}{\bf E}^{\rm T}{\bf S}_0={\bf 0}_{LL}
\end{gather}
and that the equality ${\bf S}_0{\bf E}{\bf E}^{\rm T}{\bf P}_0={\bf 0}_{LL}$ is the consequence of
${\bf S}_0{\bf E}{\bf E}^{\rm T}{\bf P}_0 + {\bf P}_0{\bf E}{\bf E}^{\rm T}{\bf S}_0={\bf 0}_{LL}$.
Both assertions have
similar proofs. Suppose that \eqref{eq:sum_SHEP2} holds and denote
$X={\bf H}{\bf E}^{\rm T}{\bf P}_0Z$ and $Y={\bf H}{\bf E}^{\rm T}{\bf S}_0Z$ for some $Z\in \mathbb{R}^L$.
Then the equality
${\bf S}_0X+{\bf P}_0Y=0_L$ follows from \eqref{eq:sum_SHEP2}.
Since ${\bf S}_0X\in \mathbb{U}_0^\perp$ and ${\bf P}_0Y \in \mathbb{U}_0$, the latter vectors are
orthogonal. Thus
${\bf S}_0{\bf H}{\bf E}^{\rm T}{\bf P}_0Z={\bf P}_0{\bf H}{\bf E}^{\rm T}{\bf S}_0Z=0_L$ for all $Z\in \mathbb{R}^L$.

3.
{\it The equivalence between  \eqref{eq:zero_pert} and the equality ${\bf P}_0^{\perp}(\delta)={\bf P}_0^{\perp}$.}
 Assume that \eqref{eq:zero_pert} holds. Then
${\bf S}_0{\bf B}(\delta){\bf P}_0={\bf P}_0{\bf B}(\delta){\bf S}_0={\bf 0}_{LL}$.
 Consider the operator
\begin{gather}
\label{eq:term_Wp}
{{\bf S}}_{0}^{(l_1)}{{\bf B}}(\delta){{\bf S}}_{0}^{(l_2)}\ldots
\,{\bf S}_{0}^{(l_{p})}{{\bf B}}(\delta){\bf S}_{0}^{(l_{p+1})},
\end{gather}
see \eqref{eq:Wp}. Since $l_1+\ldots+l_{p+1}=p$, then there exists $j$ such that $l_j\!>\!0$ and either
$l_{j-1}\!=~0$ or $l_{j+1}\!=\!0$. Suppose that $l_{j+1}\!=\!0$. (The case $l_{j-1}\!=\!0$ is studied analogously.) Since
${\bf S}_0^{(0)}\!=\!-\!{\bf P}_0$ and ${\bf S}_0^{(l_j)}\!=\!{\bf S}^{l_j}$, then
${\bf S}_0^{(l_j)}{\bf B}(\delta){\bf S}_0^{(l_{j+1})}={\bf 0}$ and
the operator
\eqref{eq:term_Wp} is the zero operator.
Hence
${\bf W}_p(\delta)={\bf 0}_{LL}$ for all $p$, and
${\bf P}_0^{\perp}(\delta)={\bf P}_0^{\perp}$.

In order to prove the converse assertion we consider the expansion  \eqref{eq:pert_0}. Since
${\bf P}_0^{\perp}(\delta)={\bf P}_0^{\perp}$ for all $\delta\in (-\delta_0, \delta_0)$, then the coefficients
${\bf V}_0^{(n)}$ of the power series \eqref{eq:pert_0} are zero operators. In particular, ${\bf V}_0^{(1)}={\bf 0}_{LL}$
and \eqref{eq:sum_SHEP2} holds due to \eqref{eq:lin_proj}.
Let us
consider the coefficient
${\bf V}_0^{(2)}$. As it is already proved, ${\bf S}_0{\bf A}^{(1)}{\bf P}_0={\bf P}_0{\bf A}^{(1)}{\bf S}_0={\bf 0}_{LL}$.
Thus \eqref{eq:sec_order} shows that
${\bf 0}_{LL}={\bf V}_0^{(2)}={\bf P}_0{\bf A}^{(2)}{\bf S}_0+{\bf S}_0{\bf A}^{(2)}{\bf P}_0$,
and we get \eqref{eq:sum_SHEP}. The latter equality is equivalent to \eqref{eq:zero_pert}.

4.
{\it The equivalence between \eqref{eq:zero_pert}  and \eqref{eq:zero_pert2}.}
Let us prove that the equality ${\bf H}{\bf E}^{\rm T}{\bf P}_0={\bf 0}_{LL}$ follows from ${\bf S}_0{\bf H}{\bf E}^{\rm T}{\bf P}_0={\bf 0}_{LL}$.
(The converse statement holds automatically.)
Suppose that ${\bf S}_0{\bf H}{\bf E}^{\rm T}{\bf P}_0={\bf 0}_{LL}$ and let
$Z\in \mathbb{R}^L$. Note that vectors ${\bf P}_\mu {\bf H}{\bf E}^{\rm T}{\bf P}_0Z$ are pairwise orthogonal  for
different $\mu\in \Sigma$. Therefore, since
\begin{gather*}
{\bf S}_0{\bf H}{\bf E}^{\rm T}{\bf P}_0Z= \sum_{\mu>0}{\mu^{-1}}\, {\bf P}_\mu {\bf H}{\bf E}^{\rm T}{\bf P}_0Z = 0_L,
\end{gather*}
then $ {\bf P}_\mu {\bf H}{\bf E}^{\rm T}{\bf P}_0Z = 0_L$ for all non-zero $\mu\in \Sigma$.
In view of the equality ${\bf H}=\sum_{\mu>0}{\bf P}_\mu {\bf H}$ we see that ${\bf H}{\bf E}^{\rm T}{\bf P}_0={\bf 0}$.

Let us prove that the equality ${\bf H}^{\rm T}{\bf E}{\bf E}^{\rm T}{\bf P}_0={\bf 0}_{LK}$ holds iff
${\bf S}_0{\bf E}{\bf E}^{\rm T}{\bf P}_0={\bf 0}_{LL}$.
If ${\bf H}^{\rm T}{\bf E}{\bf E}^{\rm T}{\bf P}_0={\bf 0}_{LK}$, then
\begin{gather*}
{\bf 0}_{LL}={\bf H}{\bf H}^{\rm T}{\bf E}{\bf E}^{\rm T}{\bf P}_0=\sum_{\mu>0}\mu {\bf P}_\mu {\bf E}{\bf E}^{\rm T}{\bf P}_0,
\end{gather*}
${\bf P}_\mu {\bf E}{\bf E}^{\rm T}{\bf P}_0={\bf 0}_{LL}$ for all non-zero $\mu\in \Sigma$, and
${\bf S}_0{\bf E}{\bf E}^{\rm T}{\bf P}_0={\bf 0}_{LL}$. The converse statement is proved in the same manner due to the
equality
$
{\bf H}^{\rm T}{\bf E}{\bf E}^{\rm T}{\bf P}_0=\sum_{\mu>0}{\bf H}^{\rm T}{\bf P}_\mu {\bf E}{\bf E}^{\rm T}{\bf P}_0.
$
This completes the entire proof.

\vspace{3mm}
{\bf Proof of Lemma \ref{lem:weak_sep}}

\vspace{1mm}
\noindent
1. Let ${\bf H}^{\rm T}{\bf E}={\bf 0}_{KK}$. Then ${\bf H}^{\rm T}{\bf H}{\bf E}X=0_L$ for any $X\in \mathbb{R}^K$.
If $0\notin \Sigma$, then ${\bf E}X=0_L$ and ${\bf E}={\bf 0}_{KL}$. Therefore, $0\in \Sigma$ and
${\bf E}X\in \mathbb{U}_0$. Thus ${\bf P}_0{\bf E}X={\bf E}X$ and ${\bf P}_0{\bf E}={\bf E}$. Analogously,
${\bf P}_\mu{\bf E}={\bf P}_\mu{\bf P}_0{\bf E}={\bf 0}_{KL}$ for all non-zero $\mu\in \Sigma$.\\
2. Since ${\bf H}=\sum_{\mu>0}{\bf P}_\mu {\bf H}$, then
$
{\bf H}^{\rm T}{\bf E}={\bf H}^{\rm T}\left(\sum_{\mu>0}{\bf P}_\mu {\bf E}\right),
$
and the second assertion is proved as well.

\vspace{3mm}
{\bf Proof of Theorem \ref{theor:gen_upper}}

\vspace{1mm}
\noindent
Let us start with the following lemma.
\begin{lemma}
\label{lem:comb}
If $0<\beta<1/4$ and $k\geq 0$, then
\begin{gather}
\label{eq:comb_noneq}
\sum_{p=k}^\infty {2p \choose p}\ \beta^p\leq C\frac{(4\beta)^k}{1-4\beta}\ ,
\end{gather}
where $C=e^{1/6}/\sqrt\pi\approx 0.667$.
\end{lemma}
\begin{proof}
Since $n!=\sqrt{2\pi}\,n^{n+1/2}e^{-n}\,e^{\theta_n/12n}$ with  $0< \theta_n<1$, then
\begin{gather*}
\begin{split}
&{2p \choose p}
\!=\!\frac{1}{\sqrt{2\pi}}\, \frac{\sqrt{2p}}{p}\, \frac{4^p p^{2p}}{p^{2p}}e^{(\theta_{2p}-\theta_p)/6p}
\!\leq\!
\frac{1}{\sqrt{\pi p}}\, 4^p\, e^{1/6p}\!\leq \!C4^p,
\end{split}
\end{gather*}
and the result becomes clear.
\end{proof}
Now let us demonstrate Theorem \ref{theor:gen_upper}.
Due to Theorem  \ref{theor:full_decomp}, we can use the expansion \eqref{eq:pert_0_B}.

In the same way as in the demonstration of Theorem \ref{theor:zero_pert},
 consider the term \eqref{eq:term_Wp}
 and take $j$ such that $l_j>0$  and $l_{j+1}=0$ (the case $l_{j-1}=0$ is quite analogous).
Since $l_1+\ldots+l_{p+1}=p$ and $\|{\bf S}^{(k)}_0\|=1/\mu_{\min}^k$ for any $k\geq 0$,
then the norm of the operator
\eqref{eq:term_Wp} can be estimated from above:
\begin{gather}
\begin{split}
&\left\|{{\bf S}}_{0}^{(l_1)}{{\bf B}}(\delta){\bf S}_{0}^{(l_2)}\ldots\,{{\bf B}}(\delta){{\bf S}}_{0}^{(l_{p+1})}\right\|
\\
&=
\left\|{{\bf S}}_{0}^{(l_1)}{{\bf B}}(\delta)\ldots\,
{{\bf S}}_{0}^{(l_{j}-1)}{\bf S}_{0}{{\bf B}}(\delta){\bf P}_{0}\ldots
{{\bf B}}(\delta){{\bf S}}_{0}^{(l_{p+1})}\right\|
\\
&\leq \big\|{\bf S}_{0}{{\bf B}}(\delta){\bf P}_{0} \big\|
\left({\big\|{\bf B}(\delta)\big\|}\big/{\mu_{\min}}\right)^{p-1}.
\end{split}
\label{eq:gen_Wp_upper}
\end{gather}

Denote $\Im_p$ the number of  nonnegative integer vectors $(l_1,\dots,l_{p+1})$
such that $l_1+\ldots+l_{p+1}=p$. Then $\Im_p=
{2p \choose p}$.
In view of \eqref{eq:pert_0_B} and \eqref{eq:gen_Wp_upper},
\begin{gather*}
\begin{split}
&\big\|{\bf P}_0^{\perp}(\delta)-{\bf P}_0^{\perp}\big\|
\\
&\leq \big\|{\bf S}_{0}{{\bf B}}(\delta){\bf P}_{0} \big\|\,
\sum_{p=1}^\infty\,
\sum_{l_1+\ldots+l_{p+1}=p,\  l_j\geq 0}\left({\|{{\bf B}}(\delta)\|}\big/{\mu_{\min}}\right)^{p-1}
\\
&=
\big\|{\bf S}_{0}{{\bf B}}(\delta){\bf P}_{0} \big\|\,\sum_{p=1}^\infty {2p \choose p}\,
\left({\|{{\bf B}}(\delta)\|}\big/{\mu_{\min}}\right)^{p-1}.
\end{split}
\end{gather*}
Now we get the result with the help of
\eqref{eq:comb_noneq}.

\vspace{3mm}
{\bf Proof of Corollary \ref{cor:coarse_proj}}

\vspace{1mm}
\noindent
Since $\|{\bf S}_0 {\bf B}(\delta){\bf P}_0\|\leq \|{\bf B}(\delta)\|/\mu_{\min}$, then \eqref{eq:coarse_proj}
automatically follows from \eqref{eq:gen_upper}.
The inequality
${\|{{\bf B}}(\delta)\|}/{\mu_{\min}}\leq {B(\delta)}/{\mu_{\min}}$ was already established.
The chain of inequalities
\begin{gather*}
\begin{split}
&\frac{\|{{\bf B}}(\delta)\|}{\mu_{\min}}\leq \frac{B(\delta)}{\mu_{\min}}\leq
2|\delta|\,\frac{\|{\bf H}{\bf E}^{\rm T}\|}{\mu_{\min}}+\delta^2\, \frac{\nu_{\max}}{\mu_{\min}}
\\
&\leq 2|\delta|\,
\frac{\|{\bf H}\|\,\|\bf E\|}{\mu_{\min}}+\delta^2\, \Theta_1^2\Theta_2=
2|\delta|\,{\Theta_1 \Theta_2}+\delta^2 \Theta_1^2\Theta_2
\end{split}
\end{gather*}
provides the result.

\vspace{3mm}
{\bf Proof of Proposition \ref{prop:min_angle}}\\
We begin with the following lemma.
\begin{lemma}
\label{lem:range_Ax}
Consider the matrix ${\bf M}: \mathbb{R}^K\mapsto \mathbb{R}^L$ and let $\mathbb{U}_M$ stand for the linear space
spanned by  columns of the matrix ${\bf M}$. Denote $\sigma_{\min}$ and $\sigma_{\max}$ the minimal and maximal singular
values of ${\bf M}$. Then
\begin{gather}
\label{eq:range_max}
\big\{{\bf M}X, \, \|X\|\leq 1/\sigma_{\max}\big\}\subset \big\{Y\in \mathbb{U}_M, \, \|Y\|\leq 1\big\}
\end{gather}
and
\begin{gather}
\label{eq:range_min}
\big\{{\bf M}X, \, \|X\|\leq 1/\sigma_{\min}\big\}\supset \big\{Y\in \mathbb{U}_M, \, \|Y\|\leq 1\big\}.
\end{gather}
\end{lemma}
\begin{proof}
The inclusion \eqref{eq:range_max} straightforwardly follows from the inequality
$\|{\bf M}X\|\leq \|{\bf M}\|\,\|X\|\leq 1$.
Let  us prove the inclusion \eqref{eq:range_min}.

Consider the SVD of the matrix ${\bf M}$: ${\bf M}=\sum _j\sigma_j U_j\,V_j^{\rm T}$. By definition,
left singular  vectors $U_j$ form the orthonormal basis of $\mathbb{U}_M$. (Analogously, right singular
vectors $V_j$  form the orthonormal basis of the linear space spanned by rows of the same matrix.)

Let us prove that for any $Y\in \mathbb{U}_M$ with $\|Y\|\leq 1$ there exists $X\in \mathbb{R}^K$ such that $\|X\|\leq 1/\sigma_{\min}$
and ${\bf M}X=Y$. Then \eqref{eq:range_min} will be proved as well.

Note that $Y=\sum_{j}c_jU_j$ with $\sum_jc_j^2\leq 1$.
In the same manner, each $X\in \mathbb{R}^K$ can be expressed as $X=\sum_j\alpha_j V_j+\beta_0 W_0$, where
$\big<V_j,W_0\big>_{K}=0$ and $\big<\,,\,\big>_K$ stands for the inner product in $\mathbb{R}^K$.
Therefore, ${\bf M}X=\sum_j\sigma_j\alpha_j U_j$.
If we put $\alpha_j=c_j/\sigma_j$, then $\sum_j\alpha_j^2 \leq \sum_jc_j^2/\sigma^2_{\min}\leq 1/\sigma^2_{\min}$ and the
proof is complete.
\end{proof}

Now we turn to the demonstration of Proposition \ref{prop:min_angle}.
By definition,
\begin{gather*}
\cos(\theta_{\min})\!=\!\max_{Y_1\in \mathbb{U}_1, \|Y_1\|\leq 1}\,\max_{Y_2\in \mathbb{U}_2, \|Y_2\|\leq 1}
\big<Y_1,Y_2\big>_{L}.
\end{gather*}
In view of Lemma \ref{lem:range_Ax},
\begin{gather*}
\begin{split}
&\|{\bf M}_1\|\,\|{\bf M}_2\|\,\cos(\theta_{\min})
\\
&\leq
\max_{X_1\in \mathbb{R}^K, \|X_1\|\leq 1}\,\max_{X_2\in \mathbb{R}^K, \|X_2\|\leq 1}
\big<{\bf M}_1X_1,{\bf M}_2X_2\big>_{L}
\\
&=
\max_{X_1\in \mathbb{R}^K, \|X_1\|\leq 1}\,\max_{X_2\in \mathbb{R}^K, \|X_2\|\leq 1}
\big<{\bf M}_1^{\rm T}{\bf M}_2X_1,X_2\big>_{K}
\\
&=\big\|{\bf M}_1^{\rm T}{\bf M}_2\big\|.
\end{split}
\end{gather*}
If we use \eqref{eq:range_min} instead of \eqref{eq:range_max}, then we get
\begin{gather*}
\begin{split}
&\sigma_1^{(\min)}\,\sigma_2^{(\min)}\,\cos(\theta_{\min})
\\
&\geq \max_{X_1\in \mathbb{R}^K, \|X_1\|\leq 1}\,\max_{X_2\in \mathbb{R}^K, \|X_2\|\leq 1}
\big<{\bf M}_1X_1,{\bf M}_2X_2\big>_{L}
\\
&=\big\|{\bf M}_1^{\rm T}{\bf M}_2\big\|.
\end{split}
\end{gather*}
This completes the proof.

\vspace{3mm}
{\bf Proof of Corollary \ref{cor:coarse_orth}}

\vspace{1mm}
\noindent
It is easy to see that \eqref{eq:coarse_orth} follows from \eqref{eq:gen_upper}. Besides,
\begin{gather}
\label{eq:norm_SB}
\big\|{\bf S}_0{\bf B}(\delta)\big\|\leq 2|\delta|\, {\|{\bf H}{\bf E}^{\rm T}\|}\big/{\mu_{\min}}+
\delta^2 \big\|{\bf S}_0 {\bf A}^{(2)}\big\|.
\end{gather}
Proposition \ref{prop:min_angle} shows that
\begin{gather*}
{\|{\bf H}{\bf E}^{\rm T}\|}\big/{\mu_{\min}}\leq \cos(\theta_{r}){\sqrt{\mu_{\max}\nu_{\max}}}/{\mu_{\min}}=
\Theta_1\Theta_2 \,\cos(\theta_{r}).
\end{gather*}

As for the term $\big\|{\bf S}_0 {\bf A}^{(2)}\big\|$, it is worth to mention that the columns of the matrix
 ${\bf S}_0={\bf S}_0^{\rm T}$
span the same linear space $\mathbb{U}_0^\perp$ as the columns of the matrix ${\bf H}$. Analogously, the columns of
${\bf A}^{(2)}$ span the same linear space as the columns of ${\bf E}$. Therefore, it follows from
Proposition \ref{prop:min_angle} that
$\big\|{\bf S}_0 {\bf A}^{(2)}\big\|\leq \cos(\theta_l){\nu_{\max}}/{\mu_{\min}}= \Theta_1^2\Theta_2\cos(\theta_l)$,
and the proof is complete.

\vspace{3mm}
{\bf Proof of Theorem \ref{theor:coarse_proj}}

\vspace{1mm}
\noindent
Acting in the same manner as in Theorem \ref{theor:gen_upper} and taking into consideration that
$\|{\bf S}_0{\bf B}(\delta)\|\leq \|{\bf B}(\delta)\|/{\mu_{\min}}$ we get the inequality
\begin{gather*}
\big\|{\bf P}_0^{\perp}(\delta)-{\bf P}_0^{\perp}-{\bf W}_1(\delta)\big\|\leq
\sum_{p=2}^\infty {2p \choose p}\, \left(\frac{\|{{\bf B}}(\delta)\|}{\mu_{\min}}\right)^p.
\end{gather*}
Applying Lemma \ref{lem:comb} and the condition $\|{\bf B}(\delta)\|<\mu_{\min}/4$
we obtain \eqref{eq:coarse_proj_mainterm}.

\vspace{3mm}
{\bf Proof of Theorem \ref{theor:precise_proj}}

\vspace{1mm}
\noindent
First of all, we can use the
decomposition \eqref{eq:pert_0_B}. Let us extract the operator ${\bf L}(\delta)$ from the
right-hand side  of \eqref{eq:pert_0_B}.

For fixed $p\geq 1$ we take the vector $(l_1,\ldots,l_{p+1})$ with $l_0=p$ and $l_j=0$ for $j\geq 1$ and
consider the related operator
\begin{gather*}
{\bf X}_1(p)=
(-1)^p\,{{\bf S}}_{0}^{(l_1)}{{\bf B}}(\delta){\bf S}_{0}^{(l_2)}\ldots\,{{\bf B}}(\delta){{\bf S}}_{0}^{(l_{p+1})}
\end{gather*}
It is clear that
\begin{gather*}
\begin{split}
&{\bf X}_1(p)=
(-1)^p\,{{\bf S}}_{0}^{(p)}{{\bf B}}(\delta){\bf S}_{0}^{(0)}\ldots\,{{\bf B}}(\delta){{\bf S}}_{0}^{(0)}
\\
&=
{{\bf S}}_{0}^{p}\,{{\bf B}}(\delta){\bf P}_{0}\ldots\,{{\bf B}}(\delta){{\bf P}}_{0}
=
{{\bf S}}_{0}^{p}\,{{\bf B}}(\delta){\bf P}_{0}\Big({\bf P}_{0}{{\bf B}}(\delta){{\bf P}}_{0}\Big)^{p-1}\,.
\end{split}
\end{gather*}
Since ${\bf P}_{0}{\bf B}(\delta){\bf P}_{0}=\delta^2{\bf A}^{(2)}_0$
and ${\bf S}_{0}^{p}=\sum_{\mu>0}{\bf P}_\mu/\mu^p$, then
\begin{gather*}
\begin{split}
&{\bf X}_1(p)=
\sum_{\mu>0}\mu^{-p}\,{\bf P}_\mu{{\bf B}}(\delta){\bf P}_{0}\left(\delta^2{\bf A}^{(2)}_0\right)^{p-1}
\\
&=\sum_{\mu>0}\mu^{-1}{\bf P}_\mu{{\bf B}}(\delta){\bf P}_{0}\,
\left({\delta^2{\bf A}^{(2)}_0}\big/{\mu}\right)^{p-1}.
\end{split}
\end{gather*}
Thus we obtain that
\begin{gather*}
\begin{split}
&\sum_{p=1}^\infty {\bf X}_1(p)=\sum_{\mu>0}\mu^{-1}{{\bf P}_\mu{{\bf B}}(\delta){\bf P}_{0}}\,
\sum_{p=1}^\infty\left({\delta^2{\bf A}^{(2)}_0}\big/{\mu}\right)^{p-1}
\\
&=
\sum_{\mu>0}\mu^{-1}{{\bf P}_\mu{{\bf B}}(\delta){\bf P}_{0}}\,
\left({\bf I}-\delta^2{\bf A}^{(2)}_0/{\mu}\right)^{-1}.
\end{split}
\end{gather*}
If we consider ${\bf X}_2(p)=
(-1)^p\,{{\bf S}}_{0}^{(0)}{{\bf B}}(\delta){\bf S}_{0}^{(0)}\ldots\,{{\bf B}}(\delta){{\bf S}}_{0}^{(p)}$ and act in
the same manner as for ${\bf X}_1(p)$, we get that
\begin{gather*}
\sum_{p=1}^\infty {\bf X}_2(p)=
\sum_{\mu>0}\left({\bf I}-\delta^2{\bf A}^{(2)}_0/{\mu}\right)^{-1}
{{\bf P}_0{{\bf B}}(\delta){\bf P}_{\mu}}\big/{\mu}\,.
\end{gather*}
Thus ${\bf L}(p)=\sum_{p\geq 1}\big({\bf X}_1(p)+{\bf X}_2(p)\big)$.

Now  consider the vector $(l_1,\ldots,l_{p+1})$ such that $\sum_{j=1}^{p+1}l_j=p$
and neither $l_1$ nor $l_{p+1}$
equals $p$. (It follows from this condition that $p>1$.)
In addition, consider $1\leq j\leq p+1$
such that  $l_j>0$ and either $l_{j-1}$ or $l_{j+1}$ equals zero. If $l_j=p$ and $1<j<p+1$, then
\begin{gather}
\begin{split}
&\big\|{\bf S}_{0}^{(0)}{\bf B}(\delta)\ldots
{\bf B}(\delta){\bf S}_{0}^{(p)}{{\bf B}}(\delta)\ldots\,{{\bf B}}(\delta){{\bf S}}_{0}^{(0)}\big\|
\\
&=\!
\big\|{\bf P}_{0}{{\bf B}}(\delta)\ldots
{\bf P}_{0}{\bf B}(\delta){\bf S}_{0}{\bf S}_{0}^{p-2}{\bf S}_{0}
{\bf B}(\delta){\bf P}_{0}\ldots {\bf B}(\delta){\bf P}_{0}\big\|
\\
&\leq\!
\|{\bf P}_{0}{\bf B}(\delta){\bf S}_{0}\|\, \|{\bf S}_0{\bf B}(\delta){\bf P}_{0}\|\,
\|{\bf S}_{0}^{p-2}\|\,\|{\bf B}(\delta)\|^{p-2}
\\
&\leq\!
\|{\bf S}_0{\bf B}(\delta){\bf P}_{0}\|\,
\|{\bf S}_0{\bf B}(\delta)\|\,\big({\|{\bf B}(\delta)\|}/{\mu_{\min}}\big)^{p-2}.
\end{split}
\label{eq:upper_nonextreme}
\end{gather}
If $l_j\neq p$, then there exists $k\neq j$ such that $l_k\neq 0$ and
we get  the same upper bound \eqref{eq:upper_nonextreme}.
For example, let $1<l_k<l_j<p+1$ and $l_{j+1}=0$. Then
\begin{gather*}
\begin{split}
&\left\|{{\bf S}}_{0}^{(l_1)}{{\bf B}}(\delta)\ldots{{\bf B}}(\delta){{\bf S}}_{0}^{(l_k)}\ldots\,
{{\bf S}}_{0}^{(l_j)}{{\bf B}}(\delta)\ldots
{{\bf B}}(\delta){{\bf S}}_{0}^{(l_{p+1})}\right\|
\\
&=
\left\|{{\bf S}}_{0}^{(l_1)}\ldots{{\bf B}}(\delta){\bf S}_0{{\bf S}}_{0}^{(l_k-1)}\ldots\,
{{\bf S}}_{0}^{(l_j-1)}{\bf S}_0{{\bf B}}(\delta){\bf P}_{0}\ldots
{{\bf S}}_{0}^{(l_{p+1})}\right\|
\\
&\leq
 \big\|{\bf S}_0{\bf B}(\delta)\big\|\,
\big\|{\bf S}_0{\bf B}(\delta){\bf P}_{0}\big\|\,\big\|{\bf B}\|^{p-2}\,\big\|{\bf S}_0\|^{l_1+\ldots+l_{p+1}-2}
\\
&=
\big\|{\bf S}_0{\bf B}(\delta)\big\|\,\big\|{\bf S}_0{\bf B}(\delta){\bf P}_{0}\big\|\,
\big({\big\|{\bf B}(\delta)\big\|}/{\mu_{\min}}\big)^{p-2}.
\end{split}
\end{gather*}
Other variants concerning $k,j$ provide the same result.
Therefore, applying \eqref{eq:pert_0_B} and Lemma \ref{lem:comb} we obtain that
\begin{gather*}
\begin{split}
&\big\|{\bf P}_0^{\perp}(\delta)-{\bf P}_0^{\perp}-{\bf L}(\delta)\big\|
\\
&=
\Big\|\sum_{p=2}^\infty (-1)^p
\sum_{\substack{l_1+\ldots+l_{p+1}=p\\  l_j\geq 0, \, l_1\neq p, \, l_{p+1}\neq p} }
{{\bf S}}_{0}^{(l_1)}{{\bf B}}(\delta){{\bf S}}_{0}^{(l_2)}\!\ldots {{\bf B}}(\delta){{\bf S}}_{0}^{(l_{p+1})}\Big\|
\\
&\leq
\sum_{p=2}^\infty\, \sum_{\substack{l_1+\ldots+l_{p+1}=p \\ l_j\geq 0, \, l_1\neq p, \, l_{p+1}\neq p}}
\|{\bf S}_0{\bf B}(\delta)\|\, \|{\bf S}_0{\bf B}(\delta){\bf P}_{0}\|\!
\left(\frac{\|{\bf B}(\delta)\|}{\mu_{\min}}\right)^{p-2}
\\
&=
\|{\bf S}_0{\bf B}(\delta)\|\,\|{\bf S}_0{\bf B}(\delta){\bf P}_{0}\|\,\sum_{p=2}^\infty\,\left({2p \choose p}-2\right)
\left(\frac{\|{\bf B}(\delta)\|}{\mu_{\min}}\right)^{p-2}
\\
&<
\|{\bf S}_0{\bf B}(\delta)\|\,\|{\bf S}_0{\bf B}(\delta){\bf P}_{0}\| \sum_{p=2}^\infty\,{2p \choose p}
\left(\frac{\|{\bf B}(\delta)\|}{\mu_{\min}}\right)^{p-2}
\\
&\leq 16C\
\frac{\|{\bf S}_0{\bf B}(\delta)\|\,\|{\bf S}_0{\bf B}(\delta){\bf P}_{0}\|}{1-4\|{\bf B}(\delta)\|/\mu_{\min}}\ .
\end{split}
\end{gather*}
The proof is complete.

\vspace{3mm}
 {\bf Proof of Proposition \ref{prop:L_delta}}
\vspace{1mm}

\noindent
First of all we have
\begin{gather*}
\left({\bf I}-\delta^2{\bf A}^{(2)}_0/\mu\right)^{-1}-{\bf I}
=
\left({\bf I}-\delta^2{\bf A}^{(2)}_0/\mu\right)^{-1} {\delta^2{\bf A}^{(2)}_0}/{\mu}\ .
\end{gather*}
Therefore, we get that
\begin{gather*}
\begin{split}
&{\bf L}(\delta)=\sum_{\mu>0}\big({\bf P}_\mu {\bf B}(\delta){\bf P}_0+
{\bf P}_0 {\bf B}(\delta){\bf P}_\mu\big)\big/\mu
\\
&+
\sum_{\mu>0}\frac{{\bf P}_\mu {\bf B}(\delta){\bf P}_0}{\mu}
\left(\left({\bf I}-\delta^2{\bf A}^{(2)}_0/{\mu}\right)^{-1}-{\bf I}\right)
\\
&+
\sum_{\mu>0}\left(\left({\bf I}-\delta^2{\bf A}^{(2)}_0/{\mu}\right)^{-1}-{\bf I}\right)
\frac{{\bf P}_0 {\bf B}(\delta){\bf P}_\mu}{\mu}
\\
&=\delta\,{\bf S}_0\big({\bf H}{\bf E}^{\rm T}+{\bf E}{\bf H}^{\rm T}\big){\bf P}_0+
\delta\,{\bf P}_0\big({\bf H}{\bf E}^{\rm T}+{\bf E}{\bf H}^{\rm T}\big){\bf S}_0
\\
&+
\delta^2\big({\bf S}_0{\bf A}^{(2)}{\bf P}_0+{\bf P}_0{\bf A}^{(2)}{\bf S}_0\big)
\\
&+\sum_{\mu>0}\frac{{\bf P}_\mu {\bf B}(\delta){\bf P}_0}{\mu}\,
\frac{\delta^2{\bf A}^{(2)}_0}{\mu}\left({\bf I}-\delta^2{\bf A}^{(2)}_0/\mu\right)^{-1}
\\
&+\sum_{\mu>0}
\left({\bf I}-\delta^2{\bf A}^{(2)}_0/\mu\right)^{-1}\frac{\delta^2{\bf A}^{(2)}_0}{\mu}
\frac{{\bf P}_0 {\bf B}(\delta){\bf P}_\mu}{\mu}\ .
\end{split}
\end{gather*}
Since ${\bf P}_0{\bf H}={\bf 0}$ and in view of the equality
\eqref{eq:first_term}, we get the result.

\vspace{3mm}
{\bf Proof of Theorem \ref{theor:precise_proj_noise}}

\vspace{1mm}
\noindent
As it is shown in the demonstration of Theorem \ref{theor:precise_proj},
\begin{gather*}
\begin{split}
&{\bf P}_0^{\perp}(\delta)-{\bf P}_0^{\perp}- {\bf L}(\delta)
\\
&=
\sum_{p=2}^\infty (-1)^{p}\ {\sum}^*_{(l_1,\ldots,l_{p+1})}\
{{\bf S}}_{0}^{(l_1)}{{\bf B}}(\delta){{\bf S}}_{0}^{(l_2)}\ldots\,{{\bf B}}(\delta){{\bf S}}_{0}^{(l_{p+1})},
\end{split}
\end{gather*}
where ${\sum}^*_{(l_1,\ldots,l_{p+1})}$ stands for the sum over all positive integers $l_1,\ldots,l_{p+1}$ such that $l_1+\ldots+l_{p+1}=p$
and $l_1,l_{p+1}\neq p$.

Let us calculate the sum
\begin{gather*}
\sum_{p=2}^{\infty}(-1)^{p}\ {\sum}^{**}_{(l_1,\ldots,l_{p+1})}\
{{\bf S}}_{0}^{(l_1)}{{\bf B}}(\delta){{\bf S}}_{0}^{(l_2)}\ldots\,{{\bf B}}(\delta){{\bf S}}_{0}^{(l_{p+1})},
\end{gather*}
where ${\sum}^{**}_{(l_1,\ldots,l_{p+1})}$ stands for the sum over all positive integers $l_1,\ldots,l_{p+1}$ with $l_1+\ldots+l_{p+1}=p$
such that for some $2\leq k\leq p$ either  the numbers $l_1,\ldots, l_k$ are positive and $l_{k+1}=\ldots=l_{p+1}=0$
or $l_1=\ldots=l_{p+1-k}=0$ and $l_{p+2-k},\ldots, l_{p+1}>0$. Both variants are treated similarly.

Denote for brevity
$
{\bf M}_j=\left({\delta^2{\bf A}_0^{(2)}}\big/{\mu_j}\right)^{l_j-1}.
$
If the numbers $l_1,\ldots, l_k$ are positive and $l_j=0$ for $j>k$, then
\begin{gather*}
\begin{split}
&
(-1)^p\,{{\bf S}}_{0}^{(l_1)}{{\bf B}}(\delta){{\bf S}}_{0}^{(l_2)}\ldots\,{{\bf B}}(\delta){{\bf S}}_{0}^{(l_{p+1})}
\\
&
=(-1)^{k-1}\!\left(\prod_{j=1}^{k-1} {\bf S}_0^{l_j}{\bf B}(\delta)\right)\!
{\bf S}_0^{l_k}{\bf B}(\delta){\bf P}_0
\left(\delta^2{\bf A}_0^{(2)}\right)^{p-k}
\\
&
=\!(-1)^{k-1}\!\!\left(\prod_{j=1}^{k-1}\! \sum_{\mu_j>0}\frac{{\bf P}_{\mu_j}}{\mu_j^{l_j}}{\bf B}(\delta)\!\right)\!\!
\sum_{\mu_k}\frac{{\bf P}_{\mu_k}{\bf B}(\delta){\bf P}_0}{\mu_k^{l_k}}\!
\left(\delta^2{\bf A}_0^{(2)}\right)^{\!p-k}
\\
&
=(-1)^{k-1}\sum_{\mu_1,\ldots,\mu_k>0} {\bf J}_k
\frac{{\bf P}_{\mu_k}{\bf B}(\delta){\bf P}_0}{\mu_k}\,
\frac{\left(\delta^2{\bf A}_0^{(2)}\right)^{p-k}}{\mu_1^{l_1-1}\ldots \mu_k^{l_k-1}}
\\
&
=(-1)^{k-1}\sum_{\mu_1,\ldots,\mu_k>0} {\bf J}_k
\frac{{\bf P}_{\mu_k}{\bf B}(\delta){\bf P}_0}{\mu_k}\ {\bf M}_k
\prod_{j=1}^{k-1} {\bf M}_j.
\end{split}
\end{gather*}
For fixed $k\geq 2$ we have
\begin{gather}
\begin{split}
&
\sum_{p\geq k}\ \sum_{\substack{l_1+\dots+l_k=p, \ l_j>0\\
\mu_1,\ldots,\mu_k>0}} {\bf J}_k
\frac{{\bf P}_{\mu_k}{\bf B}(\delta){\bf P}_0}{\mu_k}\  {\bf M}_k
\prod_{j=1}^{k-1}  {\bf M}_j
\\
&=\sum_{l_1,\ldots,l_k=1}^{\infty}
\sum_{\mu_1,\ldots,\mu_k>0} {\bf J}_k
\frac{{\bf P}_{\mu_k}{\bf B}(\delta){\bf P}_0}{\mu_k} \  {\bf M}_k
\prod_{j=1}^{k-1} {\bf M}_j
\\
&=\sum_{\mu_1,\ldots,\mu_k>0} {\bf J}_k
\frac{{\bf P}_{\mu_k}{\bf B}(\delta){\bf P}_0}{\mu_k}
\sum_{l_k=1}^{\infty}  {\bf M}_k
\prod_{j=1}^{k-1}\sum_{l_j=1}^{\infty} {\bf M}_j
\\
&=\sum_{\mu_1,\ldots,\mu_{k-1}>0} {\bf J}_k
\sum_{\mu_k}\frac{{\bf P}_{\mu_k}{\bf B}(\delta){\bf P}_0}{\mu_k}
\left({\bf I}-{\bf M}_k\right)^{-1} {\bf G}_k
\\
&=\sum_{\mu_1,\ldots,\mu_{k-1}>0} {\bf J}_k
{\bf L}_1(\delta) {\bf G}_k.
\end{split}
\label{eq:term_T1}
\end{gather}
Since the case $k=1$ corresponds to ${\bf L}_1(\delta)$, then the sum of \eqref{eq:term_T1} over $k\geq 1$ gives
us ${\bf T}_1$. The term ${\bf T}^{\rm T}_1$ corresponds to $l_1=\ldots=l_{p+1-k}=0$ and $l_{p+2-k},\ldots, l_{p+1}>0$
and is calculated in the same manner.

Let us consider the residual ${\bf P}_0^{\perp}(\delta)\!-\!{\bf P}_0^{\perp}\!-\!{\bf T}(\delta)$.
As it was mentioned in the demonstration of Theorem \ref{theor:zero_pert},
for any $l_1,\ldots,l_{p+1}$ with $l_1+\ldots+l_{p+1}=p$ there exists at least one $j$ such that $l_j>0$ and
either $l_{j-1}=0$ or $l_{j+1}=0$. It is easy to understand that  such $j$  is unique iff either $l_i>0$ for
$i<j$ and $l_i=0$ for $i>j$ or $l_i=0$ for $i<j$ and $l_i>0$ for $i>j$.

This means that  any term
\begin{gather*}
(-1)^p\,{{\bf S}}_{0}^{(l_1)}{{\bf B}}(\delta){{\bf S}}_{0}^{(l_2)}\ldots\,{{\bf B}}(\delta){{\bf S}}_{0}^{(l_{p+1})}
\end{gather*}
of the residual ${\bf P}_0^{\perp}(\delta)-{\bf P}_0^{\perp}-~{\bf T}(\delta)$
corresponds to the vector $(l_1,\ldots,l_{p+1})$ such that for some positive  $l_i, l_j$ with $i<j$  each of the triples
$(l_{i-1},l_i,l_{i+1})$, $(l_{j-1},l_j,l_{j+1})$ contains a unique zero.
Denote
\begin{gather*}
m=
\begin{cases}
i& \mbox{for} \  l_{i+1}=0,\\
i-1& \mbox{for}  \ l_{i-1}=0
\end{cases}
 \ \mbox{and} \ \
k=
\begin{cases}
j& \mbox{for} \  l_{j+1}=0,\\
j-1& \mbox{for}  \ l_{j-1}=0
\end{cases}
\end{gather*}
and consider the pairs $(m,m+1)$ and $(k,k+1)$. Note that these pairs are either  disjoint or $m+1=i+1=k=j-1$
with $l_{m+1}=l_k=0$.
If pairs are disjoint, then
\begin{gather*}
\begin{split}
&\big\|{{\bf S}}_{0}^{(l_1)}{{\bf B}}(\delta){{\bf S}}_{0}^{(l_2)}\!\ldots\,{{\bf B}}(\delta){{\bf S}}_{0}^{(l_{p+1})}\big\|
\\
&=
\big\|{{\bf S}}_{0}^{(l_1)}\!\ldots {\bf S}_0^{(l_m)}{{\bf B}}(\delta){\bf S}_0^{(l_{m+1)}}\ldots
{\bf S}_0^{(l_k)}{{\bf B}}(\delta){\bf S}_0^{(l_{k+1})}\ldots
{{\bf S}}_{0}^{(l_{p+1})}\big\|
\\
&\leq
\big\|{\bf S}_0{{\bf B}}(\delta){\bf P}_0\big\|^2\,\left({\|{\bf B}(\delta)\|}\big/{\mu_{\min}}\right)^{p-2}.
\end{split}
\end{gather*}
If $m+1=k$, then $l_m,l_{m+2}>0$ and
\begin{gather*}
\begin{split}
&\big\|{{\bf S}}_{0}^{(l_1)}{{\bf B}}(\delta){{\bf S}}_{0}^{(l_2)}\ldots\,{{\bf B}}(\delta){{\bf S}}_{0}^{(l_{p+1})}\big\|
\\
&=\!\big\|{{\bf S}}_{0}^{(l_1)}{{\bf B}}(\delta)\!\ldots
{\bf S}_0^{(l_m)}{{\bf B}}(\delta){\bf P}_0{\bf P}_0{\bf B}(\delta){\bf S}_0^{(l_{m+2})}\!\ldots
{{\bf B}}(\delta){{\bf S}}_{0}^{(l_{p+1})}\big\|
\\
&\leq\!
\big\|{\bf S}_0{{\bf B}}(\delta){\bf P}_0\big\|^2\,\left({\|{\bf B}(\delta)\|}\big/{\mu_{\min}}\right)^{p-2}.
\end{split}
\end{gather*}
The rest of demonstration is the same as in Theorem \ref{theor:precise_proj}.

\vspace{3mm}
{\bf Proof of Proposition \ref{prop:big_NSR}}

\vspace{1mm}
\noindent
If $\Theta\rightarrow 0$, then $\Theta_1\rightarrow 0$ and
\begin{gather*}
\big\|{\bf B}(\delta)\big\|/\mu_{\min}\leq 2|\delta|\Theta+\delta^2\, \Theta \Theta_1 \sim 2|\delta| \Theta\rightarrow 0
\end{gather*}
for all $\delta$.
This means that we can use inequalities \eqref{eq:coarse_proj} and \eqref{eq:res_L_delta}
for any fixed $\delta$ provided that $N$ is sufficiently big. The inequality  \eqref{eq:Diff_proj} follows from \eqref{eq:coarse_proj}, \eqref{eq:B/mu_coarse}.
In the same manner inequalities \eqref{eq:coarse_proj_mainterm}, \eqref{eq:B/mu_coarse}  provide that
\begin{gather*}
\limsup_{N}\,\Theta^{-2}\,\big\|{\bf P}_0^{\perp}(\delta)-{\bf P}_0^{\perp}-{\bf W}_1(\delta)\big\|\leq
64 C \delta^2.
\end{gather*}
Since ${\bf W}_1(\delta)=
\delta\,{\bf V}^{(1)}_0+\delta^2\big({\bf P}_{0}{\bf A}^{(2)}{\bf S}_0+{\bf S}_0{\bf A}^{(2)}{\bf P}_{0}\big)$
and $\|{\bf S}_0 {\bf A}^{(2)} {\bf P}_0\|=\big\|{\bf P}_{0}{\bf A}^{(2)}{\bf S}_0\big\|
\leq \Theta \Theta_1\leq\Theta^2$, then \eqref{eq:Diff_proj_second} is  proved as well.

\vspace{3mm}
{\bf Proof of Proposition \ref{prop:asymp_ort}}

\vspace{1mm}
\noindent
It follows from the convergence $\big\|{\bf H}{\bf E}^{\rm T}\big\|/\mu_{\min}\rightarrow 0$ that
\begin{gather*}
\begin{split}
&\limsup \|{\bf B}(\delta)\|/\mu_{\min}\\
&\!\leq\! 2|\delta|\limsup \big\|{\bf H}{\bf E}^{\rm T}\big\|/\mu_{\min}\!+\!
\delta^2 \limsup \Theta \Theta_1
\!=\!\delta^2/\Delta\!<\!1/4.
\end{split}
\end{gather*}
 Therefore, \eqref{eq:coarse_orth} holds for
$|\delta|<\delta_0$ and big $N$.
Moreover,
$1/(1-4\|{\bf B}(\delta)\|/\mu_{\min})\leq 1/(1-4\Delta)$. Thus
\eqref{eq:asymp_ort} follows from \eqref{eq:coarse_orth} and \eqref{eq:norm_SB}.
Applying \eqref{eq:res_L_delta} instead of \eqref{eq:coarse_orth}, we get \eqref{eq:asymp_ort_L}.

\vspace{3mm}
{\bf Proof of Proposition \ref{prop:sum_exp}}

\vspace{1mm}
\noindent
In view of Lemma \ref{summ:series},  $\mu_{\max}\sim T_{\max}^{(a)} a_1^{2N}$ and $\mu_{\min}\sim
T_{\min}^{(a)}\,a_{p}^{2N}$.
Therefore,
\begin{gather*}
\begin{split}
\Theta_1&=\sqrt{\nu_{\max}/\mu_{\max}}\sim \sqrt{\nu_{\max}/T_{\max}^{(a)}}\,a_{1}^{-N},\\
\Theta_2&=\mu_{\max}/\mu_{\min}\sim \big(T_{\max}^{(a)}/T_{\min}^{(a)}\big)\,(a_1/a_p)^{2N},\\
\Theta&= \Theta_1 \Theta_2\sim \sqrt{\nu_{\max}},
({a_1}/{a_p^2})^N \Xi_a
\end{split}
\end{gather*}
and both assertions follow from Proposition \ref{prop:big_NSR}.

\vspace{3mm}
{\bf Proof of Proposition \ref{prop:gen_wnoise_LeqK}}

\vspace{1mm}
\noindent
1. In view of Lemma \ref{summ:series_rand} and Remark \ref{rem:gamma_prim},
$\nu_{\max}/N\ln{N}<\gamma'S$  for $N>N_0(\gamma',\omega)$, $\gamma'>\gamma_0$,
and $\omega\in \Omega_0$. This means that
\begin{gather*}
\begin{split}
\Theta = \sqrt{\frac{\nu_{\max}}{\mu_{\max}}}\, \frac{\mu_{\max}}{\mu_{\min}}
=\sqrt{\frac{\nu_{\max}}{N\ln{N}}}\,
\frac{\sqrt{\mu_{\max}{N\ln N}}}{\mu_{\min}}
\leq \sqrt{\gamma' S} {\sqrt{\mu_{\max}{N\ln N}}}/{\mu_{\min}}
\rightarrow 0
\end{split}
\end{gather*}
on $\Omega_0$ as $N\rightarrow \infty$.   The rest follows from Proposition \ref{prop:big_NSR}.

2. Since $\nu_{\max}/N\rightarrow \sigma_{\max}$ with probability 1, see Lemma \ref{summ:series_rand}, and
 $\mu_{\max}/\mu^2_{\min}=o\big(N^{-1}\big)$, then
\begin{gather*}
\begin{split}
\Theta=\sqrt{\frac{\nu_{\max}}{\mu_{\max}}}\,\frac{\mu_{\max}}{\mu_{\min}}
=\sqrt{\frac{\nu_{\max}}{N}}\,\frac{\sqrt{N\mu_{\max}}}{\mu_{\min}}
\sim \sqrt{\sigma_{\max}} \,\frac{\sqrt{N\mu_{\max}}}{\mu_{\min}}
\rightarrow 0 \ \ a.s.
\end{split}
\end{gather*}
and all we need is to apply Proposition \ref{prop:big_NSR}.

\vspace{3mm}
{\bf Proof of Proposition \ref{prop:gen_wight_noise_coarse}}

\vspace{1mm}
\noindent
To demonstrate this assertion, we use Theorems \ref{theor:gen_upper} and \ref{theor:precise_proj_noise}
as well as the results of Lemma \ref{summ:series_rand} and Appendix 2.

Since
$\big\|{\bf B}(\delta)\big\|\leq 2|\delta|\,\big\|{\bf H}{\bf E}^{\rm T}\big\|+\delta^2\nu_{\max}$
and
\begin{gather*}
\frac{2\delta_0\,\big\|{\bf H}{\bf E}^{\rm T}\big\|+\delta_0^2\nu_{\max}}{\mu_{\min}}<1/4,
\end{gather*}
then the inequality \eqref{eq:gen_upper} holds
 for any $\delta\in (-\delta_0,\delta_0)$.
It follows from Proposition \ref{prop:LIL_osc_wn} of Appendix~2 that almost surely
\begin{gather}
\label{eq:limsup_HE}
\limsup_{N}\,(N\ln\ln N)^{-1/2}\,\big\|{\bf H}{\bf E}^{\rm T}\big\|\leq c_1
\end{gather}
with some positive constant $c_1$. Therefore, $\big\|{\bf H}{\bf E}^{\rm T}\big\|/\mu_{\min}\!\rightarrow \!0$ with
probability~1.

As it follows from Lemma \ref{summ:series_rand},
$\nu_{\max}/N\!\rightarrow\!1$ with probability 1. Therefore, there exist
$\Omega_1\in {\cal F}$ with $\mathbb{P}(\Omega_1)=1$ such that for any
$\omega\in \Omega_1$
\begin{gather*}
\begin{split}
\limsup_N \big\|{\bf B}(\delta)\big\|/\mu_{\min}
\leq \limsup_N\left(2\delta \frac{\|{\bf H}{\bf E}^{\rm T}\|}{\mu_{\min}}+\delta^2 \frac{\nu_{\max}}{\mu_{\min}}\right)
= \delta^2/(L_0\Lambda_{\min}).
\end{split}
\end{gather*}
Thus  \eqref{eq:gen_upper} is proved.
Consider the  right-hand side
\begin{gather}
\label{eq:res_gen}
\frac{\|{\bf S}_0{\bf B}(\delta){\bf P}_{0}\|}{1-4\|{\bf B}(\delta)\|/\mu_{\min}}
\end{gather}
of the inequality \eqref{eq:gen_upper}. We have already proved that
$1-4\|{\bf B}(\delta)\|/\mu_{\min}\geq 1-4\delta_0^2/(\Lambda_{\min}L_0)>0$
for $|\delta|<\delta_0$, $\omega\in \Omega_1$ and $N>N_0(\omega)$.)

Consider the numerator of the fraction \eqref{eq:res_gen}. It is easy to see that
\begin{gather*}
{\bf S}_0{\bf B}(\delta){\bf P}_0=\delta\,{\bf S}_0{\bf H}{\bf E}^{\rm T}{\bf P}_0
+\delta^2\, \sum_{\mu>0}\frac{{{\bf P}_\mu\,\bf A}^{(2)}{\bf P}_0}{\mu}\ .
\end{gather*}
In view of \eqref{eq:limsup_HE},
\begin{gather*}
\begin{split}
&\limsup_N\,\sqrt{N/\ln \ln N}\,\big\|{\bf S}_0{\bf H}{\bf E}^{\rm T}{\bf P}_0\big\|\\
&\leq
\limsup_N\,\sqrt{N/\ln \ln N}\,\big\|{\bf H}{\bf E}^{\rm T}\big\|/\mu_{\min}\leq c_1/(\Lambda_{\min}L_0)
\end{split}
\end{gather*}
for $\omega\in \Omega_1$.
On the other hand,
since ${\bf P}_{\mu}{\bf P}_{0}={\bf 0}$, then it follows from Proposition \ref{prop:white_noise} of Appendix 2 that
\begin{gather*}
\begin{split}
&\limsup_N\,\sqrt{N/\ln \ln N}\ \frac{\big\|{\bf P}_\mu\,{\bf A}^{(2)}{\bf P}_0\big\|}{\mu}\\
&\leq
\limsup_N\,\sqrt{N/\ln \ln N}\ \frac{\big\|{\bf P}_\mu\,\big({\bf A}^{(2)}-N\,{\bf I}\big){\bf P}_0\big\|}{\mu_{\min}}
\\
&=
\limsup_N\,\sqrt{N/\ln \ln N}\ \frac{\big\|{\bf A}^{(2)}-N{\bf I}\big\|}{N}\,\frac{1}{\Lambda_{\min}L_0}
\\
&=
\frac{1}{\Lambda_{\min}L_0}\ \limsup_N\,\sqrt{N}\ \frac{\big\|{\bf A}^{(2)}/N-{\bf I}\big\|}{\sqrt{\ln \ln N}}
\leq c_1
\end{split}
\end{gather*}
with probability 1.
Note that the number of positive eigenvalues of the
matrix ${\bf H}{\bf H}^{\rm T}$ does not exceed $2p$, thus \eqref{eq:gen_wight_noise_coarse1} is proved. To prove
\eqref{eq:gen_wight_noise_coarse2} it is sufficient to use \eqref{eq:res_T_delta} instead of \eqref{eq:gen_upper}.

\vspace{3mm}
{\bf Proof of Proposition \ref{prop:exp_const_proj_res}}

\vspace{1mm}
\noindent
We present only the sketch of the proof. Note that the right-hand side of \eqref{eq:sec_or_exp} is $O\big(N^2a^{-2N}\big)$ if $L/N\rightarrow \alpha \in
(0,1)$
and has the order $O\big(Na^{-2N}\big)$ in the case $\min(L,K)={\rm const}>1$.
Let us study the asymptotic behavior
of ${\bf V}_0^{(1)}$.

Accurate calculations show that
\begin{gather*}
\begin{split}
&\big\|{\bf V}_0^{(1)}\big\|\sim\big\|{\bf Z}_0^{(1)}\big\|\\
&\sim
\frac{\sqrt N}{a^{N+1}}
\begin{cases}
\displaystyle
\alpha(a+1)\sqrt{a^2-1}& \mbox{for}\ \ L/N\rightarrow \alpha\in (0,1),\\
\displaystyle
\frac{(a+1)\sqrt{a^2-1}}{1-a^{-K_0}}& \mbox{for} \ \ K=K_0>1
\end{cases}
\end{split}
\end{gather*}
and
\begin{gather*}
\begin{split}
\big\|{\bf V}_0^{(1)}\big\|
\sim \frac{(a+1)}{a}\,\frac{a^{L_0}\sqrt{L_{0}\|W_{L_0}\|^2-\beta_{L_0}^2}}{\|W_{L_0}\|^2}\ a^{-N}
\end{split}
\end{gather*}
for $L=L_0={\rm const}>1$. This completes the proof.

\vspace{3mm}
{\bf Proof of Proposition \ref{prop:ort_const_saw}}

\vspace{1mm}
\noindent
Simple calculations show that
 $\|{\bf H}{\bf E}^{\rm T}\big\|/\mu_{\min}=0$
 if $K$ is even and that $\|{\bf H}{\bf E}^{\rm T}\big\|/\mu_{\min}=1/K$ for odd $K$. Additionally,
 $\|{{\bf S}}_{0}{\bf A}^{(2)}\big\|=0$ for even $L$ and $\|{{\bf S}}_{0}{\bf A}^{(2)}\big\|=1/L$
for odd $L$.

The assertion of
Theorem \ref{theor:precise_proj} holds under the condition
$2|\delta|\,\|{\bf H}{\bf E}^{\rm T}\|+\delta^2{\nu_{\max}}<\mu_{\min}/4$.
If $K$ is even, then
this inequality is equivalent to $\delta^2<1/4$. If $K$ is odd and tends to
infinity, then we asymptotically get the same restriction  $\delta^2<1/4$.

Since the right-hand side of \eqref{eq:asymp_ort} has the form
$O\!\left(\big\|{\bf H}{\bf E}^{\rm T}\big\|/\mu_{\min}+\big\|{{\bf S}}_{0}{\bf A}^{(2)}\big\|\right)$,
the proof is complete.

\vspace{3mm}
{\bf Proof of Proposition \ref{prop_const_saw}}

\vspace{1mm}
\noindent
The first assertion follows from the equality $\|W_E\|=\sqrt{L^2-\beta_L}$.

Let us prove the second assertion.
Straightforward calculations show that
\begin{gather*}
\|{\bf S}_0{\bf B}(\delta)\|\,\|{\bf S}_0{\bf B}(\delta){\bf P}_{0}\|\leq
\delta^2\left(\frac{\beta_K}{K}+\frac{\beta_L\beta_K}{LK}+|\delta|\,\frac{\beta_L}{L}\right)^2
\end{gather*}
and
\begin{gather*}
\big\|{\bf B}(\delta)\big\|/\mu_{\min}\leq |\delta|\,\beta_k\sqrt{L^2-\beta_L}/LK+\delta^2
\leq \beta_k\,|\delta|/K+\delta^2.
\end{gather*}
Therefore, it follows from \eqref{eq:res_L_delta} that for any $|\delta|<1/2$
\begin{gather*}
 \begin{split}
&\big\|{\bf P}_0^{\perp}(\delta)-{\bf P}_0^{\perp}-{\bf L}(\delta)\big\|\\
&\!=\!
\begin{cases}
\!
O\left(K^{-2}\right)
\!&\!\mbox{for odd} \,  K\rightarrow \infty \,  \mbox{and even} \,  L, \\
O\left(L^{-2}\right)
\!&\!\mbox{for odd} \,  L\rightarrow \infty \,  \mbox{and even} \,  K, \\
O\big(L^{-2}\!+\!K^{-2}\big)\!
\!&
\begin{aligned}
&\!\mbox{for odd} \, K,L\\
&\!\mbox{and}\,   \min(L,K)\!\rightarrow\! \infty.
\end{aligned}
\end{cases}
\end{split}
\end{gather*}

Comparing ${\bf L}(\delta)$ with ${\bf M}(\delta)$ we see that ${\bf L}(\delta)={\bf M}(\delta)$ in the case when
$L$ is even and $K$ is odd. Otherwise,
\begin{gather*}
\begin{split}
&{\bf L}(\delta)={\bf M}(\delta)\\
&+
 \begin{cases}
 O(L^{-3})&\mbox{for odd} \,  L\rightarrow \infty \,  \mbox{and even} \,  K, \\
 O(L^{-2}K^{-1}+L^{-3})&
 \begin{aligned}
 &\mbox{for odd} \, K,L\\
&\mbox{and}\,   \min(L,K)\!\rightarrow\! \infty.
 \end{aligned}
 \end{cases}
 \end{split}
 \end{gather*}
This completes the proof.

\vspace{3mm}
{\bf Proof of Proposition \ref{prop:const_wnoise_main_term_L_0}}

\vspace{1mm}
\noindent
It follows  from conditions of Theorem \ref{theor:precise_proj_noise}
that
$\mathbb{P}\big(\Omega_N\big)\geq \mathbb{P}\big(|\delta|\, C(\delta)<1/4\big),$
where
\begin{gather*}
C(\delta)=
 2\,\frac{\|{\bf H}{\bf E}^{\rm T}\|}{L_0K}+|\delta|\, \frac{\|{\bf A}\|}{L_0K}\,.
\end{gather*}
Since almost surely $\|{\bf H}{\bf E}^{\rm T}\|/K\!\rightarrow\! 0$ and $\|{\bf A}\|/K\!\sim\! \|{\bf A}\|/N
\!\rightarrow\! 1$ as $N\rightarrow \infty$,
then $\mathbb{P}\big(\Omega_N\big)\rightarrow 1$.

Now let us turn to the operator ${\bf T}(\delta)$, see  \eqref{eq:res_T_delta}.
First of all,
${\bf P}_{\mu}={\bf P}_0^{\perp}=W_{L_0}W_{L_0}^{\rm T}/L_0$. Thus this projector does not depend on $N$.
We get from  \eqref{eq:Ldelta_short} and \eqref{eq:T_1}
that
\begin{gather}
\label{eq:lim_T_wn}
\sqrt{K}\,{\bf T}_1(\delta)\!=\! \sum_{i=0}^{\infty}(-1)^{i}
\left(\frac{{\bf P}_{\mu} {\bf B}(\delta)}{\mu}\right)^i
\left(\sqrt{K} \frac{{\bf P}_\mu {\bf B}(\delta){\bf P}_0}{\mu}\right)
{\bf I}_i(\delta),
\end{gather}
where
\begin{gather}
\begin{split}
&\frac{{\bf P}_\mu {\bf B}(\delta){\bf P}_0^\perp}{\mu}\!=\!\delta {\bf P}_0^\perp
\frac{{\bf H}{\bf E}^{\rm T}\!+\!{\bf E}{\bf H}^{\rm T}}{L_0K} {\bf P}_0
\!+\!\delta^2 {\bf P}_0^\perp \frac{{\bf A}^{(2)}}{L_0K} {\bf P}_0\\
&=
\delta\,
\frac{{\bf H}{\bf E}^{\rm T}}{L_0K}\,{\bf P}_0+\delta^2 {\bf P}_0^\perp\frac{{\bf A}^{(2)}}{L_0K}\,{\bf P}_0
=\frac{\delta}{L_0}\left(\frac{{\bf H}{\bf E}^{\rm T}}{K}\,{\bf P}_0+
\delta{\bf P}_0^\perp\frac{{\bf A}^{(2)}}{K}\,{\bf P}_0\right)
\label{eq:Pmu_term}
\end{split}
\end{gather}
and
\begin{gather*}
{\bf I}_i(\delta)=\left({\bf I}-\delta^2\frac{{\bf A}^{(2)}_0}{L_0K}\right)^{-(i+1)}.
\end{gather*}

Let us consider the operator $\sqrt{K}{\bf P}_\mu {\bf B}(\delta){\bf P}_0/{\mu}$,  see the right-hand side of
\eqref{eq:lim_T_wn}.
Since ${\bf H}{\bf E}^{\rm T}=\big\{\psi_{ij}\big\}_{i,j=0}^{L-1}$ with $\psi_{ij}=
\varepsilon_j+\ldots+\varepsilon_{K-j-1}$ and since
\begin{gather*}
\left(\varepsilon_j+\ldots+\varepsilon_{K-j-1}\right)/\sqrt{K}\!=\!
\left(\varepsilon_0+\ldots+\varepsilon_{K-1}\right)/\sqrt{K}+D(K)
\end{gather*}
with $D(K)\rightarrow 0$ in probability, we obtain that
\begin{gather*}
{\cal L}\left(\frac{1}{\sqrt N}\,{\bf H}{\bf E}^{\rm T}\right)\Longrightarrow
{\cal L}\big(L_0\,\xi\,{\bf P}_\mu\big),
\end{gather*}
where $\xi\in {\rm N}(0,1)$.
In view of the equality $L_0\,\xi\,{\bf P}_\mu{\bf P}_0={\bf 0}$ this means that
\begin{gather*}
\sqrt{K}\,\frac{\delta}{L_0}\,\frac{{\bf H}{\bf E}^{\rm T}}{K}\,{\bf P}_0
\stackrel{\mathbb{P}}\longrightarrow
{\bf 0}.
\end{gather*}

Consider the second term of the right-hand side of \eqref{eq:Pmu_term}.
Since $\mathbb{E}{\bf A}^{(2)}/K={\bf I}$
and
in view of the third assertion of Proposition \ref{prop:white_noise} (see Appendix 2), we have
\begin{gather*}
\begin{split}
&
{\cal L}\left(\sqrt{K}\delta^2 {\bf P}_0^\perp\frac{{\bf A}^{(2)}}{K L_0} {\bf P}_0\right)\!=\!
{\cal L}\left(\frac{\delta^2}{L_0} {\bf P}_0^\perp\sqrt{K}\,\frac{{\bf A}^{(2)}
\!-\!\mathbb{E}{\bf A}^{(2)}}{K}\,{\bf P}_0\right)
\\
&
\Longrightarrow {\cal L}\left(\frac{\delta^2}{L_0}\,{\bf P}_0^\perp {\Psi}_{L_0}\,{\bf P}_0\right).
\end{split}
\end{gather*}
Therefore,
\begin{gather}
\label{eq:PmuGP_0_CLT}
{\cal L}\big(\sqrt{K}{\bf P}_0^\perp {\bf B}(\delta){\bf P}_0/{L_0K}\big)
\Longrightarrow {\cal L}\left(\frac{\delta^2}{L_0}\,{\bf P}_0^\perp {\Psi}_{L_0}\,{\bf P}_0\right).
\end{gather}

Since
\begin{gather*}
\frac{{\bf P}_\mu {\bf B}(\delta)}{\mu}=\delta\, {\bf P}_0^\perp
\frac{{\bf H}{\bf E}^{\rm T}+{\bf E}{\bf H}^{\rm T}}{L_0K}
+\delta^2 {\bf P}_0^\perp\frac{{\bf A}^{(2)}}{L_0K}\rightarrow \frac{\delta^2}{L_0}\, {\bf P}_0^\perp \ \ \mbox{a.s.},
\end{gather*}
we see that
$\big({\bf P}_\mu {\bf B}(\delta)/\mu\big)^i\!\rightarrow\! \big(\delta^2/L_0\big)^i {\bf P}_0^\perp$ with probability 1.

In view of Proposition \ref{prop:white_noise},
${\bf I}-\delta^2{{\bf A}^{(2)}_0}/{\mu}\rightarrow {\bf I}-{\delta^2}{\bf P}_0 /{L_0}$ almost surely.
Therefore, with probability~1
\begin{gather*}
\begin{split}
&\left({\bf I}-\delta^2\frac{{\bf A}^{(2)}_0}{\mu}\right)^{-1}\rightarrow
\left({\bf I}-\frac{\delta^2}{L_0}\, {\bf P}_0\right)^{-1}\\
&={\bf I}+
\sum_{k=1}^\infty \left(\frac{\delta^2}{L_0}\right)^k{\bf P}_0^{k}=
{\bf I}+\frac{\delta^2}{L_0}\,\frac{1}{1-\delta^2/L_0}\,{\bf P}_0
\\
&= {\bf I} - {\bf P}_0 + \frac{1}{1-\delta^2/L_0}\,{\bf P}_0={\bf P}_0^{\perp}+\frac{1}{1-\delta^2/L_0}\,{\bf P}_0.
\end{split}
\end{gather*}
and
\begin{gather*}
\left({\bf I}-\delta^2\frac{{\bf A}^{(2)}_0}{\mu}\right)^{-i}\rightarrow
{\bf P}_0^\perp+\left(\frac{1}{1-\delta^2/L_0}\right)^i\,{\bf P}_0 \ \ \mbox{a.s.}
\end{gather*}
Collecting all these results, we see that $\sqrt{N}\,{\bf T}_1(\delta)$ converges in distribution to the random
matrix
\begin{gather*}
\begin{split}
&\sum_{i=0}^{\infty}(-1)^{i}\!
\left(\!\frac{\delta^2}{L_0}\!\right)^i\! {\bf P}_0^\perp
\frac{\delta^2}{L_0}{\bf P}_0^\perp {\Psi}_{L_0} {\bf P}_0
\!\left(\!
{\bf P}_0^\perp\!+\!\left(\!\frac{1}{1\!-\!\delta^2/L_0}\!\right)^{i+1}\!{\bf P}_0\!
\right)
\\
&=
\sum_{i=0}^{\infty}(-1)^{i}
\left(\frac{\delta^2}{L_0}\right)^{i+1}
\left(\frac{1}{1-\delta^2/L_0}\right)^{i+1}
{\bf P}_0^\perp {\Psi}_{L_0}\,{\bf P}_0\\
&=
\frac{\delta^2/{L_0}}{1-\delta^2/L_0} \big(1-\delta^2/L_0\big) {\bf P}_0^\perp {\Psi}_{L_0} {\bf P}_0
=\frac{\delta^2}{L_0}\,{\bf P}_0^\perp {\Psi}_{L_0}\,{\bf P}_0.
\end{split}
\end{gather*}
Since ${\bf T}_2(\delta)={\bf T}_1^{\rm T}(\delta)$ we obtain that
\begin{gather*}
{\cal L}\left(\sqrt{N}\,{\bf T}(\delta)\right)\Longrightarrow
{\cal L}\left(\frac{\delta^2}{L_0}\big({\bf P}_0^\perp
\Psi_{L_0}\,{\bf P}_0+{\bf P}_0 \Psi_{L_0}\,{\bf P}_{0}^\perp\big)\right).
\end{gather*}

In view of \eqref{eq:res_T_delta} and
under suppositions of Proposition \ref{theor:precise_proj_noise},
${\bf P}_0^{\perp}(\delta)-{\bf P}_0^{\perp}-{\bf T}(\delta)={\bf R}(\delta)$, where
\begin{gather*}
\big\|{\bf R}(\delta)\big\|\leq
C\,\frac{\|{\bf S}_0{\bf B}(\delta){\bf P}_{0}\|^2}{1-4\|{\bf B}(\delta)\|/\mu}=
C\,\frac{\|{\bf P}_\mu{\bf B}(\delta){\bf P}_{0}\|^2/(L_0K)^2}{1-4\|{\bf B}(\delta)\|/L_0K}\, .
\end{gather*}
We have already proved that $4\|{\bf B}(\delta)\|/\mu_{\min}\rightarrow 4\delta^2/L_0<1$ with probability 1.
Also, ${\bf S}_0{\bf B}(\delta){\bf P}_{0}\rightarrow 0$ with probability 1 (see demonstration of
Proposition \ref{prop:gen_wight_noise_coarse})
and in view of \eqref{eq:PmuGP_0_CLT} we obtain that
\begin{gather*}
\sqrt{K}{\bf S}_0{\bf B}(\delta){\bf P}_{0}=\sqrt{K}\,\frac{{\bf P}_0^\perp{\bf B}(\delta){\bf P}_{0}}{L_0K}
\stackrel{{\rm d}}{\longrightarrow} \delta^2\,{\bf P}_0^\perp {\Psi}_{L_0}\,{\bf P}_0/L_0.
\end{gather*}
Thus $\sqrt{N}\,\big\|{\bf S}_0{\bf B}(\delta){\bf P}_{0}\big\|\rightarrow 0$ and
$\sqrt{N}\,\big\|{\bf R}(\delta)\big\|\rightarrow 0$ a.s.
Since $L_0$ is fixed, then $\sqrt{N}{\bf R}(\delta)\rightarrow {\bf 0}$ and
\begin{gather*}
{\cal L}\!\left(\sqrt{K}\big({\bf T}(\delta)\!+\!{\bf R}(\delta)\big)\right)\!\Longrightarrow\!
{\cal L}\!\left(\frac{\delta^2}{L_0}\big({\bf P}_0^\perp \Psi_{L_0}{\bf P}_0
\!+\!{\bf P}_0 \Psi_{L_0}{\bf P}_{0}^\perp\big)\!\right)\!.
\end{gather*}
In view of the convergence $\mathbb{P}(\Omega_N)\rightarrow 1$,
conditional distributions
\begin{gather*}
{\cal L}\!\left(\sqrt{K}\big({\bf T}(\delta)+{\bf R}(\delta)\big)\,\big|\,\Omega_N\right)
\end{gather*}
 have the same weak limit as
${\cal L}\left(\sqrt{K}\big({\bf T}(\delta)+{\bf R}(\delta)\big)\right)$, and the proof is complete.

\vspace{3mm}
{\bf Proof of Proposition \ref{prop:close_R}}

\vspace{1mm}
\noindent
Note that
\begin{gather*}
\begin{split}
&R(\delta)-R= -\,{\bf G}_L \frac{1}{\|{\bf P}_0(\delta)\mathfrak{e}_L\|^2}\,
\big({\bf P}_0(\delta)-{\bf P}_0\big)\mathfrak{e}_L\\
&+
{\bf G}_L\left(\frac{1}{\|{\bf P}_0(\delta)\mathfrak{e}_L\|^2}-\frac{1}{\|{\bf P}_0\mathfrak{e}_L\|^2}\right)
{\bf P}_0\mathfrak{e}_L=J_1+J_2.
\end{split}
\end{gather*}

Since $\|{\bf G}_L\|=1$, ${\bf P}_0(\delta)-{\bf P}_0=
{\bf P}_0^\perp-{\bf P}_0^\perp(\delta)$, and
$\|{\bf P}_0(\delta)\mathfrak{e}_L\|\geq \|{\bf P}_0\mathfrak{e}_L\|-\Delta {\bf P}(\delta)$,
then
\begin{gather*}
{\big\|{\bf P}_0(\delta)\mathfrak{e}_L\big\|}\,\big/\,{\big\|{\bf P}_0\mathfrak{e}_L\big\|}\geq 1 -
{\Delta {\bf P}(\delta)}\,\big/\,{\big\|{\bf P}_0\mathfrak{e}_L\big\|}\ ,
\end{gather*}
\begin{gather*}
\|J_1\|\!\leq\! \frac{\Delta {\bf P}(\delta)}{\|{\bf P}_0\mathfrak{e}_L\|^2}\!
\left(\frac{\|{\bf P}_0\mathfrak{e}_L\|}{\|{\bf P}_0(\delta)\mathfrak{e}_L\|}\right)^2
\!\leq\!
\frac{\Delta {\bf P}(\delta)}{1\!-\!\vartheta^2}\left(\!1\!-\!
\frac{\Delta {\bf P}(\delta)}{\sqrt{1\!-\!\vartheta^2}}\right)^{-2},
\end{gather*}
and
\begin{gather*}
\begin{split}
&\|J_2\|\leq
\frac{1}{\|{\bf P}_0\mathfrak{e}_L\|^2}
\left(\frac{\|{\bf P}_0\mathfrak{e}_L\|^2}{\|{\bf P}_0(\delta)\mathfrak{e}_L\|^2}-1\right)\\
&\leq
\frac{1}{\|{\bf P}_0\mathfrak{e}_L\|^2}
\left(\frac{1}{(1-\Delta {\bf P}(\delta)/\|{\bf P}_0\mathfrak{e}_L\|)^2}-1\right)\\
&\leq
\frac{2\,\Delta {\bf P}(\delta)}{(1-\vartheta^2)^{3/2}}\,
\left(1- \frac{\Delta {\bf P}(\delta)}{\sqrt{1-\vartheta^2}}\right)^{-2}.
\end{split}
\end{gather*}
Thus \eqref{eq:rec_DeltaR} is proved.

\vspace{3mm}
{\bf Proof of Lemma \ref{lem:lin_indep}}

\vspace{1mm}
\noindent
1. Denote ${\bf U}=[U_1:\ldots:U_d]$. Since  $U_1,\ldots,U_d$ are linearly independent, then
${\bf P}{\bf U}X={\bf U}X\neq 0$ for any $X\!\neq\! 0$. Consider the vector ${\bf Q}{\bf U}X \!=\!\big({\bf Q}\!-\!{\bf P}\big){\bf U}X\!+\!{\bf U}X$. Then
\begin{gather*}
\begin{split}
&\big\|{\bf Q}{\bf U}X\big\|\geq \big\|{\bf U}X\big\|-
\big\|{\bf Q}-{\bf P}\big\|\,\big\|{\bf U}X\big\|\\
&=
\big(1-\big\|{\bf Q}-{\bf P}\big\|\big)\,\big\|{\bf U}X\big\|>0.
\end{split}
\end{gather*}
2. Consider the linear space ${\bf P}\mathbb{V}$. Evidently, ${\bf P}\mathbb{V}$ is a subspace of $\mathbb{U}$. Due to
the first assertion, the dimension of ${\bf P}\mathbb{V}$ equals $d$. Therefore, ${\bf P}\mathbb{V}=\mathbb{U}$.
This completes the proof.

\vspace{3mm}
{\bf Proof of Proposition \ref{prop:Delta_P1}}

\vspace{1mm}
\noindent
Denote for short ${\bf X}\!=\!{\bf U}^{\rm T}{\bf F}_1{\bf U}$, ${\bf Y}\!=\!{\bf U}^{\rm T}{\bf F}_2{\bf U}$,
${\bf X}(\delta)=
{\bf U}^{\rm T}{\bf P}_0^\perp(\delta){\bf F}_1{\bf P}_0^\perp(\delta){\bf U}$,
and ${\bf Y}(\delta)\!=\!{\bf U}^{\rm T}{\bf P}_0^\perp(\delta){\bf F}_2{\bf P}_0^\perp(\delta){\bf U}$.
Note that $\|{\bf Y}(\delta)\|\leq \|{\bf U}\|^2$.
Then
\begin{gather*}
\begin{split}
&\widehat{\bf D}(\delta)-{\bf D}={\bf X}^{-1}(\delta){\bf Y}(\delta)-{\bf X}^{-1}{\bf Y}\\
&=
\big({\bf X}^{-1}(\delta)-{\bf X}^{-1}\big){\bf Y}(\delta)+{\bf X}^{-1}\big({\bf Y}(\delta)-{\bf Y}\big)=J_1+J_2.
\end{split}
\end{gather*}
Since
${\bf Y}(\delta)-{\bf Y}=
{\bf U}^{\rm T}\Delta {\bf P}(\delta){\bf F}_2{\bf P}_0^\perp(\delta){\bf U}+
{\bf U}^{\rm T}{\bf P}_0^\perp{\bf F}_2\Delta {\bf P}(\delta){\bf U}$,
then $\|{\bf Y}(\delta)-{\bf Y}\|\leq 2\|{\bf U}\|^2\,\Delta {\bf P}(\delta)$ and
we obtain the inequality
\begin{gather*}
\big\|J_2\big\|\leq
\frac{2\big\|{\bf U}\big\|^2}{\|{\bf U}^{\rm T}{\bf F}_1{\bf U}\|}\ \Delta {\bf P}(\delta)
=\frac{2}{\upsilon}\ \Delta {\bf P}(\delta).
\end{gather*}
In the same manner,
$\|{\bf X}(\delta)-{\bf X}\|\leq 2\big\|{\bf U}\big\|^2\,\Delta {\bf P}(\delta)$ and
$\big\|{\bf X}^{-1}\big({\bf X}(\delta)-{\bf X}\big)\big\|\leq 2\Delta {\bf P}(\delta)/\upsilon<1.$
Therefore,
\begin{gather*}
{\bf X}^{-1}(\delta)-{\bf X}^{-1}=
\left(\sum_{j=1}^\infty\Big({\bf I}+{\bf X}^{-1}\big({\bf X}(\delta)-{\bf X}\big)\Big)^j\right){\bf X}^{-1}
\end{gather*}
and
\begin{gather*}
\begin{split}
&\big\|\big({\bf X}^{-1}(\delta)-{\bf X}^{-1}\big){\bf Y}(\delta)\big\|\\
&\leq
\frac{\|{\bf Y}(\delta)\|}{\|{\bf X}\|}\,
\frac{\big\|{\bf X}^{-1}\big({\bf X}(\delta)-{\bf X}\big)\big\|}
{1-\big\|{\bf X}^{-1}\big({\bf X}(\delta)-{\bf X}\big)\big\|}
\leq
\frac{2}{\upsilon^2}\,
\frac{\Delta {\bf P}(\delta)}{1-2\Delta {\bf P}(\delta)/\upsilon}\, .
\end{split}
\end{gather*}
The proof is complete.

\vspace{3mm}
{\bf Proof of Proposition \ref{prop:LK_even_odd}}

\vspace{1mm}
\noindent
The first assertion follows from Remark \ref{rem:recon_zero}. Let us prove the assertions 2 -- 4.
Since $\|{\bf H}\|=\|{\bf E}\|=\sqrt{LK}$ and in view of \eqref{eq:res_prec_const_saw},
\begin{gather*}
\begin{split}
&\left\|\left({\bf P}_0^{\perp}(\delta)-{\bf P}_0^{\perp}-{\bf M}(\delta)\right){\bf H}(\delta)\right\|\\
&=
\begin{cases}
O\left(\!\sqrt{L/K^3}\right)&
\begin{aligned}
&\mbox{\rm for odd} \  K, \  \mbox{\rm even}\  L,\\
&\mbox{\rm  and}  \ L=o(K^3),
\end{aligned}\\
O\left(\!\sqrt{K/L^3}\right)&
\begin{aligned}
&\mbox{\rm for odd} \  L, \  \mbox{\rm even}\  K,\\
&\mbox{\rm  and}  \ K=o(L^3),\\
\end{aligned}\\
O\left(\!\sqrt{L/K^3}+\sqrt{K/L^3}\right)&
\begin{aligned}
&\mbox{\rm for} \ K  \ {\rm and}\  L  \ \mbox{\rm both odd  and}\\
&\min(L/K^3,\,K/L^3) \rightarrow \infty.
\end{aligned}
\end{cases}
\end{split}
\end{gather*}
This means that we must study  the entries of the matrix
$\mathcal{S}\big({\bf M}(\delta){\bf H}(\delta)+\delta {\bf P}_0^{\perp}{\bf E}\big)$.
Calculations show
that under denotation $V_E=W_LE_K^{\rm T}$ and
$W_E(\delta)= E_L W_K^{\rm T}  +\delta\, V_E$,
\begin{gather*}
{\bf M}(\delta){\bf H}(\delta)+\delta{\bf P}_0^{\perp}{\bf E}\\
\!=\!
\begin{cases}
{\displaystyle
\frac{\delta}{1\!-\!\delta^2} \frac{W_E(\delta)}{K}}&
\begin{aligned}
&\mbox{for odd} \,  K \\
&\mbox{and even} \,  L,
\end{aligned}\\
{\displaystyle
\frac{\delta^2}{1\!-\!\delta^2} \frac{W_E(\delta)}{{L}}\!+\!
\delta \frac{V_E}{L}}+O(L^{-2})&
\begin{aligned}
&\mbox{for odd} \,  L\\
&\mbox{and even} \,  K,
\\
\end{aligned}\\
\displaystyle
\frac{\delta}{1\!-\!\delta^2}\!\left(\!\frac{1}{K}\!+\!\frac{\delta}{L}\!\right)\!
\left(W_E(\delta)\!+\!O(L^{-1})\right)
\!+\!\delta \frac{V_E}{L}&
\begin{aligned}
&\mbox{for} \,  K,L\\
&\mbox{both  odd}.
\end{aligned}
\end{cases}
\end{gather*}
Now the result becomes clear.

\vspace{3mm}
{\bf Proof of Proposition \ref{prop:exp_const_hankel}}

\vspace{1mm}
\noindent
Let us consider the case $L/N\rightarrow \alpha\in (0,1)$ (other cases are studied in the same manner).
It is easy to see that
\begin{gather*}
\begin{split}
&\Lambda\big(\delta {\bf V}_0^{(1)}\big)
\!=\!
\delta\big({\bf V}_0^{(1)}{\bf H}+{\bf P}_0^\perp{\bf E}\big)\!+\!\delta^2{\bf V}_0^{(1)}{\bf E}\\
&=\Lambda_1\big(\delta {\bf V}_0^{(1)}\big)+O(N^{3/2}a^{-N})
\end{split}
\end{gather*}
and
\begin{gather*}
\Lambda_1\big(\delta {\bf V}_0^{(1)}\big)=
\!\delta\!\left(-\frac{\beta_L\beta_K{\bf H}}{\|W_L\|^2\|W_K\|^2}\!+\!\frac{\beta_K E_LW_K^{\rm T}}{\|W_K\|^2}
\!+\!\frac{\beta_L  W_L E_K^{\rm T}}{\|W_L\|^2}\!\right).
\end{gather*}
Therefore,
\begin{gather*}
\big\|\widetilde{\mathrm{F}}_N(\delta)-\mathrm{F}_N\big\|_{\max}=\lambda_1\big(\delta{\bf V}_0^{(1)}\big)+
O\big(N^2a^{-N}\big),
\end{gather*}
where $\lambda_1\big(\delta{\bf V}_0^{(1)}\big)=\big\|\mathcal{S}\Lambda_1\big(\delta{\bf V}_0^{(1)}\big)\big\|_{\max}$.

Consider the series $g_{0},\ldots,g_{N-1}$ corresponding to the Hankel matrix
$\mathcal{S}\Lambda_1\big(\delta{\bf V}_0^{(1)}\big)$.
Since hankelization does not modify extreme entries $a_{00}$ and $a_{n-1\,m-1}$ of any $n\times m$
matrix ${\bf A}=\{a_{ij}\}_{i=0,\,j=0}^{n,\, m}$, then it is easy to see that
\begin{gather*}
\begin{split}
&g_{N-1}=\delta\Big(-\frac{\beta_L\beta_K}{\|W_L\|^2\|W_K\|^2}\ a^{N-1}\\
&+
\frac{\beta_K}{\|W_K\|^2}\ a^{K-1}+\frac{\beta_L}{\|W_L\|^2}\ a^{L-1}\Big)\rightarrow \delta(a^2-1)/a^2\neq 0.
\end{split}
\end{gather*}
Thus the proof is complete.

\section{Appendix 2. Trajectory matrices of stationary processes}
\label{sect:Appendix}

Let the ``noise series'' $\mathrm{E}$ be a stationary random series with $\mathbb{E}e_n=0$ and
$\mathbb{D}e_n=1$.(Then
$\delta^2$ serves as the variance of the perturbation series $\delta \mathrm{E}$.)
Denote $\mathbf{E}=\mathbf{E}_{L,K}$ the corresponding random
trajectory $L\times K$ matrix.

In this section we present an upper bound for $\|\mathbf{E}_{L,K}\|$ under supposition that $e_n$ is a {\it linear
stationary process}. More precisely, we suppose that
\begin{gather}
\label{eq:linear_process}
e_n=\sum_{j=-\infty}^\infty c_{j}\,\varepsilon_{j+n},
\end{gather}
where $\varepsilon_j$ is a sequence of independent random variables with $\mathbb{E}\varepsilon_j=0$,
$\mathbb{D}\varepsilon_j=1$ and $\sup_n\mathbb{E}|\varepsilon_{n}|^{3}<\infty$.
In addition we assume that $\sum_j|c_j|<\infty$. Since $\mathbb{D}e_n=1$, then
$\sum_jc_j^2=1$.

Note that this
model of a stationary process is widely used in statistics (for example, see \cite{Anderson71}).
In particular, if we put
\begin{gather}
\label{eq:autoregr_1}
e_n=\sum_{j=-\infty}^n a\rho^{n-j}\varepsilon_j
\end{gather}
with $|\rho|<1$ and $a=\sqrt{1-\rho^2}$, then we get the first-order autoregressional process
$e_n=\rho e_{n-1}+a\varepsilon_n$ with the covariance function $R_e(n)=\rho^n$.
Thus the important ``red noise'' example is  the particular case of \eqref{eq:linear_process}.

Section \ref{sssect:LP} is devoted to general upper bounds for norms of trajectory matrices of linear processes.
In Section \ref{ssect:WN_gen} the special case of i.i.d. random variables $e_n=\varepsilon_n$ is considered in a bit different manner.

\subsection{Linear processes}
\label{sssect:LP}
Let $\big(\Omega, {\cal F}, \mathbb{P}\big)$ be a certain probability space and consider a sequence $\varepsilon_n$
$(n=0,\pm1,\pm2,\ldots)$ of random independent variables defined on $\big(\Omega, {\cal F}, \mathbb{P}\big)$.
Assume that
$\mathbb{E}\varepsilon_n=0$, $\mathbb{D}\varepsilon_n=1$,
 and $\sup_n\mathbb{E}|\varepsilon_n|^{3}<\infty$. Denote for $n=0,1,2,\ldots$ random
 variables $e_n$ by formula \eqref{eq:linear_process} under condition $\sum_jc_j^2=1$.
 Lastly, let $N\rightarrow \infty$, $L=L(N)$, $K=N+1-L$ and denote ${\bf E}_{L,K}=
 \big\{e_{j+k-2}\big\}_{1\leq j\leq L, \, 1\leq k\leq K}$ the $L\times K$ trajectory
 matrix of the series $\mathrm{E}=(e_0,\ldots,e_n,\ldots)$.

\begin{proposition}
\label{prop:LP}
If $S\stackrel{\rm def}=\sum_j|c_j|<\infty$, then there exists $\Omega^{(0)}\in {\cal F}$ with
$\mathbb{P}\big(\Omega^{(0)}\big)=1$ such that
for any $\omega\in \Omega^{(0)}$
\begin{gather}
\label{eq:Upper_lim_sup}
\limsup_N\,\frac{\big\|{\bf E}_{L,K}(\omega)\big\|}{\sqrt{N\ln N}}\leq \gamma_0 S,
\end{gather}
where $\gamma_0>0$ is an absolute constant.
\end{proposition}
\begin{proof}
Our proof is essentially based on results and ideas published in \cite{Adamczak08}. Since we need only an upper bound for
$\big\|{\bf E}_{L,K}\big\|$, we do not use final results of \cite{Adamczak08}, dedicated to the precise rate of growth
for norms of square Hankel (and Toeplitz) matrices with i.i.d. entries. Still
the  line of our considerations is borrowed from \cite{Adamczak08}.

We start with mentioning that
\begin{gather}
 \label{eq:Hankel_norm_ineq}
\big\|{\bf E}_{L\wedge K}\big\|\leq \big\|{\bf E}_{L,K}\big\|\leq \big\|{\bf E}_{L\vee K}\big\|,
\end{gather}
where we use denotation ${\bf E}_{M}$ instead of ${\bf E}_{M,M}$. (As usual, $L\wedge K=\min(L,K)$ and
$L\vee K=\max(L,K)$.) Indeed, this follows from the remark that
${\bf E}_{L\wedge K}$ can be considered as a submatrix of
${\bf E}_{L,K}$ and that ${\bf E}_{L,K}$ is a submatrix of ${\bf E}_{L\vee K}$.

The second step is to
consider the infinite Laurent matrix
\begin{gather*}
{\bf L}_n={\bf L}_n(\mathrm{E})=
\big\{e_{|i-k|}{\bf{1}}_{|j-k|\leq n-1}\big\}_{j,k=0,\pm 1,\pm 2,\ldots}.
\end{gather*}
(Note that  ${\bf L}_n$ depends on
$e_j$ only for $0\leq j\leq n-1$.)
Then
\begin{gather}
\label{eq:norm_Laurent_exact}
\big\|{\bf L}_n\big\|=\max_{0\leq t\leq 1}\Big|e_0+\sum_{j=1}^{n-1}e_j\cos(2\pi tj)\Big|.
\end{gather}
This fact is a variant of basic relationship of Toeplitz matrices to multiplicative operators in $\mathbb{L}^2[0,1]$
(see \cite{BottcherSilbermann99}, the short outline of the proof can be found in \cite{Meckes07}).
Note that here ${\bf L}_n$ is considered as the operator on $\ell^2(\mathbb{Z})$,
where $\ell^2(\mathbb{Z})$ is treated as a set of two-side infinite complex sequences
${\bf b}=\{b_j\}_{j\in \mathbb{Z}}$ such that $||{\bf b}||^2=\sum_{j\in \mathbb{Z}}|b_j|^2<\infty$.

Laurent matrices ${\bf L}_{n}$ are useful for our goals since
\begin{gather}
\label{eq:Hankel_Laurent_ineq}
\big\|{\bf E}_n\big\|\leq \big\|{\bf L}_{2n-1}(\mathrm{E})\big\|.
\end{gather}

To prove
\eqref{eq:Hankel_Laurent_ineq},  consider the $L\times K$ matrix ${\bf R}_{L,K}=
\big\{e_{j-k+K-1}\big\}_{1\leq j\leq L, \, 1\leq k\leq K}$.
This matrix is obtained by reversing of columns of Hankel matrix ${\bf E}_{L,K}$, therefore singular values and the
corresponding left singular vectors of matrices ${\bf E}_{L,K}$ and ${\bf R}_{L,K}$ coincide. In particular,
$\|{\bf R}_{L,K}\|=\|{\bf E}_{L,K}\|$.

If we take $L=K=n$, then we see that ${\bf R}_{n,n}$  is a submatrix of the Laurent
matrix $L_{2n-1}$. Indeed,
the matrix $R_{n,n}$ is placed on the intersection of the set
$I=\big\{-(n-1),-(n-2),\ldots,-1,0\big\}$ of rows and the same set $J=I$ of columns of the Laurent matrix
$L_{2n-1}$.
Thus $\big\|{\bf E}_n\big\|=\big\|{\bf R}_{n,n}\big\|\leq \big\|{\bf L}_{2n-1}\big\|$.
Note that this procedure was discussed in \cite{BrycDJ06}.

Suppose now that $\mathrm{E}=(\varepsilon_0,\varepsilon_1,\ldots,\varepsilon_n,\ldots)$ is the sequence of independent
random variables with zero means and unit variances. Then it follows from \cite[theor. 4]{Adamczak08} that
there exists  an absolute constant $\gamma_1>0$ such that
\begin{gather}
\label{eq:exp_Laurent_norm_ineq}
\mathbb{E}\big\|{\bf L}_n(\mathrm{E})\big\|\leq \gamma_1 \sqrt{n \ln n}.
\end{gather}
If we assume additionally that $\sup_n\mathbb{E}|\varepsilon_n|^{3}<\infty$, then
\begin{gather}
\label{eq:Laurent_norm_ineq}
\limsup_n\frac{\|{\bf L}_n(\mathrm{E})\|}{\sqrt{n\ln n}}\leq
\limsup_n\frac{\mathbb{E}\|{\bf L}_n(\mathrm{E})\|}{\sqrt{n\ln n}}
\end{gather}
with probability 1. This fact is proved in the demonstration of \cite[theor. 5]{Adamczak08} under more weak moment
conditions on $\varepsilon_n$. We have simplified these conditions to make them more transparent.

Collecting inequalities \eqref{eq:Hankel_norm_ineq}, \eqref{eq:Hankel_Laurent_ineq}, \eqref{eq:exp_Laurent_norm_ineq},
and \eqref{eq:Laurent_norm_ineq} we come to inequality \eqref{eq:Upper_lim_sup} for the particular case
$e_n=\varepsilon_n$. Indeed, let $n\stackrel{\rm def}=\max(L,K$). Then $\lceil N/2 \rceil\leq n\leq N$ and
\begin{gather}
\begin{split}
&\limsup_N \frac{\|{\bf E}_{L,K}\|}{\sqrt{N\ln N}}\leq \limsup_N \frac{\|{\bf E}_{n}\|}{\sqrt{N\ln N}}
\leq \limsup_N \frac{\|{\bf L}_{2n-1}\|}{\sqrt{N\ln N}}
\\
&\leq\limsup_N \frac{\|{\bf L}_{2n-1}\|}{\sqrt{(2n-1)\ln (2n-1)}}
\frac{\sqrt{(2n\!-\!1)\ln (2n\!-\!1)}}{\sqrt{N\ln N}}\\
&\leq
\limsup_{2n-1} \frac{\sqrt{2}\|{\bf L}_{2n-1}\|}{\sqrt{(2n-1)\ln (2n-1)}}\leq \sqrt{2}\gamma_1=\gamma_0,
\label{eq:main_ineq WN}
\end{split}
\end{gather}
where $\gamma_1$ is described in \eqref{eq:exp_Laurent_norm_ineq} and the last inequality holds with probability 1.

The last step is to return to the linear  process under consideration. Since
${\bf E}_{L,K}=\sum_l c_m {\bf E}_{L,K}^{(m)}$, where ${\bf E}_{L,K}^{(m)}= \big\{e_{m+j+k-2}\big\}_{1\leq j\leq L, \, 1\leq k\leq K}$, then
\begin{gather*}
\big\|{\bf E}_{L,K}\big\|\leq \sum_m |c_{m}|\, \big\|{\bf E}_{L,K}^{(m)}\big\|.
\end{gather*}

Of course, inequality
\eqref{eq:main_ineq WN} holds for operators
${\bf E}_{L,K}^{(m)}$. Taking into account the fact that $\gamma_0$ is an absolute
constant, we  come to \eqref{eq:Upper_lim_sup}.
\end{proof}
%

\begin{remark}
\label{rem:rand_hankel_norm}
1.
It follows from Proposition \ref{prop:LP} that for $\omega\in \Omega^{(0)}$
\begin{gather*}
\max_N \frac{\|{\bf E}_{L,K}(\omega)\|}{\sqrt{N\ln N}}\leq c(\omega)
\end{gather*}
with $\gamma_0 S\leq c(\omega)<\infty$. On the other hand, for any $\alpha'>\alpha_0$ there exists $N_0=
N_0(\omega,\alpha')$ such that
\begin{gather*}
\max_{N>N_0} \frac{\|{\bf E}_{L,K}(\omega)\|}{\sqrt{N\ln N}}\leq \alpha'\,S.
\end{gather*}
2. Of course, $S=1$ for the white-noise variant of $e_n$.
It is easy to check that $S=\sqrt{(1+|\rho|)/(1-|\rho|)}$
for the autoregression process \eqref{eq:autoregr_1} with $a=\sqrt{1-\rho^2}$. This means that the upper bound in the
right-hand side of \eqref{eq:Upper_lim_sup} tends to infinity as $|\rho|\rightarrow 1$.\\
3. Inequality \eqref{eq:Upper_lim_sup} holds for the arbitrary behavior of $L=L(N)$. Yet in some cases
this inequality can be refined. For example, if $n=\max(L,K)/N\rightarrow \beta<1$
as $N\rightarrow \infty$, then the right-hand side of \eqref{eq:Upper_lim_sup} can be multiplied by $\sqrt{\beta}$.
Since $\beta\geq 1/2$, the advantage is not very big.\\
4. The case when $L=L_0={\rm const}$ is much more simple.
Of course,
$\big\|{\bf E}_{L_0,K}\big\|^2$ is equal to largest eigenvalue of the $L_0\times L_0$ matrix ${\bf E}_{L_0,K}{\bf E}_{L_0,K}^{\rm T}$.

Entries of the matrix ${\bf E}_{L_0,K}{\bf E}_{L_0,K}^{\rm T}/K$
have the form
\begin{gather}
\label{eq:phi_ij(e)}
\phi_{ij}^{(K)}=\frac{1}{K}\sum_{m=0}^{K-1}e_{i+m}e_{j+m}, \quad i,j=0,\ldots,L_0-1.
\end{gather}
If the series $\mathrm{E}$ is strongly stationary and  ergodic (see \cite[ch. 5]{Shiryaev} for precise definitions),
then the series $\mathrm{E}_{ij}^{(2)}=(e_ie_j,e_{i+1}e_{j+1},\ldots,e_{i+m}e_{j+m},\ldots)$ is also strongly stationary and
ergodic and therefore $\phi_{ij}^{(K)}$ converges almost surely to $R_e(i-j)=\mathbb{E}e_ie_j$ as $K\rightarrow \infty$.

This means that $\|{\bf E}_{L_0,K}{\bf E}_{L_0,K}^{\rm T}\|/N$ tends a.s. to the maximal eigenvalue $\lambda_{\max}$
of the covariance matrix
\begin{gather*}
{\bf R}_{L_0}=\big\{R_e(i-j)\big\}_{0\leq i,j\leq L_0-1}.
\end{gather*}
Therefore,
$\|{\bf E}_{L_0, K}\|\sqrt{N}\rightarrow \sqrt{\lambda_{\max}}$ with probability 1. In particular, if $e_n$ are i.i.d. random
variables with $\mathbb{E}e_n=0$ and $\mathbb{D}e_n=1$, then ${\bf E}_{L_0,K}{\bf E}_{L_0,K}^{\rm T}/N$ converge a.s.
to the identity matrix ${\bf I}_{L_0}$ and $\|{\bf E}_{L_0,K}\|\sqrt{N}\rightarrow 1$ with probability 1.\\
5. These examples show that normalizing factor $1/\sqrt{N \ln N}$
in the lefthand side of \eqref{eq:Upper_lim_sup} is not uniformly precise. On the other hand, it is proved in
\cite[corr. 3 and sect. 3.4]{Adamczak08} that
in the case of i.i.d. random variables with zero mean and unit variances,
\begin{gather*}
\limsup_n\frac{\|{\bf E}_{n}\|}{\sqrt{n\ln n}}<\infty
\end{gather*}
with probability 1. Applying inequality \eqref{eq:Hankel_norm_ineq}, we see that
  the same result is valid for the sequence of rectangular matrices ${\bf E}_{L,K}$ under supposition that
  $\min(L,K)/N\rightarrow \alpha>0$.
\end{remark}

\subsection{The case of i.i.d. random variables}
\label{ssect:WN_gen}
Consider a probability space $(\Omega, {\cal F}, \mathbb{P})$ and a sequence of i.i.d. random variables
$\varepsilon_0,\varepsilon_1,\ldots$ defined on this space. Assume that
$\mathbb{E}\varepsilon_i\!=\!0$, $\mathbb{D}\varepsilon_i\!=\!1$,
and $\sup_i\mathbb{E}\varepsilon_i^3\!<\!\infty$.

Let us present several auxiliary assertions related to the series $\mathrm{E}$. For convenience we collect these
assertions in the following propositions.

Denote
\begin{gather*}
x_n=\sum_{k=1}^p\gamma_k \cos(2\pi\omega_k n +\phi_k)
\end{gather*}
with $\gamma_k\neq 0$ and $\omega_j\neq \omega_k$ for $j\neq k$.

\begin{proposition}
\label{prop:LIL_osc_wn}
For any $m,l\geq 0$
\begin{gather}
\label{eq:LIL_osc_wn}
\limsup_n \frac{\big|\sum_{j=0}^{n-1}x_{m+j}\varepsilon_{l+j}\big|}{\sqrt{n\ln \ln n}}=\sqrt{\sum_{k=1}^p\gamma_l^2}\,.
\end{gather}
with probability 1.
\end{proposition}
\begin{proof}
For fixed $m$ and $l$
set $\xi_j=x_{m+j}\varepsilon_{l+j}$. It can be easily checked that
\begin{gather*}
B_n^2 \stackrel{\rm def}=\mathbb{D}(\xi_0+\ldots+\xi_{n-1})=\sum_{j=0}^{n-1}x_{m+j}^2=
n\sum_{k=1}^p \gamma_k^2/2 +O(1).
\end{gather*}
Since $\sup_j \mathbb{E}|\xi_j|^3<\infty$, then it follows from \cite[ch. 17]{Borovkov98} that almost surely
\begin{gather*}
\limsup_n \frac{\big|\sum_{j=0}^{n-1}x_{m+j}\varepsilon_{l+j}\big|}{B_n\sqrt{2 n\ln \ln n}}=1.
\end{gather*}
Thus \eqref{eq:LIL_osc_wn} is proved.
\end{proof}

Denote $\phi_{ij}^{(n)}$ in the manner of  \eqref{eq:phi_ij(e)} changing $e_n$ in this formula for $\varepsilon_n$.
\begin{proposition}
\label{prop:white_noise}
Suppose that $\mathbb{E}\varepsilon_j^4<\infty.$\\
1. If $i=j$ and $n\rightarrow \infty$, then $\mathbb{E}\phi_{ii}^{(n)}=1$ and
$n\,\mathbb{D}\phi_{ii}^{(n)}=\Delta^2\stackrel{\rm def}=\mathbb{E}\varepsilon^4-1$.
Besides, $\phi_{ii}^{(n)}\rightarrow 1$ almost surely.
 If $i\neq j$ and $n\rightarrow \infty$, then $\mathbb{E}\phi_{ij}^{(n)}=0$ and
$n\,\mathbb{D}\phi_{ij}^{(n)}\rightarrow 1$.
Besides, $\phi_{ij}^{(n)}\rightarrow 0$ almost surely.\\
2. Denote $\rho^{(n)}(i,j;s,t)$ the correlation coefficient between $\phi_{ij}^{(n)}$ and $\phi_{st}^{(n)}$.
If $|j-i|=|t-s|$, then $\rho^{(n)}(i,j;s,t)\rightarrow 1$ as $n\rightarrow \infty$. If
$|j-i|\neq |t-s|$, then $\rho^{(n)}(i,j;s,t)=0$.\\
3. Let $M$ be fixed and $n\rightarrow\infty$. Denote $\Phi_n=\big\{\phi_{ij}^{(n)}\big\}_{i,j=0}^{M-1}$.
Then $\Phi_n\rightarrow {\bf I}$ with probability 1 and
 the limit distribution of the sequence of random matrices $\sqrt{n}(\Phi_n-{\bf I})$
coincides with the distribution of the random $M\times M$ matrix $\Psi_M=\big\{\psi_{ij}\big\}_{i,j=1}^{M}$ such that\\
a) $\psi_{ij}=\psi_{st}$ for $|i-j|=|s-t|$;\\ b) random variables $\psi_{11}, \psi_{12},\ldots,  \psi_{1M}$
are independent;\\
 c) $\psi_{11}\in \mathrm{N}(0,\Delta^2)$ and $\psi_{1i}\in \mathrm{N}(0,1)$ for $i>1$.\\
4. If $i=j$, then
\begin{gather*}
\limsup_n \sqrt n\, \frac{\big|\phi_{ii}^{(n)}-1\big|}{\sqrt{\ln\ln n}}=\sqrt{\mathbb{E}\varepsilon^4}
\end{gather*}
with probability 1. If $i\neq j$, then almost surely
\begin{gather*}
\limsup_n \sqrt n\, \frac{\big|\phi_{ij}^{(n)}\big|}{\sqrt{\ln\ln n}}=1.
\end{gather*}
\end{proposition}
\begin{proof}
1. The  assertion follows from the fact that random variables $\varepsilon_0^2,\varepsilon_1^2,\ldots$ are i.i.d. with
$\mathbb{E}\varepsilon^2=1$ and $\mathbb{D}\varepsilon^2=\Delta^2$. The a.s. convergences of $\phi_{ij}^{(n)}$ are already
discussed in Remark \ref{rem:rand_hankel_norm}.
Other assertions are checked straightforwardly.\\
2. Let $i\leq j$ and $i\leq s\leq t$. (Other variants are treated in the same manner.)
If $j-i=t-s=h$, then $\mathbb{E}\phi_{ij}^{(n)}=\mathbb{E}\phi_{st}^{(n)}$,
$\mathbb{D}\phi_{ij}^{(n)}=\mathbb{D}\phi_{st}^{(n)}$, and  under denotation
\begin{gather*}
A_n= \frac{1}{n}\sum_{m=0}^{s-i-1}\varepsilon_{i+m}\varepsilon_{i+h+m},
\end{gather*}
for  $n$ big enough
\begin{gather*}
\begin{split}
&\phi_{ij}^{(n)}=\frac{1}{n}\,\sum_{m=0}^{n-1}\varepsilon_{i+m}\varepsilon_{i+h+m}
=
\frac{1}{n}\sum_{m=s-i}^{n-1}\varepsilon_{i+m}\varepsilon_{i+h+m}+A_n
\\
&
=
\frac{1}{n}\sum_{p=0}^{n-1-s+i}\varepsilon_{s+p}\varepsilon_{t+p}+A_n\\
&=\frac{1}{n}\sum_{p=0}^{n-1}\varepsilon_{s+p}\varepsilon_{t+p}\!+\!A_n
\!-\!\frac{1}{n}\sum_{p=n-s+i}^{n-1}\varepsilon_{s+p}\varepsilon_{t+p}
\!=\!\phi_{st}^{(n)}\!+\!B_n
\end{split}
\end{gather*}
with $\mathbb{E}B^2_n=O(1/n^2)$.
Therefore, $\mathbb{E}\phi_{ij}^{(n)}\phi_{st}^{(n)}= \mathbb{E}\big(\phi_{ij}^{(n)}\big)^2
+\mathbb{E}\big(\phi_{ij}^{(n)}B_n\big)$,
\begin{gather*}
\big|\mathbb{E}\big(\phi_{ij}^{(n)}B_n\big)\big|
\leq \sqrt{\mathbb{E}\big(\phi_{ij}^{(n)}\big)^2\mathbb{E}B^2_n}\rightarrow 0,
\end{gather*}
and $\rho^{(n)}(i,j;s,t)\!\rightarrow\! 1$.
If $j\!-\!i\!\neq\! t\!-\!s$, then $\mathbb{E}\phi_{ij}^{(n)}\mathbb{E}\phi_{st}^{(n)}\!=\!0$,
 $\mathbb{E}\varepsilon_{i+m}\varepsilon_{j+m}\varepsilon_{s+p}\varepsilon_{t+p}=0$ for any nonnegative $m,p$, and
$\rho^{(n)}(i,j;s,t)=0$.\\
3. The convergence $\Phi_n\rightarrow {\bf I}$ immediately follows from the first assertion of the proposition.

If $j=|t-s|$, then
\begin{gather*}
\begin{split}
\sqrt{n}\big(\phi_{ts}^{(n)}-\mathbb{E}\phi_{ts}^{(n)}\big)=
\sqrt{n}\big(\phi_{st}^{(n)}-\mathbb{E}\phi_{st}^{(n)}\big)
=
\sqrt{n}\big(\phi_{0j}^{(n)}-\mathbb{E}\phi_{0j}^{(n)}\big)+C_n,
\end{split}
\end{gather*}
where $C_n\stackrel{\mathbb{P}}\rightarrow 0$ as $n\rightarrow \infty$. This means that all we need is to prove
that  coordinates of the random vector
\begin{gather}
\label{eq:CLT_vector_phi}
\sqrt{n}\Big(\phi_{00}^{(n)}-\mathbb{E}\phi_{00}^{(n)},\ldots,\phi_{0,M-1}^{(n)}-\mathbb{E}\phi_{0,M-1}^{(n)}\Big)
\end{gather}
are asymptotically  independent and asymptotically normal  with  proper variances.

Denote for $0\leq j\leq M-1$
\begin{gather*}
\tau_j^{(n)}=\frac{1}{n-j}\sum_{m=0}^{n-j-1}\varepsilon_{m}\varepsilon_{m+j}.
\end{gather*}
Then $\mathbb{E}\tau_j^{(n)}=\mathbb{E}{\phi_{0,j}^{(n)}}$ and
the difference  $\sqrt{n}\big(\phi_{0,j}^{(n)}-\tau_j^{(n)}\big)$ tends to zero in probability as $n\rightarrow \infty$.
This means that  the asymptotic distribution of the random vector
\begin{gather}
\label{eq:CLT_covar}
\sqrt{n}\Big(\tau_{0}^{(n)}-\mathbb{E}\tau_{0}^{(n)},\ldots,\tau_{M-1}^{(n)}-\mathbb{E}\tau_{M-1}^{(n)}\Big)
\end{gather}
(provided that this distribution exists)
coincides with the asymptotic distribution of the vector \eqref{eq:CLT_vector_phi}.

Note that $\tau_{j}^{(n)}$ is the usual unbiased estimate of the covariance function of the stationary
series $\mathrm{E}$. Due to  \cite[th. 8.4.2]{Anderson71}, distribution of the vector \eqref{eq:CLT_covar}
weakly tends to $M$-dimensional Gaussian distribution with zero mean and covariance matrix ${\bf W}=
\big\{w_{kl}\big\}_{k,l=0}^{M-1}$, where
\begin{gather*}
w_{kl}=\frac{1}{\pi} \int_{-\pi}^{\pi}\cos(k\lambda)\,\cos(l\lambda)\,d\lambda+
(\mathbb{E}\varepsilon^4-3)\delta_0(k)\delta_0(l)
\end{gather*}
and $\delta_0(j)$ is either 1 or 0 depending on the conditions $j=0$ and $j\neq 0$. Since
\begin{gather*}
\int_{-\pi}^{\pi}\cos(k\lambda)\,\cos(l\lambda)\,d\lambda=
\begin{cases}
0& \mbox{for} \ \ k\neq l,\\
\pi& \mbox{for} \ \ k=l\neq 0,\\
2\pi& \mbox{for}\ \ k=l=0,
\end{cases}
\end{gather*}
we see that $w_{00}=\mathbb{E}\varepsilon^4-1$, $w_{kk}=1$ for $k\neq 0$, and $w_{kl}=0$ for $k\neq l$.\\
4. The first equality is the particular case of the standard Law of Iterated Logarithms (LIL) for i.i.d. random variables.
To demonstrate the second equality it is sufficient to put $i=1$. Then under denotation $N=\lfloor n/j \rfloor$
\begin{gather*}
n\phi_{1j}^{(n)}=\sum_{k=1}^{j}\sum_{l=0}^{N-1}\varepsilon_{k+jl}\varepsilon_{k+j(l+1)-1}+r_n
=\sum_{k=1}^{j}S_{k}^{(j)}+r_n.
\end{gather*}
Since $S_{k}^{(j)}$ is the sum of i.i.d. variables we can apply the mentioned LIL to this sum and asymptotically change
$N$ for $n/j$. As $\mathbb{D}r_n$ is bounded above, this completes the proof.
\end{proof}


\end{document}